\documentclass[11pt]{article}

\usepackage{CJKutf8}
\usepackage{empheq}
\usepackage{amsmath}
\usepackage{amsthm}
\usepackage{amssymb}
\usepackage{mathtools}
\usepackage[cal=euler,scr=boondoxupr]{mathalpha}
\usepackage{upgreek}
\usepackage[T1]{fontenc}
\usepackage{lmodern}
\usepackage{mleftright}
\usepackage[sf,bf]{titlesec}
\usepackage{bm}
\usepackage{euscript}
\usepackage{url}            
\usepackage{booktabs}       
\usepackage{amsfonts}       
\usepackage{nicefrac}       
\usepackage{microtype}      
\usepackage{color}
\usepackage{tabu}
\usepackage{multirow}
\usepackage{float}
\usepackage{booktabs}
\usepackage{algorithm}
\usepackage{algorithmic}
\usepackage{color}
\usepackage{tabu}
\usepackage{comment}
\usepackage{graphicx}
\usepackage{subcaption}
\usepackage{framed}
\usepackage{enumitem}
\usepackage{textcomp}
\usepackage{setspace}
\usepackage{hyperref}
\usepackage[capitalize,noabbrev]{cleveref}
\usepackage{latexsym}
\usepackage{tcolorbox}
\usepackage[numbers,sort&compress]{natbib}
\usepackage{alltt}
\usepackage{stmaryrd}
\usepackage{xspace}
\usepackage{euscript}
\usepackage{xfrac}    
\usepackage{nicefrac} 
\usepackage{tikz}
\usepackage{tikz-qtree}
\usepackage{todonotes}
\usepackage{inconsolata} 

\usepackage[margin=1in]{geometry}
\usepackage{fancyhdr}
\fancypagestyle{plain}{\fancyhf{} 
	
	\fancyfoot[C]{\textsf{\thepage}}
}
\definecolor{mydarkblue}{rgb}{0,0.08,0.45}
\hypersetup{ %
	colorlinks=true,
	linkcolor=mydarkblue,
	citecolor=mydarkblue,
	filecolor=mydarkblue,
	urlcolor=mydarkblue}

\newcommand{\calA}{\mathcal{A}}

\newcommand{\calC}{\mathcal{C}}
\newcommand{\calD}{\mathcal{D}}
\newcommand{\calE}{\mathcal{E}}

\newcommand{\calM}{\mathcal{M}}

\newcommand{\scrO}{\mathscr{O}}
\newcommand{\scrP}{\mathscr{P}}

\newcommand{\euE}{\EuScript{E}}

\newcommand{\Ex}{\mathbb{E}}

\newcommand{\RR}{\mathbb{R}}

\newcommand{\OO}{\mathbb{O}}
\newcommand{\Rp}{\RR_+}
\newcommand{\Rpp}{\RR_{++}}
\newcommand{\bbS}{\mathbb{S}}

\newcommand{\NN}{\mathbb{N}}
\newcommand{\VV}{\mathbb{V}}

\DeclareMathOperator{\proj}{proj}

\DeclareMathOperator*{\argmin}{argmin}

\DeclareMathOperator*{\minimize}{minimize}

\DeclareMathOperator*{\argmax}{argmax}

\DeclareMathOperator*{\Diag}{Diag}

\newcommand{\diag}{\mathrm{diag}}
\newcommand{\sgn}{\mathrm{sgn}}
\DeclareMathOperator*{\dom}{dom}

\newcommand{\sumK}{\sum_{k=1}^K}

\newcommand{\sumd}{\sum_{i=1}^d}

\newcommand{\tr}{\mathrm{tr}}

\renewcommand{\vec}{\mathrm{vec}}

\newcommand{\sfP}{\mathsf{P}}

\newcommand{\sff}{\mathsf{f}}

\newcommand{\sfm}{\mathsf{m}}

\newcommand{\sfv}{\mathsf{v}}

\newcommand{\oRR}{\overline{\RR}}

\newcommand{\midd}{\,|\kern-0.25ex|\,}
\newcommand{\setn}{\llbracket n\rrbracket}

\newcommand{\setK}{\llbracket K\rrbracket}
\newcommand{\setT}{\llbracket T\rrbracket}

\newcommand{\dotp}[2]{\langle #1, #2\rangle}
\newcommand{\dotpF}[2]{\left\langle #1, #2\right\rangle_{\rm F}}

\newcommand{\sfd}{\mathsf{d}}

\newcommand{\half}{\sfrac12}

\DeclareMathAlphabet\rsfscr{U}{rsfso}{m}{n}

\let\le\leqslant
\let\ge\geqslant

\let\hat\widehat
\let\tilde\widetilde
\let\bar\overline

\makeatletter
\DeclareFontFamily{OMX}{MnSymbolE}{}
\DeclareSymbolFont{MnLargeSymbols}{OMX}{MnSymbolE}{m}{n}
\SetSymbolFont{MnLargeSymbols}{bold}{OMX}{MnSymbolE}{b}{n}
\DeclareFontShape{OMX}{MnSymbolE}{m}{n}{
	<-6>  MnSymbolE5
	<6-7>  MnSymbolE6
	<7-8>  MnSymbolE7
	<8-9>  MnSymbolE8
	<9-10> MnSymbolE9
	<10-12> MnSymbolE10
	<12->   MnSymbolE12
}{}
\DeclareFontShape{OMX}{MnSymbolE}{b}{n}{
	<-6>  MnSymbolE-Bold5
	<6-7>  MnSymbolE-Bold6
	<7-8>  MnSymbolE-Bold7
	<8-9>  MnSymbolE-Bold8
	<9-10> MnSymbolE-Bold9
	<10-12> MnSymbolE-Bold10
	<12->   MnSymbolE-Bold12
}{}

\let\llangle\@undefined
\let\rrangle\@undefined
\DeclareMathDelimiter{\llangle}{\mathopen}%
{MnLargeSymbols}{'164}{MnLargeSymbols}{'164}
\DeclareMathDelimiter{\rrangle}{\mathclose}%
{MnLargeSymbols}{'171}{MnLargeSymbols}{'171}
\makeatother

\theoremstyle{plain}
\newtheorem{theorem}{Theorem}[section]
\newtheorem{proposition}[theorem]{Proposition}

\theoremstyle{definition}
\newtheorem{definition}{Definition}[section]
\newtheorem{assumption}{Assumption}[section]
\newtheorem{example}{Example}[section]

\theoremstyle{remark}
\newtheorem{remark}{Remark}[section]

\crefname{assumption}{Assumption}{Assumptions}
\Crefname{assumption}{Assumption}{Assumptions}
\crefname{problem}{Problem}{Problems}
\Crefname{problem}{Problem}{Problems}
\crefname{example}{Example}{Examples}
\Crefname{example}{Example}{Examples}

\let\le\leqslant
\let\ge\geqslant

\let\hat\widehat
\let\tilde\widetilde
\let\bar\overline

\newcommand{\lrangle}[1]{\left\llangle #1 \right\rrangle}
\renewcommand{\dotpF}[2]{\lrangle{#1, #2}_{\rm F}}
\newcommand{\norm}[1]{\left\lVert#1\right\rVert}
\newcommand{\euclidnorm}[1]{\left\lVert#1\right\rVert_2}

\newcommand{\vecnorm}[2]{\left\lVert #1 \right\rVert_{{#2}}}
\newcommand{\matsnorm}[2]{\lvert\kern-0.25ex\lvert\kern-0.25ex\lvert #1 \rvert\kern-0.25ex\rvert\kern-0.25ex\rvert_{#2}}

\newcommand{\fronorm}[1]{\matsnorm{#1}{\mathrm{F}}}
\newcommand{\infnorm}[1]{\vecnorm{#1}{\mbox{\tiny{$\infty$}}}}

\newcommand{\onenorm}[1]{\vecnorm{#1}{1}}

\newcommand{\nucnorm}[1]{\matsnorm{#1}{\mathrm{nuc}}}

\newcommand{\specnorm}[1]{\matsnorm{#1}{\mathrm{S}}}

\renewcommand{\left}{\mleft}
\renewcommand{\right}{\mright}

\newcommand{\PL}{P\L{}\xspace}

\newcommand{\rank}{\mathrm{rank}}

\newcommand{\Adam}{\textsc{Adam}\xspace}
\newcommand{\AdamW}{\textsc{AdamW}\xspace}
\newcommand{\SGD}{\textsc{SGD}\xspace}
\newcommand{\SGDM}{\textsc{SGDM}\xspace}
\newcommand{\RMSprop}{\textsc{RMSprop}\xspace}

\newcommand{\Unif}{\mathsf{Unif}}

\newcommand{\AdaGrad}{\textsc{AdaGrad}\xspace}

\newcommand{\Adadelta}{\textsc{Adadelta}\xspace}
\newcommand{\Adafactor}{\textsc{Adafactor}\xspace}

\newcommand{\Muon}{\textsc{Muon}\xspace}
\newcommand{\Shampoo}{\textsc{Shampoo}\xspace}
\newcommand{\Lion}{\textsc{Lion}\xspace}
\newcommand{\Sophia}{\textsc{Sophia}\xspace}

\newcommand{\signSGD}{\textsc{signSGD}\xspace}
\newcommand{\Signum}{\textsc{Signum}\xspace}
\newcommand{\RProp}{\textsc{RProp}\xspace}

\newcommand{\AdaBelief}{\textsc{AdaBelief}\xspace}

\newcommand{\NAdamW}{\textsc{NAdamW}\xspace}

\newcommand{\Scion}{\textsc{Scion}\xspace}

\newcommand{\PolarGrad}{\textsc{PolarGrad}\xspace}

\newcommand{\PolarHB}{\textsc{PolarHB}\xspace}
\newcommand{\PolarGradM}{\textsc{PolarGradM}\xspace}

\newcommand{\PolarSGDM}{\textsc{PolarSGDM}\xspace}
\newcommand{\PolarSGD}{\textsc{PolarSGD}\xspace}

\newcommand{\PolarMuon}{\textsc{PolarMuon}\xspace}
\newcommand{\AltGD}{\textsc{AltGD}\xspace}
\newcommand{\Gluon}{\textsc{Gluon}\xspace}
\newcommand{\MuonAll}{\textsc{MuonAll}\xspace}

\newcommand{\AlgoPerf}{\textsc{AlgoPerf}\xspace}

\newcommand{\msgn}{\mathrm{msgn}}
\DeclareMathOperator{\SVD}{SVD}
\newcommand{\polar}{\mathrm{polar}}
\newcommand{\sfp}{\mathsf{p}}
\newcommand{\mat}{\mathrm{mat}}

\newcommand{\fstar}{f^\star}
\newcommand{\hG}{\widehat{G}}
\newcommand{\hU}{\widehat{U}}
\newcommand{\hH}{\widehat{H}}
\newcommand{\hnu}{\widehat\nu}
\newcommand{\tU}{\widetilde{U}}
\newcommand{\tH}{\widetilde{H}}
\newcommand{\tnu}{\widetilde\nu}
\newcommand{\tC}{\widetilde{C}}
\newcommand{\polarhat}{\widehat{\mathrm{polar}}}

\renewcommand{\setK}{\{1,\ldots, K\}}
\renewcommand{\setT}{\{1,\ldots, T\}}
\renewcommand{\setn}{\{1,\ldots, n\}}

\begin{document}
    \title{\PolarGrad: A Class of Matrix-Gradient Optimizers\\ from a Unifying Preconditioning Perspective}
    \author{
        Tim Tsz-Kit Lau%
        \thanks{University of Pennsylvania, Philadelphia, PA 19104, USA. Emails: 
        \href{mailto:timlautk@gmail.com}{\texttt{timlautk@gmail.com}}, 
        \href{mailto:qlong@upenn.edu}{\texttt{qlong@upenn.edu}}, 
        \href{mailto:suw@wharton.upenn.edu}{\texttt{suw@wharton.upenn.edu}}. 
        }
        \and
        Qi Long%
        \footnotemark[1]
        \and 
        Weijie Su%
        \footnotemark[1]
    }
    
    \maketitle

    \begin{abstract}
            The ever-growing scale of deep learning models and training data underscores the critical importance of efficient optimization methods. While preconditioned gradient methods such as \Adam and \AdamW are the de facto optimizers for training neural networks and large language models, structure-aware preconditioned optimizers like \Shampoo and \Muon, which utilize the matrix structure of gradients, have demonstrated promising evidence of faster convergence. In this paper, we introduce a unifying framework for analyzing ``matrix-aware'' preconditioned methods, which not only sheds light on the effectiveness of \Muon and related optimizers but also leads to a class of new structure-aware preconditioned methods. A key contribution of this framework is its precise distinction between preconditioning strategies that treat neural network weights as vectors (addressing curvature anisotropy) versus those that consider their matrix structure (addressing gradient anisotropy). This perspective provides new insights into several empirical phenomena in language model pre-training, including \Adam's training instabilities, \Muon's accelerated convergence, and the necessity of learning rate warmup for \Adam. Building upon this framework, we introduce \PolarGrad, a new class of preconditioned optimization methods based on the polar decomposition of matrix-valued gradients. As a special instance, \PolarGrad includes \Muon with updates scaled by the nuclear norm of the gradients. We provide numerical implementations of these methods, leveraging efficient numerical polar decomposition algorithms for enhanced convergence. Our extensive evaluations across diverse matrix optimization problems and language model pre-training tasks demonstrate that \PolarGrad outperforms both \Adam and \Muon. 
        \end{abstract}
        
        \section{Introduction} 
        Gradient-based optimization methods are the cornerstone for the success of modern large-scale machine learning and deep learning \citep{bottou2018optimization}. However, training very large deep neural networks remains a highly intricate task, often attributed to nonconvexity and nonsmoothness of the loss landscape of complex network architectures, as well as nonstationary data distribution. Motivated and guided by the neural scaling law \citep{kaplan2020scaling,hoffmann2022training}, we are able to achieve better model performance by scaling both model and data sizes given a certain level of compute. As the size of models scales, gigantic computational costs have been incurred. Consequently, more efficient model training algorithms have been sought relentlessly in recent years by the deep learning community. Despite more than a decade of effort, \Adam \citep{kingma2015}---the test of time award winner from the International Conference on Learning Representations (ICLR) 2025---and its decoupled weight decay variant \AdamW \citep{loshchilov2019decoupled} are still predominantly the default optimizers for training neural networks. 
        
        When designing optimizers for deep learning, a mostly overlooked fact is that neural networks are often composed of parameters of different algebraic structures---scalars in normalization layers, (bias) vectors in fully connected layers, matrices in fully connected and attention layers, and tensors in convolution layers. In traditional optimization problems, optimization variables usually have only one of the above structures (otherwise block coordinate methods are usually used; see e.g., \citep{zeng2019global,lau2018proximal} for deep learning), and they necessitate different optimization methods to solve, leading to a wide range of vector, matrix and tensor optimization methods. However, when training neural networks, elementwise optimizers such as \SGD \citep{robbins1951}, \SGDM \citep{sutskever2013importance}, \AdaGrad \citep{duchi2011adagrad,mcmahan2010adaptive}, and \Adam \citep{kingma2015} are often employed, which is equivalent to flattening and concatenating all parameters into a single vector. This treatment implicitly ignores the underlying algebraic structures of the higher-order parameters and also forgoes the optimization methods developed specifically for matrix and tensor parameters. Previous works have also pursued the direction of developing deep learning optimizers that respect the algebraic structures of different network parameters, with \Shampoo \citep{gupta2018shampoo,anil2020scalable} being the most notable example. More recently, in \citep{bernstein2024modular,bernstein2024old,large2024scalable} the use of proper norms is suggested for the design of optimizers for deep learning. This has led the introduction of \Muon \citep{jordan2024muon,bernstein2025deriving}, which has recently emerged as an empirically competitive optimizer to train transformers for both image classification and language generation, with its scalability justified in \citep{liu2025muon,su2025muon} for pre-training a Mixture-of-Experts (MoEs) model with 15.29B total parameters. However, our understanding of its working principle remains largely limited. For instance, the underlying reason for using orthogonalized gradient for the updates of \Muon and why it outperforms \Adam remains elusive. 
        
        \paragraph{Contributions.}
        In this work, we provide theoretical insights into the effectiveness of \Muon and \Adam optimizers through a unifying lens of preconditioning. While \Muon and \Adam can be interpreted as steepest descent with respect to non-Euclidean norms, we instead suggest an alternative view built upon preconditioning. In particular, we explicitly point out two different types of preconditioning for vector and matrix optimization methods: while typical preconditioning aims to reduce the condition number of the Hessian mostly for vector optimization problems, matrix optimization problems indeed can make use of preconditioning that minimizes the condition number of the gradient. Due to such a discrepancy, we argue that the preconditioning of \Adam is mainly derived from the principle of curvature preconditioning mainly for strongly convex vector optimization problems, whereas orthogonalized gradient methods like \Muon perform gradient preconditioning as orthogonal matrices are the best conditioned matrices with condition numbers of 1 \citep{turing1948rounding}. In practical implementation, this preconditioning view also justifies the use of different optimizers for vector and matrix parameters as in the \texttt{modded-nanogpt} repository \citep{modded_nanogpt_2024} where \Muon is used for matrices (except for the embedding and head layers) and \Adam is used for vectors and scalars. We also make various algorithmic contributions which improve \Muon in several aspects. We formulate a class of matrix optimization methods called polar gradient methods (\PolarGrad), which is based on the polar decomposition of the gradient or the momentum with a nuclear norm scaling term derived from steepest descent unlike \Muon and make various comparisons with \Muon. We also propose the use of better numerical polar decomposition algorithms, namely the QDWH \citep{nakatsukasa2010optimizing} and ZOLO-PD \citep{nakatsukasa2016computing} algorithms, which require almost no tuning, unlike the Newton--Schulz iteration in \Muon, and study how the choice of different numerical polar decomposition algorithms affects the efficacy of \PolarGrad through convergence analysis. This makes \PolarGrad a generally applicable class of matrix optimization algorithms for different matrix optimization problems including structured problems like low-rank matrix factorization as well as optimizers for matrix parameters in neural networks. 
        
        \paragraph{Notation.} 
        The $\ell_p$-norm of a vector $x = (x_i)_{1\le j\le d}\in\RR^d$ with $d\in\NN^*\coloneqq\NN\setminus\{0\}$ is denoted by $\norm{x}_p \coloneqq ( \sumd|x_i|^p)^{1/p}$, where $p\in[0,\infty]$. For any $S\in\RR^{d\times d}$, $\tr(S)$ is its trace and $\diag(S)\in\RR^d$ denotes the vector of its diagonal entries. For any $x\in\RR^d$, $\Diag(x)\in\RR^{d\times d}$ is the diagonal matrix with diagonal entries equal to the entries of $x$. For any $A, B\in\RR^{m\times n}$ with $m,n\in\NN^*$, we denote the Frobenius inner product of $A$ and $B$ by $\dotpF{A}{B} \coloneqq \tr(A^\top B)$. For any $A\in\RR^{m\times n}$, we denote its Frobenius norm by $\fronorm{A}$, its nuclear norm by $\nucnorm{A}$, its spectral norm by $\specnorm{A}$, and its ($2$-)condition number by the ratio between its largest and smallest positive singular values $\kappa_2(A)\coloneqq\sigma_{\max}(A)/\sigma_{\min}(A)$. We also denote the set of $m\times n$ semi-orthogonal matrices by $\OO^{m\times n} \coloneqq \{A\in\RR^{m\times n} : A^\top A = I_n \text{ or } AA^\top=I_m\}$, where $I_n$ is the $n\times n$ identity matrix.  Let $\calE$ be a Euclidean space endowed with an inner product $\dotp{\cdot}{\cdot}$ and the induced norm $\|\cdot\|$. The domain of a function $f\colon\calE\to\oRR\coloneqq \RR\cup\{+\infty\}$ is $\dom f \coloneqq \{x\in\calE : f(x)<\infty\}$.     
        The projection of $x$ onto a nonempty closed convex set $\calC$ is denoted by $\proj_{\calC}(x)$.

        \section{Related Work}
        We outline related work on optimizers for deep learning and first-order optimization methods. 
        
        \subsection{Recent Development on Optimizers for Deep Learning}    
        	Distributed \Shampoo \citep{shi2023distributed} achieved the fastest speed-ups among all optimizers in the 2023 \AlgoPerf competition \citep{dahl2023benchmarking,kasimbeg2025accelerating} under the external tuning ruleset, while the self-tuning ruleset is dominated by variants of \AdamW \citep{loshchilov2019decoupled} such as \NAdamW \citep{dozat2016incorporating,medapati2025training} and the winning submission belongs to \textsc{ScheduleFreeAdamW} \citep{defazio2024the}. 
            The discrepancy of the base optimizers for these two rulesets leaves us with a doubt regarding the choice of the most efficient optimizers for neural network training. However, it is also noteworthy that the neural network training tasks in the competition do not include very large foundation models such as large autoregressive decoder-only language models and multi-modal models, which are of more significant interest nowadays. 
            
            The recent success of \Muon \citep{jordan2024muon} has motivated numerous recent variants such as SWAN \citep{ma2024swan}, \Scion \citep{pethick2025training}, COSMOS \citep{liu2025cosmos} and \Gluon \citep{riabinin2025gluon}. While the original development of \Muon \citep{jordan2024muon} is motivated by steepest descent w.r.t.~the spectral norm \citep{bernstein2024old}, it also possesses various interpretations or coincidences with other related methods, including \emph{stochastic spectral descent} \citep{carlson2015stochasticRBM,carlson2016stochastic,carlson2015preconditioned} and  \emph{orthogonalized gradient methods} \citep{tuddenham2022orthogonalising}. It can also be viewed as the \Signum optimizer \citep{bernstein2018signsgd} for matrix parameters where the elementwise sign function is replaced by the matrix sign function.
            
            In addition to the interpretation of \Muon as steepest descent w.r.t.~the spectral norm, the recent work \citep{pethick2025training} interprets gradient orthogonalization as non-Euclidean trust-region optimization, followed by the same interpretation in the work \citep{kovalev2025understanding}. Furthermore, the work \citep{chen2025muon} establishes that \Muon implicitly solves an optimization problem with a spectral norm constraint on weight matrices. These works establish convergence rates for \Muon but are still unable to explain the discrepancy between \Muon and \Adam. That said, a recent work \citep{su2025isotropic} unveils the benefits of gradient orthogonalization as employed in \Muon within one iteration, though it stops short of establishing a convergence rate. We emphasize that it is indeed a matrix preconditioned gradient method that addresses gradient anisotropy. Usually, the condition of the update direction (e.g., the gradient or the momentum) in an iterative algorithm governs its convergence speed (see e.g., Chapter 5 of \citep{bach2024learning}), leading to various preconditioned methods in solving linear systems and iterative algorithms \citep{jambulapati2020fast,qu2025optimal} in order to improve the condition. Adopting this unifying preconditioning viewpoint, we emphasize the substantial difference in the characteristics of the preconditioning of update directions for vector and matrix parameters in neural networks. For vector parameters, preconditioning is usually performed via the multiplication of a matrix preconditioner. For instance, adaptive gradient methods such as \AdaGrad \citep{duchi2011adagrad,mcmahan2010adaptive}, \RMSprop \citep{tieleman2012} and \Adam \citep{kingma2015} can all be viewed as preconditioned methods with diagonal matrix preconditioners, mainly motivated by addressing curvature (or Hessian) anisotropy by approximating the inverse square root of the Hessian through a diagonal matrix. Understanding the Hessian structure of neural networks has been an active area of research that helps understand neural network training; see e.g., \citep{zhang2024transformers,kunstner2024heavy,dong2025towards}. In contrast, preconditioning for matrix parameters is more intricate. Explicit preconditioners for matrix optimization problems might come in pairs, namely left and right preconditioners which are both square matrices, e.g., \Shampoo \citep{gupta2018shampoo} and its variants CASPR \citep{duvvuri2024combining} and SOAP \citep{vyas2024soap}. It turns out that matrix orthogonalization (or semi-orthogonal projection) performs preconditioning without explicit preconditioners. To see this, let us recall that a standard convention in matrix analysis to measure the ``condition'' of a matrix is the ($2$-)condition number, given by $\kappa_2(X)\coloneqq\sigma_{\max}(X)/\sigma_{\min}(X)$, where $\sigma_{\max}(X)$ and $\sigma_{\min}(X)$ are the largest and smallest positive singular values of $X$ respectively. If the update direction has a large condition number, it is called \emph{ill-conditioned} and could lead to slow convergence. Orthogonalization (or more rigorously, a semi-orthogonal projection) of the update direction indeed reduces its condition number to accelerate convergence, since ``the best conditioned matrices are the orthogonal ones, which have condition numbers of 1'' \citep{turing1948rounding}. Taking this preconditioning perspective of matrix parameters for accelerated convergence into account, it is no surprising that Distributed \Shampoo \citep{shi2023distributed} won the external tuning ruleset of the \AlgoPerf competition \citep{dahl2023benchmarking,kasimbeg2025accelerating}, since \Shampoo without preconditioner accumulations is equivalent to \Muon \citep{bernstein2024old}. In contrast, adaptive gradient methods such as \Adam applied to matrix parameters might not enjoy this gradient/momentum preconditioning effect (i.e., might not reduce the condition number of the update direction) since they are derived based on curvature preconditioning via approximating the inverse Hessian and might even lead to undesirable effects such as training instability and loss divergence, illustrated in numerical experiments in \cref{sec:expt}.

            \subsection{Related First-Order Optimization Methods}
            We also go over various related optimization methods, including a very general discussion on steepest descent, as well as matrix optimization methods and various classes of optimizers for deep learning. 
            
            \subsubsection{Steepest Descent Methods}
            We first give a brief overview of the steepest descent method which is at the heart of many first-order methods in mathematical optimization. Let $\calE$ be a Euclidean space endowed with an inner product $\dotp{\cdot}{\cdot}$ and the induced norm $\|\cdot\|$. Let us consider the optimization problem with an objective function $f\colon\calE\to\oRR\coloneqq \RR\cup\{+\infty\}$. Most first-order optimization algorithms can be subsumed as (constrained) steepest descent with respect to a \emph{distance-like function} $\sfd(\cdot, \cdot)$ (see Chapter 9.4 of \citep{boyd2004convex}): 
            \begin{equation}\label{eqn:steepest_descent}
               	(\forall k\in\NN)\quad x_{k+1} \in \argmin_{x\in\calC}\,\tilde{f}(x) \coloneqq f(x_k) + \dotp{\nabla f(x_k)}{x - x_k} + \frac{1}{2\gamma_k}\sfd(x, x_k),  
            \end{equation}
            where $\calC\subseteq\calE$ is a constraint set. 
            Notable examples include gradient descent (GD), preconditioned gradient descent, mirror descent \citep{beck2003mirror,teboulle2018simplified}, proximal splitting algorithms \citep{condat2023proximal} and many others (see e.g., \citep{auslender2006interior} for detailed exposition). 
            However, most adaptive gradient optimizers popular in deep learning cannot be directly expressed in the form of \eqref{eqn:steepest_descent}, including \Adam \citep{kingma2015} and \AdamW \citep{loshchilov2019decoupled}.     	
            While the distance-like function $\sfd$ is mainly chosen to be Euclidean norms in most algorithms, non-Euclidean vector and matrix norms have aroused much attention in recent algorithmic design. For instance, stochastic and preconditioned spectral descent \citep{carlson2015stochasticRBM,carlson2016stochastic,carlson2015preconditioned,hsieh2018non}
            all make use of the spectral norm. The use of non-Euclidean norms in steepest descent can be also found in \citep{flynn2017duality,kelner2014almost}.

            \subsubsection{Matrix Optimization Methods}
            There is a rich literature of matrix optimization algorithms and spectral methods in the field of mathematical optimization, such as eigenvalue optimization \citep{lewis1996eigenvalue,lewis2003mathematics} and proximal methods \citep{benfenati2020proximal}, targeting a wide range of applications in data science \citep{chen2021spectral,chi2019nonconvex,combettes2021fixed}, e.g., structured covariance and precision matrix estimation. However, their applications to deep neural network training remain very limited. In particular, most optimizers for deep learning are based on coordinatewise updates, entailing a vectorization treatment of higher-dimensional parameters (i.e., matrices and tensors) and applications of vector optimization methods. This implies an ignorance of the difference between their underlying algebraic structures. This also leaves a large gap in understanding the proper choice of optimizers for training neural networks consisting of parameters of different algebraic structures.

            \subsubsection{Optimizers for Deep Learning}
            Optimizers for deep learning based on stochastic (sub)gradients  are mainly derived from or at least motivated by various principles from convex optimization theory and algorithms. One main class of such optimizers is viewed as accelerated first-order methods, in which the acceleration is performed via momentum, as well as adaptive learning rate. Another class of popular optimizers belongs to approximate second-order methods, which mainly involve Hessian approximation or Fisher information matrix approximation for natural gradient descent \citep{amari1998natural}. The first class of optimizers is much more popular than the second one, especially for large-scale applications, due to their use of coordinatewise updates which incur much cheaper computational and memory costs. 
        
            \paragraph{Momentum acceleration methods.}
            In classical convex optimization algorithms, the use of momentum, including Polyak's heavy ball \citep{polyak1964some} and Nesterov's accelerated method, \citep{nesterov1983method} is able to accelerate the convergence of gradient descent for convex objectives. Incorporating stochastic gradients with Robbins--Monro's method \citep{robbins1951}, \SGD with Polyak's momentum and Nesterov's accelerated gradient \citep{sutskever2013importance} are developed respectively and are widely used. It is believed that momentum-based methods converge slower than adaptive gradient methods which also consider adaptive learning rates but might generalize better in tasks like image classification.     
            
            \paragraph{Adaptive gradient methods.}
            Adaptive gradient methods are a large class of first-order methods which attempt to adapt learning rates, with a view to achieving better preconditioning. The earliest adaptive gradient method that appeared in the literature is probably \RProp \citep{riedmiller1993direct}, which has motivated other adaptive gradient methods, including \AdaGrad \citep{duchi2011adagrad,mcmahan2010adaptive}, \Adadelta \citep{zeiler2012adadelta}, \RMSprop \citep{tieleman2012}, \Adam \citep{kingma2015}, \Adafactor \citep{shazeer2018adafactor}, \AdaBelief \citep{zhuang2020adabelief}, \Lion \citep{chen2023symbolic}, \Sophia \citep{liu2024sophia}, etc.
            However, the interpretation of adaptive learning rates for adaptive gradient methods is not the only one in the literature. For instance, \Adam \citep{kingma2015} can be viewed as a form of smoothed sign descent (\signSGD and \Signum) \citep{bernstein2018signsgd,balles2018dissecting}, which is equivalent to (normalized) steepest descent with respect to the $\ell_\infty$-norm. 
            
            \paragraph{Approximate second-order methods.}
            Motivated by second-order optimization methods which converge much faster than first-order methods on strongly convex problems, various deep learning optimizers were developed based on the principle of Hessian approximation or preconditioner approximation, particularly with layerwise Kronecker-factored preconditioners, including K-FAC \citep{martens2015optimizing}, \Shampoo \citep{gupta2018shampoo,anil2020scalable}, BFGS and L-BFGS \citep{goldfarb2020practical}, CASPR \citep{duvvuri2024combining} and SOAP \citep{vyas2024soap}, as well as learned preconditioners in preconditioned SGD (PSGD) \citep{li2017preconditioned,pooladzandi2024curvature}. While the inverse Hessian is understood as a good preconditioner for strongly convex optimization problems, it remains elusive to understand the performance of optimizers based on its diagonal approximations and layerwise Kronecker-factored preconditioners for nonconvex problems other than purely technical convergence analysis. We point out the insufficiency of the Hessian approximation and Kronecker-factored preconditioning viewpoints of \Shampoo \citep{morwani2024new}, since diagonal approximations might worsen the preconditioning effect and the Kronecker-factored structure might not hold at all for most neural networks. In contrast, we advocate for an understanding of deep learning optimizers via the intrinsic working principle of these preconditioned gradient methods---reducing the ill-conditionedness of the Hessian or the anisotropy of the gradient. 
            
            One-sided \Shampoo \citep{xie2025structured,an2025asgo} only uses the left preconditioner which potentially saves memory, whereas preconditioned Riemannian gradient descent (RPGD) \citep{bian2024preconditioned} further replaces the left and right preconditioners with their diagonal approximations. CASPR \citep{duvvuri2024combining} and SOAP \citep{vyas2024soap} are two other notable improved variants of \Shampoo that also have explicit preconditioners. The left and right preconditioners in \Shampoo take a total memory requirement of $\scrO(m^2+n^2) \gg \scrO(mn)$ for large $m$ and $n$, which are prohibitive for training very large layers in large-scale pre-training. Besides, without more advanced numerical linear algebra algorithms, \Shampoo and its variants with explicit preconditioners in such form cannot be easily parallelized and require high precision due to the involved matrix inverse roots. In contrast, \Muon and its variants based on semi-orthogonal projections do not involve any explicit preconditioners and matrix inverse operations, making them suitable for parallelization. As model size grows, we are often more memory-bound than compute-bound, making implicit preconditioners more plausible.

        \section{Polar Gradient Methods}
        \label{sec:polar_grad}
        Our development of polar gradient methods is largely motivated by \Muon and related orthogonalized gradient methods, which we detail below. 
        
        \subsection{\Muon and Orthogonalized Gradient Methods}
        We first recover the connection between the steepest descent and the matrix sign descent interpretations \citep{jordan2024muon,bernstein2024old,su2024muon} of orthogonalized gradient methods \citep{tuddenham2022orthogonalising}. The matrix sign function on real rectangular matrices can be defined through its singular value decomposition (SVD). If $U\Sigma V^\top = \SVD(X)$ is the SVD of $X\in\RR^{m\times n}$, then the \emph{matrix sign function} of $X$ is defined by $\msgn(X) \coloneqq UV^\top$. 
        
        Note that this definition is a slight abuse of notion and is different from that in the numerical linear algebra literature such as the one in Chapter 5 of \citep{higham2008functions}, which is only defined for square matrices. The matrix sign function defined above should be better referred to as the \emph{orthogonal polar factor} arising from the \emph{polar decomposition} (see \Cref{sec:polar_grad}).     
        It turns out that under the above definition, the matrix sign function of $X\in\RR^{m\times n}$ is equivalent to the projection of $X$ onto the space of $m\times n$ semi-orthogonal matrices $\OO^{m\times n}$ in any unitarily invariant norm $\matsnorm{\cdot}{}$, i.e., $\proj_{\OO^{m\times n}}(X)\coloneqq \argmin_{O\in\OO^{m\times n}}\,\matsnorm{O - X}{}$ (see \Cref{thm:polar}). 
        \Muon without momentum or stochastic spectral descent (SSD) can be interpreted as (resp., normalized and unnormalized) stochastic steepest descent w.r.t.~the spectral norm, as illustrated below.     
        Let $\sff\colon\RR^{m\times n}\to\oRR$ be a (possibly nonconvex) objective function and consider the stochastic optimization problem minimizing $\sff(X) \coloneqq \Ex_{\xi\sim\sfP}[\sff(X, \xi)]$. We then denote a stochastic gradient of $\sff$ at $X_k$ with the sample $\xi_k$ by $G_k = \nabla \sff(X_k, \xi_k)$. Then, \Muon without momentum \citep{jordan2024muon} or stochastic spectral descent \citep{carlson2015stochasticRBM,carlson2016stochastic,bernstein2024old} can be derived by solving the following subproblem at every iteration: 
        \begin{equation}\label{eqn:ssd}
           	(\forall k\in\NN)\quad X_{k+1} \in \argmin_{X\in\RR^{m\times n}} \, \left\{\dotpF{G_k}{X - X_k} + \frac{1}{2\gamma_k}\specnorm{X-X_k}^2\right\}. 
        \end{equation}
        Note that, since the spectral norm is non-differentiable (nonsmooth), the right-hand side of \eqref{eqn:ssd} might not be a singleton. Indeed, the subdifferential is a singleton if and only if $G_k$ is of full rank. Then, \eqref{eqn:ssd} takes the following closed-form update: 
        \begin{equation}\label{eqn:mat_sign_descent}
           	(\forall k\in\NN)\quad X_{k+1} = X_k - \gamma_k \nucnorm{G_k}\cdot\msgn(G_k). 
        \end{equation}
        If $G_k$ is not of full rank, \eqref{eqn:mat_sign_descent} becomes a stochastic subgradient method.     
        \Muon can be derived by simply introducing the momentum either in the form of $M_k = \mu M_{k-1} + G_k$ with $\mu>0$ or $M_k = \beta M_{k-1} + (1-\beta)G_k$ with $\beta\in(0,1)$ and replacing $G_k$ in \eqref{eqn:mat_sign_descent} by $M_k$. Note that the nuclear norm term in \eqref{eqn:mat_sign_descent} does not appear in \Muon, which is investigated in detail in \Cref{sec:polar_grad}.

        \subsection{Connection to Polar Decomposition}        
        The term ``orthogonalized gradient'' could be confusing since the matrix sign function is not equivalent to the orthonormal matrix obtained from the QR decomposition, but its semi-orthogonal projection instead (see also \Cref{thm:polar}). To avoid this confusion and the proper use of terminology, we now introduce the \emph{polar decomposition} of matrices \citep{autonne1902sur,higham1986computing}. 
        \begin{definition}[Polar decomposition]\label{def:polar_decomp}
        	Any matrix $A\in\RR^{m\times n}$ with $m\ge n$ (resp.~$m<n$) has a polar decomposition $A = U_\sfp H$ (resp.~$A = HU_\sfp$), where the \emph{orthogonal polar factor} $U_\sfp\in\OO^{m\times n}$ has orthonormal columns (resp.~rows) and the \emph{symmetric polar factor} $H\in\bbS_+^n$ (resp.~$H\in\bbS_+^m$) is a symmetric positive semidefinite matrix. The matrix $H$ is unique, and $U_\sfp$ is unique if $A$ has full rank. We write $U_\sfp H = \polar(A)$ as the polar decomposition of $A$. 
        \end{definition}
        Note that, if $U\Sigma V^\top=\SVD(A)$, then  $U_\sfp H = \polar(A)$ (resp.~$HU_\sfp = \polar(A)$) can also be represented by $U_\sfp = UV^\top = \msgn(A)$ and $H = V\Sigma V^\top$ (resp.~$H = U\Sigma U^\top$). Therefore, we can compute the matrix sign function of $A$ using its orthogonal polar factor \citep{higham1994matrix}.     
        Since the orthogonal polar factor $U_\sfp = \msgn(A)$ can almost be uniquely determined for any matrix $A\in\RR^{m\times n}$, we coin this class of matrix optimization methods based on the polar decomposition of the gradient as \emph{polar gradient methods} (\PolarGrad). Despite its similarities to orthogonalized gradient methods such as \Muon, we emphasize that \PolarGrad also makes use of the symmetric polar factor $H$ and the potential usage of more advanced numerical polar decomposition algorithms than the Newton--Schulz iteration, hence necessitating its own name to refer to a broader class of matrix optimization methods based on the polar decomposition of the gradient or momentum.

    	\subsection{Polar-Decomposed Gradient with Nuclear Norm Scaling} 
    	\label{subsec:nuc_norm}
        Recall that the the orthogonal polar factor of the gradient performs gradient-anisotropy preconditioning (cf.~\cref{subsec:precond}).  
        However, (almost) perfect gradient preconditioning via orthogonal polar factors preserves only directional information via singular vectors and removes curvature adaptation provided by singular values, which is crucial for fast optimization. 
        In the original implementation of \Muon \citep{jordan2024muon}, a scaling factor of $\sqrt{\max\{1, m/n\}}$ is used, while in \citep{liu2025muon} a scaling factor of $\sqrt{\max\{m, n\}}$ is used. \Scion \citep{pethick2025training}, a close variant of \Muon, instead adopts a scaling factor of $\sqrt{m/n}$ and leads to hyperparameter transfer. However, these choices can only address the sizes of different weight matrices in neural networks, hence not being adaptive across different iterations. 
        
        In contrast, the learning rate should be scaled adaptively based on the actual gradient magnitude using the nuclear norm of the gradient as in \eqref{eqn:mat_sign_descent}, as opposed to the original form of \Muon \citep{jordan2024muon} and the analysis of \Muon as a non-Euclidean trust-region gradient method \citep{pethick2025training,kovalev2025understanding}. 
        Such methods would converge faster than pure polar gradient updates, as in \Muon by providing curvature sensitivity via nuclear norm scaling while retaining isotropy advantages. 
        Here, we also mention an intimate relationship between the nuclear norm $\nucnorm{G}$ and the symmetric polar factor $H$. Without loss of generality, we assume that the gradient $G\coloneqq\nabla\sff(X)\in\RR^{m\times n}$ with $m\ge n$. If $U\Sigma V^\top=\SVD(G)$, then $\nucnorm{G}=\tr(\Sigma)$. We also recall that for $U_\sfp H=\polar(G)$, $H=V\Sigma V^\top$, so $\tr(H) = \tr(V\Sigma V^\top) = \tr(V^\top V\Sigma) = \tr(\Sigma)$ since $V$ is orthogonal.     
        Therefore, the unnormalized matrix sign descent \eqref{eqn:mat_sign_descent} can be explicitly written in terms of the two polar factors of the gradient, leading to vanilla \PolarGrad:     
        \begin{equation}\label{eqn:polargrad}
        	U_kH_k = \polar(G_k),\quad X_{k+1} = X_k - \gamma_k \,\tr(H_k)\, U_k, 
        \end{equation}
        where $G_k$ represents a deterministic gradient $\nabla\sff(X_k)$ or a stochastic gradient $ \nabla \sff(X_k, \xi_k)$ with a sample $\xi_k$, and $\gamma_k>0$ is a learning rate independent of $X_k$, $H_k$ and $U_k$. 
        \PolarGrad with exponential moving average (EMA) momentum and decoupled weight decay (\textsc{PolarGradM(W)}), similar to \Muon (henceforth \PolarMuon), is given by: 
        \begin{equation*}
        	M_k = \beta M_{k-1} + (1-\beta)G_k,\quad U_kH_k = \polar(M_k),\quad X_{k+1} = (1-\lambda\gamma_k)X_k - \gamma_k \,\tr(H_k)\, U_k. 
        \end{equation*}
        \PolarMuon is only one of the possible ways to introduce EMA momentum to \PolarGrad, which performs a momentum update before the polar decomposition of momentum (henceforth \emph{momentum-first}). We can also perform the polar decomposition of the gradient and perform a momentum update afterward (henceforth \emph{polar-first}) as follows: 
        \begin{equation*}
        	U_kH_k = \polar(G_k),\quad M_k = \beta M_{k-1} + (1-\beta)U_k,\quad X_{k+1} = (1-\lambda\gamma_k)X_k - \gamma_k \,\tr(H_k)\, M_k. 
        \end{equation*}
        As we will see in the next subsection, the inclusion of the nuclear norm scaling term is able to improve the convergence rate from sublinear to linear for deterministic strongly convex objectives. We also observe this empirically for a nonconvex low-rank matrix completion example in \Cref{subsec:mat_fac}.

    	\subsection{Comparison with \Muon}
    	\label{subsec:comparison_muon}
    	The nuclear norm scaling factor, $\tr(H_k)$, in the \PolarGrad update \eqref{eqn:polargrad} leads to a pivotal distinction from the original \Muon optimizer. As shown in \Cref{subsec:nuc_norm}, this scaling arises naturally from the steepest descent formulation with respect to the spectral norm. Beyond this derivation, the inclusion of $\tr(H_k)$ confers a crucial property that we term \emph{null-gradient consistency}.

        \subsubsection{Null-Gradient Consistency}
        We define \emph{null-gradient consistency} below. 
    	\begin{definition}[Null-gradient consistency]
    	\label{def:null-grad_consistency}
    	An optimization algorithm exhibits null-gradient consistency if the magnitude of its update step tends to zero as the effective gradient term approaches zero.
    	\end{definition}
    	
    	While not a strict mathematical prerequisite for all optimization methods, null-gradient consistency is a desirable characteristic. It ensures that the optimizer's parameter changes diminish as the gradient indicating the direction of descent vanishes. This behavior is conducive to identifying convergence to stationary points and for maintaining a consistent interpretation of the learning rate's role throughout the optimization process.
    	
    	Now, consider the behavior of \Muon and \PolarGrad in the vicinity of a point where the effective gradient $G_k$ (or $M_k$, if momentum is used) is very small, i.e., $G_k \approx 0$. In the standard \Muon update, the step is proportional to $\msgn(G_k)$. Even as $G_k \to 0$ (but $G_k \neq 0$), $\msgn(G_k)$ remains a semi-orthogonal matrix, whose magnitude does not diminish to zero as $G_k$ itself vanishes. Consequently, the magnitude of the update direction provided by $\msgn(G_k)$ does not tend to zero. Thus, \Muon---at least in its original formulation---does not satisfy the null-gradient consistency property. This can lead to persistent updates or oscillations around an optimum where the true gradient is negligible, unless the learning rate is meticulously adjusted or decayed.
    	
    	In contrast, for the \PolarGrad update, the scaling factor $\tr(H_k)$ is equivalent to the nuclear norm of $G_k$. As $G_k \to 0$, its nuclear norm, and therefore $\tr(H_k)$, also tends to zero. Thus, the entire update term $\gamma_k \,\tr(H_k) U_k$ vanishes as $G_k \to 0$, ensuring that \PolarGrad satisfies the null-gradient consistency property. The satisfaction of this property by \PolarGrad suggests more stable behavior, particularly in later stages of optimization where true gradients are typically small. 
    	
    	It is worth emphasizing that we present the property of null-gradient consistency in a conceptual, rather than a mathematically formal, manner. When evaluating whether an optimizer satisfies this property, we exclude exogenous terms such as decoupled weight decay. Furthermore, the effective gradient should be understood as the modified gradient that ultimately dictates the update magnitude---for instance, the momentum gradient, rather than the raw gradient.

        \subsubsection{Recovering \PolarGrad from \Muon with Armijo's Backtracking Line Search}
        \label{subsubsec:line_search}
        In most deep learning applications, we emphasize that learning rate sequences (or schedules) are usually independent of the iterates and specified prior to model training. As a consequence, it is almost impossible to handpick learning rate sequences that absorb the nuclear norm scaling of the matrix gradient or momentum without any iterate-dependent information. This entails a noted difference from optimizers based on the Linear Minimization Oracle (LMO) optimization framework, such as \Muon \citep{jordan2024muon}, \Scion \citep{pethick2025training} and \Gluon \citep{riabinin2025gluon}, whose learning rates could be dimension-dependent but iterate-independent. 
        
        On the other hand, popularly used for gradient descent for (unconstrained) convex optimization, Armijo's backtracking line search \citep{armijo1966minimization} is a line search method to find the (iterate-dependent) learning rate of each iteration, requiring that the objective function is differentiable and its gradient is available. Let us recall that Armijo's backtracking line search determines the learning rate $\alpha_k>0$ of \Muon without momentum such that 
        \[f(X_k - \alpha_k U_k) \le f(X_k) - c\alpha_k \dotpF{G_k}{U_k} = f(X_k) - c\alpha_k\nucnorm{G_k}, \]
        where $c\in(0,1)$ is a selected control parameter, $G_k\coloneqq\nabla f(X_k)$ is the gradient and $U_k$ is the orthogonal polar factor of $G_k$. Furthermore, if $f$ is $L$-Lipschitz smooth (see \Cref{def:Lipschitz_smooth} below), then we have
        \[f(X_k - \alpha_k U_k) \le f(X_k) - \alpha_k\nucnorm{G_k} + \frac{L}{2}\alpha_k^2 r_k, \]
        where $r_k\coloneqq\rank(G_k) = \fronorm{U_k}^2$ (see Proof of \Cref{thm:polargrad_strcvx} in \Cref{sec:proofs} for its proof). The Armijo's condition and the $L$-Lipschitz smoothness assumption together yield 
        \[\alpha_k \le \frac{2(1-c)}{Lr_k} \nucnorm{G_k}. \]
        Consequently, the backtracking line search procedure picks $\alpha_k$ so that $\alpha_k/\nucnorm{G_k}$ stays in a stable range, so it turns out that the nuclear norm scaling term will be recovered. We however make the nuclear norm scaling term explicit in \PolarGrad as opposed to \Muon or \Scion since backtracking line search procedures for learning rates are almost never used in deep learning potentially due to the extra computation and implementation complication required. We also remark that when $c=1/2$, we obtain $\alpha_k\le\nucnorm{G_k}/(Lr_k)$ which recovers the choice of $\gamma_k=1/(Lr_k)$ in \Cref{thm:polargrad_strcvx} in the following subsection.

    	\subsection{Convergence Analysis of \PolarGrad with Exact Polar Factors}
    	\label{subsec:conv}
        To better characterize the convergence behavior of the optimizers in the \PolarGrad family, we derive their convergence rates in terms of the gradient condition number $\kappa_G$ and the Hessian condition number $\kappa_H$. Without loss of generality, we assume that the optimization variable $X\in\RR^{m\times n}$ has dimensions $m\ge n$. We do not consider any weight decay. We emphasize that there are several works \citep{li2025muon,an2025asgo,kovalev2025understanding,pethick2025training,shen2025convergence} that analyze the convergence of \Muon, but we emphasize the difference between \PolarGrad and \Muon---the inclusion of the nuclear norm term. We first derive the convergence rates of \PolarGrad with deterministic gradients for Lipschitz smooth and strongly convex functions.   
        In what follows, we denote the deterministic or full gradient $G_k \coloneqq \nabla f(X_k)$ and the stochastic gradient $\hG_k \coloneqq \nabla f(X_k, \xi_k)$. 
        We first recall the following standard results for functions satisfying $L$-Lipschitz smoothness and $\mu$-strong convexity. 
        \begin{definition}[$L$-Lipschitz smoothness]\label{def:Lipschitz_smooth}
              	Let $f\colon\RR^{m\times n}\to\oRR$ be $L$-Lipschitz smooth, i.e., there exists a constant $L\in(0,\infty)$ such that 
              	\[(\forall(X,Y)\in\RR^{m\times n}\times\RR^{m\times n})\quad\fronorm{\nabla f(X) - \nabla f(Y)} \le L \fronorm{X-Y}. \]
              	Then, equivalently, we have
              	\[(\forall(X,Y)\in\RR^{m\times n}\times\RR^{m\times n})\quad f(Y) \le f(X) + \dotpF{\nabla f(X)}{Y- X} + \frac{L}{2}\fronorm{Y - X}^2. \]
              	Furthermore, we also have 
              	\begin{equation*}
              		(\forall X\in\RR^{m\times n})\quad\fronorm{\nabla f(X)}^2 \le 2L\left(f(X) - \fstar\right).  
              	\end{equation*}
        \end{definition}
        
        We also state the following result that strong convexity implies the Polyak--Łojasiewicz (\PL) condition \citep{karimi2016linear}. 
        \begin{proposition}[$\mu$-strong convexity]\label{thm:strong_cvx}
              	Let $f\colon\RR^{m\times n}\to\oRR$ be $\mu$-strongly convex, i.e., there exists a constant $\mu\in(0,\infty)$ such that 
              	\[(\forall(X,Y)\in\RR^{m\times n}\times\RR^{m\times n})\quad\dotpF{\nabla f(X) - \nabla f(Y)}{X - Y} \ge \mu\fronorm{X - Y}^2, \]
              	or equivalently, 
              	\[(\forall(X,Y)\in\RR^{m\times n}\times\RR^{m\times n})\quad f(Y) \ge f(X) + \dotpF{\nabla f(X)}{Y- X} + \frac{\mu}{2}\fronorm{Y - X}^2. \]
              	Note that $\mu$-strong convexity implies the $\mu$-Polyak--Łojasiewicz (\PL) condition or inequality: 
              	\begin{equation}\label{eqn:strong_cvx}
              		(\forall X\in\RR^{m\times n})\quad\fronorm{\nabla f(X)}^2 \ge 2\mu\left(f(X) - \fstar\right),  
              	\end{equation}
              	where $\fstar\coloneqq \min f$. 
              	Functions satisfying \eqref{eqn:strong_cvx} are called $\mu$-Polyak--Łojasiewicz (\PL) functions. Therefore, the \PL condition is a more relaxed condition than strong convexity (without assuming any convexity). 
        \end{proposition}

        Indeed, in the following convergence analysis, it suffices to assume the \PL condition instead of strong convexity. We now make the following assumption, defining some related notions. 
        \begin{assumption}\label{assum:base}
        	We assume that the objective function $f\colon \RR^{m \times n} \to \oRR$ is $L$-Lipschitz smooth and a $\mu$-\PL function. Let $\fstar = \min f$, and we define $r_k \coloneqq \rank(\nabla f(X_k))$ and $r_{\max}\coloneqq\max_{k\in\setK} r_k \le\min\{m,n\}$. Let $\sigma_{1_k} \ge \dots \ge \sigma_{r_k} > 0$ be the singular values of $\nabla f(X_k)$. We also define the gradient condition number $\kappa_{G_k} \coloneqq \sigma_1(\nabla f(X_k)) / \sigma_{r_k}(\nabla f(X_k))$ and the (global) Hessian condition number $\kappa_H\coloneqq L/\mu$.    
        \end{assumption}
        
        Under the above assumptions, we state our first theoretical result. 
        \begin{theorem}[\PolarGrad]\label{thm:polargrad_strcvx}
        	Suppose that \Cref{assum:base} holds. For a learning rate sequence $\gamma_k = 1/(Lr_k)$, the iterates of \PolarGrad \eqref{eqn:polargrad} satisfy $f(X_{k+1}) - \fstar \le \left(1 - 1/( r_k\kappa_H)\right)(f(X_k) - \fstar)$ and $f(X_{k+1}) - \fstar \le \left(1 - 1/(\kappa_{G_k}^2\kappa_H) \right)(f(X_k) - \fstar)$, respectively. 
        \end{theorem}

        Consequently, this implies that the gradient-based rate can significantly outperform the Hessian-based rate when $\kappa_{G_k}^2 \ll r_k$, i.e., when the gradient is well-conditioned even if the Hessian is poorly conditioned. This situation could arise in structured matrix problems (e.g., matrix factorization). While the rank $r_k$ is usually not known in practice, we can use $r_{\max}$ at each iteration and obtain a uniform rate of convergence of $\scrO(\exp(-k/(r_{\max}\kappa_H)))$ with a constant learning rate. 
        In such case, the convergence rate also becomes dimension-dependent. 
        
        To distinguish the algorithms with deterministic gradients, we use \PolarSGD to refer to the stochastic gradient counterpart of \PolarGrad. We now derive the convergence rates of \PolarSGD under the following additional bounded gradient variance assumptions on the stochastic gradient. 
        \begin{assumption}\label{assum:grad_noise}
        	For any $X\in\RR^{m\times n}$ and sample $\xi\sim\calD$, the stochastic gradient $\nabla f(X, \xi)$ is unbiased, i.e., $\Ex_{\xi\sim\calD}[\nabla f(X, \xi)] = \nabla f(X)$, and has bounded variance, i.e., $\Ex_{\xi\sim\calD}[\fronorm{\nabla f(X, \xi) - \nabla f(X)}^2] \le \varsigma^2$ for some $\varsigma\in(0,\infty)$.  
        \end{assumption}
        
        \begin{theorem}[\PolarSGD]\label{thm:polarsgd_strcvx}
        	Suppose that \Cref{assum:base,assum:grad_noise} hold. For a constant learning rate $\gamma \in\left(0,  1/(Lr_{\max}^2)\right]$, the iterates of \PolarSGD satisfy $\Ex[f(X_k) - \fstar] \le \scrO\left(\exp(-C_1k) + C_2\varsigma^2\right)$, where $C_1$ and $C_2$ are constants depending on $L$, $\mu$, $\gamma$ and $r_{\max}$. 
        \end{theorem}

        Since \PolarSGD is similar to matrix \signSGD except for the inclusion of the nuclear norm scaling term, we are also interested in how their convergence rates compare, as well as those of their deterministic gradient counterpart \PolarGrad and matrix sign descent. 
        \begin{theorem}[Matrix sign descent and matrix \signSGD]\label{thm:mat_signSGD_strcvx}
        	Suppose that \Cref{assum:base} holds. 
        	With a constant learning rate $\gamma>0$, the iterates of matrix sign descent $X_{k+1} = X_k - \gamma U_k$ with $U_kH_k = \polar(\nabla f(X_k))$ satisfy a nonlinear recursion $\Delta_{k+1} \le \Delta_k - \gamma\sqrt{2\mu\Delta_k} + \frac{L}{2}\gamma^2r_{\max}$ which converges at most sublinearly at a floor, where $\Delta_k\coloneqq f(X_k) - \fstar$ is the optimality gap. 
        	On the other hand, for a general $L$-Lipschitz smooth but possibly nonconvex objective function $f\colon \RR^{m \times n} \to \oRR$, the iterates of matrix sign descent $(X_k)_{k\in\setK}$ satisfy $\min_{k\in\setK} \fronorm{\nabla f(X_k)} \le \scrO(1/(\gamma K) + L\gamma r_{\max}/2)$, and the iterates of matrix \signSGD $X_{k+1} = X_k - \gamma\hU_k$ with $\hU_k\hH_k = \polar(\nabla f(X_k, \xi_k))$ satisfy $\min_{k\in\setK}\Ex\fronorm{\nabla f(X_k)} \le \scrO\left(1/(\gamma K) + L\gamma r_{\max}/2 + \varsigma\sqrt{r_{\max}}\right)$ if \Cref{assum:grad_noise} also holds. 
        \end{theorem}
        Thus, if the learning rate is constant, convergence plateaus at a floor, implying that learning rate decay is necessary for \PolarSGD, matrix sign descent and matrix \signSGD even for strongly convex objectives.    	
        Similar results for \PolarSGDM, \Muon and non-\PL objectives are more technically involved and left for future work, but we empirically evaluate them in \Cref{sec:expt}.

        \subsection{Improving \Muon with Better Numerical Polar Decomposition Algorithms}
        Computing the nuclear norm from scratch requires a full SVD and could be computationally expensive, but it can be computed via the identity $\nucnorm{G_k} \equiv \dotpF{G_k}{\msgn(G_k)}$ due to the dual-norm relationship of the spectral norm and the nuclear norm (see \Cref{prop:subdiff_dual}). 
        The practical performance of \PolarGrad and \Muon highly relies on the involved numerical polar decomposition algorithm. 
        \Muon uses the Newton--Schulz (NS) iteration \citep{higham2008functions} to compute the orthogonal polar factor, but it requires careful choice of the matrix iterative polynomial coefficients for fast convergence. A dynamic coefficient schedule is used in GPT-2 Medium in the \texttt{modded-nanogpt} repository \citep{modded_nanogpt_2024,su2025msign}, different from the fixed coefficients used for GPT-2 Small. Tedious coefficient tuning would thus be needed for training different neural networks, preventing \Muon from being a general drop-in replacement of \Adam for any matrix parameters in neural networks. 
        
        Developing efficient polar decomposition algorithms has been a crucial research area in numerical linear algebra; see e.g., \citep{higham1986computing,higham1990fast,higham2008functions} for earlier works. In a series of work, Nakatsukasa and co-authors \citep{nakatsukasa2010optimizing,nakatsukasa2013stable,nakatsukasa2016computing} have developed various polar decomposition algorithms which provably converge much faster than the NS iteration and other standard approaches (in terms of the number of iterations) and are more numerically stable, namely the QR-based Dynamically Weighted Halley (QDWH) algorithm \citep{nakatsukasa2010optimizing} and the ZOLO-based Polar Decomposition (ZOLO-PD) algorithm \citep{nakatsukasa2016computing}, basically developed based on dynamic coefficient schedules with rational approximations as opposed to the fixed coefficients with polynomial approximations in the NS iteration. In particular, when the matrix is very ``fat'' such as the embedding and the classification head weights in language models, the NS iteration might fail to converge due to ill-conditioned initializations, thus prohibiting the use of \Muon. The implementations of these two algorithms are lacking in deep learning libraries except for the QDWH algorithm in JAX \citep{jax2018github}, despite their high-performance CPU implementation \citep{ltaief2019massively}. More recently, the work \citep{amsel2025polar} introduces the \textsc{Polar Express}, which is a new GPU-efficient numerical polar decomposition algorithm, inspired by the works \citep{chen2014stable,nakatsukasa2016computing}. Likewise, the work \citep{grishina2025accelerating} proposes CANS, both of which attempt to accelerate the NS iteration by optimizing the coefficients of the matrix iterative polynomial in the NS iteration. Further discussion on numerical polar decomposition algorithms is given in \Cref{subsec:NLA}.

        \subsection{Convergence Analysis of \PolarGrad with Inexact Polar Oracles}
        \label{subsec:conv_inexact}
        The convergence analysis in \Cref{subsec:conv} implicitly assumes that the orthogonal polar factor $U_k$ is exact at each iteration $k\in\NN$, which is unrealistic in practice since it is obtained by a numerical algorithm with incurred inaccuracy. Almost all existing theoretical analyses of \Muon such as \citep{li2025muon,shen2025convergence,chen2025muon} are also established under the same assumption. We now relax this assumption, only assuming access to an \emph{inexact polar oracle} $\polarhat$ which provides approximate orthogonal and symmetric polar factors $(\tU_k, \tH_k)$ at each iteration $k\in\NN$, in order to better characterize the convergence behavior of the realized optimizers in the \PolarGrad family. Now, the realized algorithm for \PolarGrad \eqref{eqn:polargrad} becomes
        \begin{equation}\label{eqn:polargrad_practical}
           	\tU_k \tH_k = \polarhat(G_k),\quad X_{k+1} = X_k - \gamma_k \tnu_k\tU_k, 
        \end{equation}
        where the nuclear norm scaling is computed using $\tnu_k \coloneqq\dotpF{\tU_k}{G_k}$ instead of $\nu_k=\dotpF{U_k}{G_k}$. The realized algorithms for other optimizers in the \PolarGrad family are likewise defined. 
        
        We first study the convergence rates of these algorithms with access to a general inexact polar oracle which satisfies the following assumption.     
        \begin{assumption}\label{assum:inexact_polar_oracle}
            At each iteration of the optimizers in the \PolarGrad family, we only assume access to an \emph{inexact polar oracle} $\polarhat$ which provides a pair of approximate orthogonal and symmetric polar factors $(\tU_k, \tH_k)$ of the (deterministic or stochastic) gradient $G_k$ at each iteration $k\in\NN$, satisfying the following conditions: (i) $\specnorm{\tU_k - U_k} \le \varepsilon_k$ for some $\varepsilon_k\in\left[0,1\right)$; (ii) $\specnorm{\tU_k^\top\tU_k - I} = \scrO(\delta_k)$ for some $\delta_k\ge0$ where $r_k\coloneqq\rank(G_k)$. We also define $\varepsilon_{\max}\coloneq\sup_{k\in\setK} \varepsilon_k$ and $\delta_{\max} \coloneqq \sup_{k\in\setK} \delta_k$, and recall that $r_{\max}\coloneqq\max_{k\in\setK} r_k \le\min\{m,n\}$. 
        \end{assumption}
        The first condition is an error bound of the approximate orthogonal polar factor in the spectral norm, while the second condition implies that $\fronorm{\tU_k}^2 \le r_k(1+\delta_k)$ and is closely related to the backward stability of the concerned polar decomposition provided by the inexact polar oracle. 
        With the above additional assumption, we obtain a convergence rate for \PolarGrad with general inexact polar oracles as follows. 
        
        \begin{theorem}[\PolarGrad with general inexact polar oracles]\label{thm:polargrad_strcvx_inexact}
            Suppose that \Cref{assum:base,assum:inexact_polar_oracle} hold. For a constant learning rate $\gamma \coloneqq c/(Lr_{\max}(1 + \delta_{\max}))$ for some $c\in\left(0,1\right]$ for all $k\in\NN$, the iterates of realized \PolarGrad \eqref{eqn:polargrad_practical} satisfy 
            \[
            f(X_{k+1}) - \fstar \le \left(1 - \frac{2c}{r_{\max}\kappa_H}\left(1-\frac{c}2\right) \frac{(1-\varepsilon_{\max})^2}{1+\delta_{\max}} \right)(f(X_k) - \fstar).
            \]
        \end{theorem}
        From the above theorem, if we set $c=1$, we can deduce that the convergence rate of \PolarGrad with general inexact polar oracles is slowed down by a factor of $(1+\delta_{\max})/(1-\varepsilon_{\max})^2$ compared to that of the exact \PolarGrad in \Cref{thm:polargrad_strcvx}. 
    	
    	For stochastic gradient $\hG_k\coloneqq\nabla f(X_k, \xi_k)$, we use alternative notation for \Cref{assum:inexact_polar_oracle}. We write $\hU_k\hH_k = \polar(\hG_k)$ for the exact polar decomposition of $\hG_k$, and $\tU_k \tH_k = \polarhat(\hG_k)$ for its inexact counterpart. Then \Cref{assum:inexact_polar_oracle} becomes (i) $\specnorm{\tU_k - \hU_k} \le \hat\varepsilon_k$ for some $\hat\varepsilon_k\in\left[0,1\right)$; (ii) $\specnorm{\tU_k^\top\tU_k - I} = \scrO(\hat\delta_k)$ for some $\hat\delta_k\ge0$ where $\hat{r}_k\coloneqq\rank(\hG_k)$. The constants $\hat\varepsilon_{\max}$, $\hat\delta_{\max}$ and $\hat{r}_{\max}$ are defined similarly. 
        \begin{theorem}[\PolarSGD with general inexact polar oracles]\label{thm:polarsgd_strcvx_inexact}
            Suppose that \Cref{assum:base,assum:grad_noise,assum:inexact_polar_oracle} hold. For a constant learning rate $\gamma \in\left(0, (1-\hat\varepsilon_{\max})^2/(L\hat{r}_{\max}^2(1+\hat\delta_{\max})^2)\right]$, the iterates of \PolarSGD satisfy 
            $\Ex[f(X_k) - \fstar] \le \scrO\left(\exp(-\tilde{C}_1k) + \tilde{C}_2\varsigma^2\right)$, where $\tilde{C}_1$ and $\tilde{C}_2$ are constants depending on $L$, $\mu$, $\gamma$, $\hat\varepsilon_{\max}$, $\hat\delta_{\max}$ and $\hat{r}_{\max}$. 
        \end{theorem}
        Likewise, from the above theorem, the convergence rate of \PolarSGD with general inexact polar oracles is slowed down by a factor of $(1+\hat\delta_{\max})^2/(1-\hat\varepsilon_{\max})^4$ compared to that of the exact \PolarSGD in \Cref{thm:polarsgd_strcvx}. 
            
        We also derive the corresponding convergence rates without the nuclear norm scaling, i.e., matrix sign descent and matrix \signSGD without assuming the $\mu$-\PL condition. 
        \begin{theorem}[Matrix sign descent and matrix \signSGD with general inexact polar oracles]\label{thm:mat_signSGD_strcvx_inexact}
            Suppose that \Cref{assum:base,assum:inexact_polar_oracle} holds. 
            With a constant learning rate $\gamma>0$, the iterates of matrix sign descent $X_{k+1} = X_k - \gamma \tU_k$ satisfy a nonlinear recursion $\Delta_{k+1} \le \Delta_k - \gamma(1-\varepsilon_{\max})\sqrt{2\mu\Delta_k} + \frac{L}{2}\gamma^2r_{\max}(1+\delta_{\max})$ which converges at most sublinearly at a floor, where $\Delta_k\coloneqq f(X_k) - \fstar$ is the optimality gap. 
            On the other hand, for a general $L$-Lipschitz smooth but possibly nonconvex objective function $f\colon \RR^{m \times n} \to \oRR$, the iterates of matrix sign descent $(X_k)_{k\in\setK}$ satisfy $\min_{k\in\setK} \fronorm{\nabla f(X_k)} \le \scrO\left(\frac1{\gamma (1-\varepsilon_{\max})K} + \frac{L}{2}\gamma r_{\max}\frac{1+\delta_{\max}}{1-\varepsilon_{\max}}\right)$, and the iterates of matrix \signSGD $X_{k+1} = X_k - \gamma\tU_k$ with $\tU_k\tH_k = \polarhat(\nabla f(X_k, \xi_k))$ satisfy 
            \[
            \min_{k\in\setK}\Ex\fronorm{\nabla f(X_k)} \le \scrO\left(\frac1{\gamma (1-\hat\varepsilon_{\max})K} + \frac{L}2 \gamma \hat{r}_{\max}\frac{1+\hat\delta_{\max}}{1-\hat\varepsilon_{\max}} + \varsigma\frac{\sqrt{\hat{r}_{\max}(1+\hat\delta_{\max})}}{1-\hat\varepsilon_{\max}}\right)
            \]
            if \Cref{assum:grad_noise} also holds. 
        \end{theorem}
        
        We are interested in what the above convergence rates would become when specific numerical polar decomposition algorithms are used in practice. Let us denote the number of inner steps used within any inexact polar oracles by $T\in\NN^*$. We now provide results specific to inexact polar oracles used in practice including the NS iteration and the QDWH algorithm \citep{nakatsukasa2010optimizing}, by determining the orders of $\varepsilon_{\max}$ and $\delta_{\max}$ (resp.~$\hat\varepsilon_{\max}$ and $\hat\delta_{\max}$) in terms of $T$. From this we are also able to determine the order of the number of inner steps $T$ required for different numerical polar decomposition algorithms given a desired level of accuracy. For simplicity, we only detail the results for the deterministic case under \PL condition (\Cref{assum:base}). The stochastic and nonconvex cases are more involved but can be obtained by plugging in the corresponding values of $\hat\varepsilon_{\max}$ and $\hat\delta_{\max}$ in \Cref{thm:polarsgd_strcvx_inexact,thm:mat_signSGD_strcvx_inexact}.

        \begin{theorem}[Newton--Schulz]\label{thm:NS}
            Running the Newton--Schulz iteration with quintic polynomials $\tU_{k,j+1} = a\tU_{k,j} + b\tU_{k,j}\tU_{k,j}^\top \tU_{k,j} + c\tU_{k,j}(\tU_{k,j}^\top \tU_{k,j})^2$ and $\tU_{k,0} = G_k/\fronorm{G_k}$ with coefficients $(a,b,c) = (15/8, -5/4, 3/8)$ for $T$ inner steps so that $\tU_k = \tU_{k,T}$, we have the oracle error bounds $\varepsilon_{\max}(T) = \scrO(e_0^{3^T})$ and $\delta_{\max}(T) = \scrO(e_0^{3^T})$, where $e_{k,j} \coloneqq \specnorm{\tU_{k,j}^\top \tU_{k,j} - I}$ for $k\in\{0,\ldots,K\}$ and $j\in\{0,\ldots,T\}$, and $e_0 \coloneqq \max_{k\in\{0,\ldots,K\}} e_{k,0}$. 
            Therefore, when running realized \PolarGrad \eqref{eqn:polargrad_practical} with the quintic polynomial Newton--Schulz iteration under \Cref{assum:base,assum:inexact_polar_oracle}, to stay within $1-\eta$ of the exact rate in \Cref{thm:polargrad_strcvx} for some $\eta\in(0, 1)$, it requires at least $\left\lceil \scrO(\log(\log\eta / \log e_0) \right\rceil$ inner steps. 
        \end{theorem}
        The above theorem says that both oracle errors decay triply exponentially and helps us determine the number of inner steps required if we specify $\eta$ and have the knowledge of $e_0$ (usually through initialization). 
        Since the \textsc{Polar Express} \citep{amsel2025polar} is an improved variant of the NS iteration with quintic polynomials with dynamically optimized polynomial coefficients, the above corollary also applies to the \textsc{Polar Express} with potentially better error bound constants (cf.~Theorem 4.3 of \citep{amsel2025polar}).     
        \begin{remark}
        	The default coefficients $(a,b,c)=(3.4445, -4.775, 2.0315)$ in the quintic matrix iterative polynomial in \Muon\footnote{Also the default coefficients in PyTorch's \texttt{torch.optim.Muon} \citep{pytorch_muon2025} and Optax's \texttt{optax.contrib.muon} \citep{optax_muon2025}. } \citep{jordan2024muon} do not lead to \emph{convergent} polar decomposition \citep{amsel2025polar}, especially for ill-conditioned matrices. The coefficients chosen in \Cref{thm:NS} are determined by solving the conditions $\varphi(1)=1$, $\varphi'(1) =0$ and $\varphi''(1)=0$ of the quintic polynomial $\varphi(t) = t(a+bt+ct^2)^2$. 
        \end{remark}
    
        The QDWH algorithm \citep{nakatsukasa2010optimizing} also has oracle errors decay triply exponentially, so it has similar oracle error bounds to those of the NS iteration. 
            
        \begin{theorem}[QDWH]\label{thm:QDWH}
            Running the QDWH algorithm or its equivalent DWH iteration $\tU_{k, j+1} = \tU_{k,j}(a_j I + b_j \tU_{k,j}^\top \tU_{k,j})(I + c_j\tU_{k,j}^\top \tU_{k,j})^{-1}$ and $\tU_{k,0}= G_k/\specnorm{G_k}$ with dynamic weighting parameters $(a_j, b_j, c_j)$ for $T$ inner steps, we have the error bounds $\varepsilon_{\max}(T)=\scrO((1-\ell_0^2)^{3^T})$ and $\delta_{\max}(T)=\scrO((1-\ell_0^2)^{3^T})$ where $\ell_0$ is the smallest lower bound on the singular value of $\tU_{k,0}$ over all iterations $k\in\setK$. Therefore, when running realized \PolarGrad \eqref{eqn:polargrad_practical} with the QDWH algorithm under \Cref{assum:base,assum:inexact_polar_oracle}, to stay within $1-\eta$ of the exact rate in \Cref{thm:polargrad_strcvx} for some $\eta\in(0, 1)$, it requires a number of $T \ge \left\lceil \scrO(\log(\log\eta / \log (1-\ell_0^2)) ) \right\rceil$ inner steps. 
        \end{theorem}

        \begin{remark}[Inexactness in numerical polar decomposition algorithms]
            The recent work \citep{shulgin2025beyond} is the first known work which studies the inexact orthogonalized update for \Muon by introducing a realistic additive error model within the general framework of LMO-based optimization. Our analysis has two major noted differences from theirs.   
            In Theorem 3 of \citep{shulgin2025beyond} where an adaptive learning rate involving the dual norm of the gradient is considered, the inexactness arising from the computation of the dual norm is not accounted for in their analysis, implicitly assuming that its computation is readily available. For the case of the spectral norm whose dual norm is the nuclear norm, the computation of the nuclear norm is known to be essentially as expensive as that of the full SVD. Omitting the inexactness of the computation of the nuclear norm in practice overlooks another source of inexactness in the realized algorithms. In contrast, we compute the nuclear norm of the (deterministic or stochastic) gradient $G_k$ using the approximation $\dotpF{\tU_k}{G_k}$ where $\tU_k$ is obtained from an inexact polar oracle, and include this source of inexactness in our analysis. Furthermore, we also provide specific results for inexact polar oracles, namely the NS iteration and the QDWH algorithm. We however do not provide results for \PolarGradM and leave it for future work. 
        \end{remark}

    	\begin{remark}[Comparing Newton--Schulz and QDWH]
    		\label{remark:NS_QDWH}
    		While \Cref{thm:NS,thm:QDWH} inform us of similar oracle error bounds, their error constants are indeed vastly different since the NS iteration is a \emph{polynomial} iteration whereas QDWH is a \emph{rational} iteration. To see this, we recall that both the NS iteration and QDWH give cubic convergence of orthogonality error $e_{j+1} \le \zeta e_j^3$ for some $\zeta>0$, where $e_j\coloneqq \specnorm{\tU_j^\top\tU_j - I}$ for $j\in\setT$. Since the NS iteration is a polynomial iteration, its local error constant $\zeta_{\mathrm{NS}}$ depends strongly on $e_0 = 1 - \ell^2$, where $\ell = \sigma_{\min}(G)/\sigma_{\max}(G) $. The initial error is close to $1$ when $\ell$ is small, so the iteration enters its cubic regime later or \emph{never}. Therefore, the NS iteration loses its cubic convergence behavior and may even diverge without additional rescaling if $G$ is so ill-conditioned that its local error constant $\zeta_{\mathrm{NS}}$ could be unbounded. On the other hand, QDWH's local error constant $\zeta_{\mathrm{QDWH}}$ is bounded and does not blow up as $\ell\to0$ because its rational part $(I + c_j M_j)^{-1}$ compresses large singular values and stretches small ones, and keeps the iteration centered at the optimal cubic fixed point. QDWH is indeed \emph{provably stable} and \emph{cubically convergent} even when $\kappa_2(G)= 10^{16}$ \citep{nakatsukasa2010optimizing}. 
    	\end{remark}

        \begin{remark}[Choice of polar oracles in \PolarGrad]
            From the above comparison of convergence results, we can roughly determine the number of inner steps of each of the considered polar oracles for a desired level of accuracy. However, for practical usage, when choosing a suitable polar oracle, there are various factors for consideration, such as computational cost, required precision, numerical stability, hardware consideration such as GPU-friendliness of involved operations (e.g., matrix multiplications, scalar multiplications and their linear combinations), as well as the complexity of the operations involved. While the NS iteration and the \textsc{Polar Express} would be better suited for deep learning due to their lower FLOPS and GPU-friendliness, the QDWH algorithm could be more desirable for ill-conditioned gradient/momentum matrices and when solving smaller-scale matrix optimization problems on CPUs and higher precision is desired. 
        \end{remark}

        \begin{remark}[Optimizers for embedding and head layers]
        	\label{remark:optim_embed}
        	While the input embedding and head layers also have matrix parameters, the current training protocols of \Muon still use \textsc{Adam}(W) for these two layers \citep{jordan2024muon}. There is indeed a mismatch between this practical choice of optimizers and the corresponding choice of norms for steepest descent as suggested in Example 6 of \citep{bernstein2024modular}. Here we provide a principled explanation based on the choice of numerical polar decomposition algorithms and the corresponding choice of optimizers. Let us consider an input embedding matrix $E\in\RR^{V\times d}$ and the head matrix $W\in\RR^{V\times d}$ where $V$ is the vocabulary size and $d$ is the embedding dimension with $V\gg d$. For the input embedding, its gradient has the form $G_E = S^\top H$, where $S\in\RR^{b\times V}$ is a sparse token-selection or count matrix (one-hot),  $H\in\RR^{b\times d}$ is a dense backpropagated signal and $b$ is the batch size. Consequently, the gradient is rank-deficient since $\rank(G) \le \min\{b, d\}\ll d$ and fluctuates with batch composition. For very large vocabulary size $V$, many rows are never ``touched'' in a batch so the lower bound $\ell\coloneqq \sigma_{\min}(G_E)/\sigma_{\max}(G_E) \approx 0$. In the case of stochastic gradient, the small singular values are thus dominated by stochastic noise, not signal. Thus, for the input embedding, polynomial polar oracles such as the NS iteration or the \textsc{Polar Express} all have an initial orthogonality defect of $e_0=1-\ell^2\approx1$, so \Muon or \PolarGrad updates based on these polar oracles become weak, noisy or unstable. As discussed in \Cref{remark:NS_QDWH}, \Muon or \PolarGrad updates based on rational approximations such as QDWH is still stable and convergent. The head layer is even worse than token embeddings since its gradient $G_W$ is driven by softmax logits with highly skewed distributions where rare tokens get near-zero signal, leading to an even more ill-conditioned spectrum. In short, input embedding and head layers operate in an extreme ill-conditioning regime where polynomial polar oracles lose their theoretical guarantees. In contrast, rational approximation methods such as QDWH are therefore structurally better suited for these layers. \PolarGrad with cheap polar oracles is most effective on well-conditioned blocks such as attention and linear layers. We also empirically demonstrate in \Cref{subsec:qwen} that QDWH-\PolarGrad optimizers can still be used for these two layers, instead of \Muon or NS-\PolarGrad.     	
        	While QDWH works well for layers with ill-conditioned gradients, it does come with a cost of expensive QR decomposition and is especially problematic for huge $V\times d$ matrices. Through the same lens, even though \Adam does not compute a polar direction, it implicitly applies a \emph{diagonal rational preconditioner} whose directions with tiny singular values are heavily damped and suppresses small-singular-value noise when viewed spectrally. However, the diagonal structure does not capture correlations across the $d$-dimensional embedding space and completely ignores the matrix geometry. It can also have very different implicit bias and scaling behavior from polar or spectral gradient methods. Consequently, QDWH-\PolarGrad could be more desired if the embedding dimension $d$ is small or moderate, or QDWH is performed infrequently and cheaper updates are kept in between. 
        	
        \end{remark}
        
        \section{A Unifying Preconditioning View of Adaptive Gradient Optimizers}    
        \label{sec:precond}
        Adaptive gradient optimizers are a family of stochastic gradient methods which are usually understood to accelerate convergence by employing adaptive learning rates. In this section, we use $x\in\RR^d$ or $X\in\RR^{m\times n}$ to denote the optimization variable. The stochastic gradient is denoted by $g_k = \nabla f(x_k, \xi_k)$ with the sample $\xi_k$. Most adaptive gradient optimizers can be written as 
        \begin{equation}\label{eqn:ago}
        	(\forall k\in\NN)\quad x_{k+1} = x_k - \gamma_k\cdot \sfm_{k-1}(g_k)/\sfv_{k-1}(g_k^2))^{\half}, 
        \end{equation}
        where $\sfm_{k-1}\colon\RR^d\to\RR^d$ and $\sfv_{k-1}\colon\RR^d\to\RR^d$ are functions of the gradient and the coordinate-wise squared gradient conditioned on the past iterates and gradients $\{x_0, g_0, \ldots, x_{k-1}, g_{k-1}\}$, respectively. Here the division and addition operations are performed coordinatewise. The quantity $\gamma_k/(\sfv_{k-1}(g_k^2))^{\half}$ can be viewed as an adaptive learning rate of the adaptive gradient optimizer.     
        This subsumes adaptive gradient optimizers commonly used in deep learning including \AdaGrad \citep{duchi2011adagrad,mcmahan2010adaptive}, \Adadelta \citep{zeiler2012adadelta}, \RMSprop \citep{tieleman2012} and \Adam \citep{kingma2015}, as well as their many variants. 
       
        While this adaptive learning rate view has been widely accepted by the deep learning community for the success of adaptive gradient optimizers, its intrinsic motivation is indeed to approximate quasi-Newton (or second-order) methods or (inverse) Hessian approximation \citep{yao2021adahessian,su2024hessian}. However, there is still a gap in understanding whether approximate second-order methods can still accelerate convergence for highly nonconvex problems such as neural network training. 
        We emphasize that they can and should be viewed as preconditioned gradient methods (see e.g., Chapter 5 of \citep{bach2024learning}). 
        To better understand such issues, we provide a more detailed exposition of these views below. 
        
        \subsection{Three Views of Adaptive Gradient Optimizers}
        \label{subsec:adaptive_grad}
        
        \paragraph{Adaptive learning rate.}
        Using the general formulation \eqref{eqn:ago}, a coordinatewise adaptive learning rate in the form of $\gamma_k/(\sfv_{k-1}(g_k^2))^{\half}$ is generally used. For instance, in \AdaGrad, the adaptive learning rate is given by $\gamma_k/(\sum_{t=1}^k g_t^2 + \varepsilon)^{\half}$, where $\varepsilon>0$ is a small constant for ensuring numerical stability. The main advantage of adaptive learning rates is that they allow different magnitudes of updates in different coordinates.  
        
        \paragraph{Diagonal inverse Hessian approximation.}
        Motivated by quasi-Newton methods such as BFGS, adaptive gradient optimizers can also be viewed as approximating the inverse (square root) of the Hessian $H_k\coloneqq\nabla^2 f(x_k, \xi_k)$. To see this, let us denote a stochastic gradient by $g_k \coloneqq \nabla f(x_k, \xi_k)$. Then the Gauss--Newton method approximates the Hessian by $H_k \approx g_kg_k^\top$ (dropping a factor of $2$ for more coherent discussion). To save memory, it is further approximated only by its diagonal, which is thus given by $\Diag(g_k^2)$. To ensure that this Hessian approximation is invertible when $f$ is nonconvex, a constant diagonal matrix is added to it, i.e., $\Diag(g_k^2+ \varepsilon)$, where $\varepsilon>0$ is a small positive constant.  Since it is a diagonal matrix, its inverse is simply $\Diag(1/(g_k^2+\varepsilon))$, where the division and addition operations are performed coordinatewise. In most adaptive gradient optimizers such as \Adam and \RMSprop, the exponential moving average of the squared historical gradients with a coordinatewise square root are used instead. While we can apply this directly for matrix optimization problems by vectorizing all matrices and performing coordinatewise updates as in \Adam, there remains a large gap in justifying that this is still technically correct as a diagonal inverse square root Hessian approximation for matrices since we want to maintain their original matrix structures. 
        
        \paragraph{Preconditioning and preconditioned gradient methods.}
        In addition to the above two views, we emphasize the importance of employing a preconditioning view. Borrowing from the details of the above Hessian approximation view, the preconditioner of adaptive gradient optimizers can be further generalized as the diagonal approximation of the inverse Hessian with the exponential moving average of the squared gradients, given by $\Diag\left(1/(\sfv_{k-1}(g_k^2) + \varepsilon)^{\half}\right)$. We now turn to the inner workings of preconditioning in preconditioned gradient methods. In general, preconditioning via the inverse Hessian or its approximation achieves accelerated convergence by minimizing the condition number of the objective function (see e.g., Chapter 5.2 of \citep{bach2024learning}). 
        There are however two separate notions of condition numbers arising in matrix analysis and optimization theory, one being the condition number of a matrix defined through the ratio of its largest and smallest positive singular values, while the other is the condition number of an optimization problem given by the ratio of the Lipschitz smoothness constant and the strong convexity constant of the objective function. We will draw the connection between these two notions of condition numbers below. 
        
        \begin{remark}
           	All the above three views consider optimization variables as vectors. When adaptive gradient optimizers are applied to matrix parameters in neural networks, as all operations in adaptive gradient optimizers are coordinatewise, it is equivalent to applying these optimizers to the vectorized (i.e., flattened) matrix parameters. We emphasize that the treatment of preconditioning for matrix-valued updates is very different from that for their vectorized counterparts. 
        \end{remark}

        \subsection{Vector Preconditioned Gradient Methods}
        We first consider the vector optimization problem $\minimize_{x\in\RR^d} f(x)$ with the objective function $f\colon\RR^d\to\oRR$ with $d\in\NN^*$. The vector preconditioned gradient method can be written as 
        \begin{equation}\label{eqn:pgd}
        	(\forall k\in\NN)\quad x_{k+1} = \argmin_{x\in\RR^d} \,\left\{ \dotp{g_k}{x - x_k} + \frac{1}{2\gamma_k}\norm{x - x_k}_{P_k^{-1}}^2 \right\} = x_k - \gamma_k P_k g_k, 
        \end{equation}
        where $g_k=\nabla f(x_k)$ and $P_k\in\bbS_{++}^d$ is a \emph{preconditioning matrix} or \emph{preconditioner}. Let us suppose that $f$ is $L$-Lipschitz smooth and $\mu$-strongly convex, i.e., $\mu\euclidnorm{x-y}\le \euclidnorm{\nabla f(x) - \nabla f(y)} \le L\euclidnorm{x-y}$ for any $(x,y)\in\RR^d\times\RR^d$, where $0<\mu\le L<\infty$. 
        More specifically, we explicitly denote these constants for the objective function $f$, i.e., $L=L_\vec(f)$ and $\mu=\mu_\vec(f)$. 
        Assuming that $f$ is twice continuously differentiable, then there is an intimate relationship between these constants and the spectrum of the Hessian of $f$, given by
        $L_\vec(f) = \sigma_{\max}(\nabla^2 f)$ and $\mu_\vec(f) = \sigma_{\min}(\nabla^2 f)$. Then, the condition number of the objective $f$ can be defined as $\kappa_\vec(f) \coloneqq L_\vec(f) / \mu_\vec(f)=\kappa_2(\nabla^2 f)$.     
        For most loss functions in deep learning, the constants $L_\vec$ and $\mu_\vec$ are global constants that are expensive to evaluate in general or do not exist. We can however define their corresponding local versions (at each iterate). 
        The local condition number of $f$ at $x\in\dom f$ can be defined by $\kappa_\vec(f)(x) \coloneqq L_\vec(f)(x) / \mu_\vec(f)(x)=\kappa_2(\nabla^2 f(x))$.     
        As a result, at each iteration $k\in\NN^*$, the equality $\kappa_\vec(f)(x) =\kappa_2(\nabla^2 f(x))$ imply that the inverse Hessian $P_k = (\nabla^2 f(x_k))^{-1}\in\bbS_{++}^d$ is the best local preconditioner as $\kappa_2(\nabla^2 f(x)^{-1})=\kappa_2(\nabla^2 f(x))^{-1}$, also explaining the fast convergence of Newton's method for strongly convex and Lipschitz Hessian objectives. 
    
        \paragraph{Adaptive gradient optimizers as vector preconditioned gradient methods.}
        In general, the objective function $f$ is nonconvex in deep learning. Adaptive gradient optimizers thus attempt to approximate preconditioners that are positive definite. For memory and computational efficiency, a diagonal preconditioner $P_k = \Diag(p_k)\in\bbS_{++}^d$ with positive diagonal entries $p_k\in\Rpp^d$ is often used, rather than the full inverse Hessian matrix. For instance, in \RMSprop and \Adam, the diagonal preconditoner is given by $p_k = 1/(\hat{v}_k^{\odot\half}+\varepsilon)$ with $\hat{v}_k = (1-\beta_2)\sum_{t=0}^k\beta_2^{k-t}g_t^2/(1-\beta_2^{k+1})$. 
        
        \paragraph{Issues with diagonal approximations of inverse Hessian.}
        While diagonal approximations of explicit preconditioners are more memory- and compute-efficient than the inverse Hessian, it might lead to declined preconditioning effect or could even be detrimental even for simple nonconvex objectives, potentially leading to divergence of such diagonally preconditioned gradient methods. See \Cref{subsec:mat_fac} for example. We could potentially attribute the training instabilities of LLMs using \textsc{Adam(W)} to this diagonal approximation.

        \subsection{Matrix Preconditioned Gradient Methods}
        We now consider the matrix optimization problem $\minimize_{X\in\RR^{m\times n}} \sff(X)$ with the objective function\footnote{Note that we use $f$ and $\sff$ to denote the vector and matrix optimization problem objectives respectively in this section. } $\sff\colon\RR^{m\times n}\to\oRR$, where $m$ and $n$ are positive integers both strictly greater than one\footnote{Here we omit the cases of $X$ being reduced to a vector or a scalar. }. A general matrix preconditioned gradient method can be written as $X_{k+1} = X_k - \gamma_k \scrP_k(G_k)$, where $G_k=\nabla\sff(X_k)$ is the gradient\footnote{In the case of fitting neural networks, $G_k$ represents the partial derivative of the loss function with respect to the matrix-valued parameter $X$ of a single layer at $X_k$, not the parameters of all layers. } of $\sff$ with respect to $X$ at $X_k$ and $\scrP_k\colon\RR^{m\times n}\to\RR^{m\times n}$ is a \emph{preconditioning function}. Such a preconditioning function can be very general. 
        The local condition number of the objective $\sff$ at $X\in\dom\sff$ can be defined by $\kappa_\mat(\sff)(X) \coloneqq L_\mat(\sff)(X) / \mu_\mat(\sff)(X)$. 
        If $\sff$ is twice continuously differentiable, then we also have $L_\mat(\sff)(X) = \sigma_{\max}(\nabla^2\sff(X))$ and $\mu_\mat(\sff)(X) = \sigma_{\min}(\nabla^2\sff(X))$, with the Hessian $\nabla^2\sff(X)\in\RR^{mn\times mn}$. These notions are indeed defined equivalently to their vector counterparts through vectorization. While most existing vector preconditioned gradient methods are curvature-aware and aim to reduce the (local) condition number of the Hessian, it is generally very hard to compute or approximate the Hessian w.r.t.~matrix parameters without assuming specific structures such as Kronecker-factored in K-FAC \citep{martens2015optimizing}. However, the matrix structure of the optimization variable $X$ and its gradient has led us to introduce another preconditioning concept for matrix optimization problems called \emph{gradient-anisotropy preconditioning}, which instead minimizes the condition number of the matrix-valued gradient. Before detailing this concept, we first introduce the \Muon optimizer \citep{jordan2024muon,su2024muon} which indeed performs this kind of preconditioning.

        \subsection{Curvature-Anisotropy Preconditioning vs.~Gradient-Anisotropy Preconditioning}
        \label{subsec:precond}
        While the interpretation of \Muon as stochastic steepest descent w.r.t.~the spectral norm or matrix sign descent can be derived directly by solving the subproblems as in \eqref{eqn:ssd}, we gain further insights into the orthogonalization step in \Muon. In what follows, we advocate for a gradient preconditioning view due to (semi-)orthogonal projections of the gradient or momentum. 
        
        While almost all existing preconditioned gradient methods in the literature consider preconditioning that addresses curvature anisotropy through reducing the Hessian condition number, more recently the class of orthogonalized gradient methods such as \Muon \citep{jordan2024muon} indeed perform preconditioning that addresses gradient or momentum anisotropy. 
        Gradient anisotropy implies discrepancy of the strength of the gradient magnitude in different directions, captured by the ($2$-)condition number of the gradient matrix, $\kappa_G(X) \coloneqq\kappa_2(\nabla\sff(X))$. 
        In contrast, curvature anisotropy is captured by the condition number of the Hessian, $\kappa_H\coloneqq\kappa_2(\nabla^2\sff)$. 
        The Hessian condition number $\kappa_H$ informs us of how distorted gradient directions are globally, whereas the gradient condition number $\kappa_G$ concerns local distortion of gradient directions at each iteration. We will see how these quantities govern the convergence rates of related algorithms in \Cref{subsec:conv}. 
        
        To address gradient anisotropy, (semi-)orthogonal projections of the gradients are usually performed, since ``the best conditioned matrices are the orthogonal ones, which have condition numbers of 1'' \citep{turing1948rounding}. They allow capturing only the gradient directions but ignore its magnitude---effectively removing all curvature information. Contrarily, the full inverse Hessian preconditioner corrects anisotropy proportionally and adjusts the gradient directions for local geometry with curvature-awareness, which can be more beneficial. However, in large-scale applications and stochastic nonconvex problems such as neural network training, the former method is more stable and easier to implement without the need of approximation, and has a lower computational cost. 
        From this angle, we can establish an interpretation that \Adam and most other adaptive gradient optimizers are vector curvature-anisotropy preconditioned gradient methods, whereas \Muon, \Shampoo and their variants are matrix gradient-anisotropy preconditioned gradient methods. Dropping all curvature information with isotropic updates could however be detrimental to the optimization process; see \Cref{subsec:nuc_norm} for a related discussion on its mitigation. To better understand the similarities and differences between these two approaches to preconditioning, we give the following simple matrix quadratic regression example. 
        
        \begin{example}[Matrix quadratic regression]\label{example:quad}
        	Let us consider a matrix quadratic regression objective $\sff(X) \coloneqq \frac12\fronorm{AXB-C}^2$, where $X\in\RR^{m\times n}$, $A\in\RR^{p\times m}$, $B\in\RR^{n\times q}$ and $C\in\RR^{p\times q}$. Then its gradient is $\nabla\sff(X) = A^\top(AXB-C)B^\top$, its Hessian is $\nabla^2\sff(X) = (BB^\top)\otimes(A^\top A)\in\RR^{mn\times mn}$, and the inverse-Hessian preconditioned gradient is given by $G_{\mathsf{pre}}(X)\coloneqq (A^\top A)^{-1}\nabla\sff(X)(BB^\top)^{-1}$. If we define $E\coloneqq AXB-C$, then the gradient condition number is given by $\kappa_2(\nabla\sff(X)) = \kappa_2(A^\top(AXB-C)B^\top) \le \kappa_2(A)\cdot\kappa_2(B)\cdot\kappa_2(E)$. 
        	The Hessian condition number is given by $\kappa_2(\nabla^2\sff(X)) = \kappa_2(A)^2\cdot\kappa_2(B)^2$, while the condition number of the preconditioned gradient is given by $\kappa_2(G_{\mathsf{pre}}(X)) = \kappa_2(A^\dagger E (B^\dagger)^\top)\le \kappa_2(A)\cdot\kappa_2(B)\cdot\kappa_2(E)$, where $A^\dagger \coloneqq (A^\top A)^{-1}A^\top$, $B^\dagger\coloneqq (BB^\top)^{-1}B$. Thus, we can use $\kappa_2(E)$ to understand the convergence behavior of different optimizers since $\kappa_2(A)$ and $\kappa_2(B)$ are constant. The preconditioned gradient using a (semi-)orthogonal projection always has a condition number $\kappa_2(\msgn(\nabla\sff(X)))=1$, hence discarding all curvature information brought by the residual $E$. Numerical studies can be found in \Cref{subsec:mat_quad_reg}. 
        \end{example}

        \subsection{Explicit Preconditioners vs.~Implicit Preconditioners}
        \label{subsec:explicit_implicit}
        Adopting such a unifying preconditioning view, most popular deep learning optimizers can also be categorized into those with explicit and implicit preconditioners respectively. Implicit preconditioners are often derived from steepest descent w.r.t.~non-Euclidean norms, while explicit preconditioners are often derived from steepest descent w.r.t.~preconditioned Euclidean norms as in \eqref{eqn:pgd} or Kronecker-factored preconditioners in K-FAC \citep{martens2015optimizing}. Detailed exposition can be found in \Cref{subsec:explicit_implicit_2}. 
            
        \paragraph{Vector preconditioned gradient methods.}
        While most vector preconditioned gradient methods such as \Adam and \RMSprop have explicit preconditioners $P_k$, preconditioning can also be performed using implicit preconditioners (in the form of a preconditioning function). A notable example is (unnormalized) \signSGD \citep{bernstein2018signsgd,balles2018dissecting} (see also \citep{xie2024adam,xie2024implicit}): 
        \begin{equation}\label{eqn:signSGD}
        	(\forall k\in\NN)\quad x_{k+1} = \argmin_{x\in\RR^d} \,\left\{ \dotp{g_k}{x - x_k} + \frac{1}{2\gamma_k}\infnorm{x - x_k}^2 \right\} = x_k - \gamma_k \onenorm{g_k}\cdot\sgn(g_k). 
        \end{equation}
        Let us also recall that \Adam \citep{kingma2015} with $\beta_1=\beta_2=0$ recovers \signSGD so \Adam can be viewed as a form of smoothed \signSGD with an explicit preconditioner. \signSGD thus has the elementwise sign function $\sgn$ as an implicit preconditioner. 
        
        \paragraph{Matrix preconditioned gradient methods.}    
        An analogous viewpoint also holds for matrix preconditioned gradient methods. In particular, due to the connection between \Muon and \Shampoo \citep{gupta2018shampoo} given in \citep{bernstein2024old,jordan2024muon}, \Shampoo can be viewed as a matrix preconditioned gradient method with explicit left and right preconditioners $L_k\in\RR^{m\times m}$ and $R_k\in\RR^{n\times n}$, whose update rules are given by 
        \[
        L_k = \beta L_{k-1} + (1-\beta)G_k G_k^\top, \; R_k = \beta R_{k-1} + (1-\beta)G_k^\top G_k, \; X_{k+1} = X_k - \gamma_k L_k^{-\sfrac14} G_k R_k^{-\sfrac14}.  
        \]

        \subsection{Vector Preconditioned Gradient Methods vs.~Matrix Preconditioned Gradient Methods}
        \label{subsec:vec_mat_precond}
        Let us recall the equivalence between preconditioned gradient methods with implicit preconditioners and steepest descent methods w.r.t.~non-Euclidean norms for both vector and matrix optimization problems. Indeed, leveraging this preconditioning perspective, we are able to develop an explanation of the potential inappropriateness of adaptive gradient optimizers like \Adam for matrix parameters in neural networks.     
        Again, considering \signSGD as a particular instance of \Adam with $\beta_1=\beta_2=0$, the update \eqref{eqn:signSGD} for matrices becomes 
        \begin{equation}
        	(\forall k\in\NN) \quad X_{k+1} \in \argmin_{X\in\RR^{m\times n}}\,\left\{\dotpF{G_k}{X - X_k} + \frac{1}{2\gamma_k}\matsnorm{X - X_k}{\max}^2 \right\}, 
        \end{equation}
        where $G_k=\nabla\sff(X_k)$ and $\matsnorm{X}{\max} \coloneqq \max_{1\le i\le m, 1\le j\le n} |x_{i,j}| = \infnorm{\vec(X)}$ is the max norm of $X\in\RR^{m\times n}$, where $m$ and $n$ are positive integers both strictly greater than one.     
        Unlike the spectral norm, the max norm $\matsnorm{\cdot}{\max}$ is neither a matrix norm (see Chapter 5.7, Example 5 of \citep{horn2012matrix}) nor a unitarily invariant norm \citep{mirsky1960symmetric}. Comparing the elementwise sign function and the matrix sign function imposed on matrix gradients, the former only takes the sign of each entry whereas the latter sets all singular values to one while maintaining the directions of the original gradients characterized by singular vectors. The preconditioning effect of the elementwise sign function on matrix gradients is inconclusive, which might even change the singular vectors and thus the original update direction provided by the gradient, or even worsen the gradient and/or Hessian condition numbers. This might potentially lead to pre-training instabilities of language models using \AdamW. 
        
        \subsubsection{signSGD on Matrices is SSD on The Diagonal Matrization of Its Vectorization}
        \label{subsubsec:elemenwise_matrix}
        Motivated by the recent \MuonAll optimizer \citep{page2025muonall}, we show that we can indeed recover unnormalized \signSGD from stochastic spectral descent (SSD) by embedding a vector variable as a diagonal matrix, drawing another connection between these two classes of optimizers. 
        
        We now consider the vector variable $x\in\RR^d$ and ``matrize'' it as the diagonal matrix $D \coloneqq \Diag(x)\in\RR^{d\times d}$. Now we define $F\colon\RR^{d\times d}\to\oRR$ such that $F(D) = f(\diag(D)) = f(x)$, where $\diag$ is the adjoint of $\Diag$ which extracts the diagonal of a matrix into a vector. Since the map $x\mapsto\Diag(x)$ is linear with adjoint $\diag$, we have $G\coloneqq\nabla F(D) = \Diag(\nabla f(x)) = \Diag(g)$, where $g\coloneqq\nabla f(x)$. Then, with $G\coloneqq\Diag(g)$, the orthogonal polar factor of $G$ is equal to $G(G^\top G)^{-\half} = (\Diag(g_i/|g_i|))_{1\le i\le d} = \Diag(\sgn(g))$. Moreover, the nuclear norm of $G$ also reduces to the $\ell_1$-norm of $g$, i.e., $\nucnorm{G} = \sumd |g_i| = \onenorm{g}$. Hence, running SSD on $D$ takes the form $D_{k+1} = D_k - \gamma_k\onenorm{g_k}\Diag(\sgn(g_k))$, which essentially amounts to running unnormalized \signSGD in its vector form $x_{k+1} = x_k - \gamma_k\onenorm{g_k}\sgn(g_k)$. Similar arguments hold when momentum is also considered, recovering \Signum from \Muon for instance.     
        
        Consequently, two further conclusions can be drawn here: (i) The \MuonAll optimizer is indeed equivalent to \Signum for vector or scalar parameters and \Muon for matrix parameters; (ii) Running unnormalized \signSGD (an instance of \Adam) on a matrix parameter $X\in\RR^{m\times n}$ elementwise is equivalent to running SSD on the diagonal matrization of its vectorization $\Diag(\vec(X))\in\RR^{mn\times mn}$ and then flattening it back to $\RR^{m\times n}$. It is not hard to see that the gradients $\nabla\sff(X)$ and $\Diag(\vec(\nabla\sff(X)))$ have different spectral properties including polar decomposition (see \Cref{def:polar_decomp}).

        \subsubsection{Reduction of Matrices to Vectors or Scalars in SSD and \Muon}
        \label{subsubsec:m_n_1}
        In the above discussion, we deliberately exclude the corner case of the matrix variable $X$ being reduced to a vector or a scalar, i.e., the cases of $m=1$ or $n=1$ and $m=n=1$. This is consistent with the practical use of \Muon where elementwise optimizers such as \textsc{Adam}(W) (or \Lion in \citep{ahn2025dion}) is used for vector and scalar parameters in a neural network. 
        
        We now give a potential explanation for this choice. When $X$ is a (row or column) vector ($m=1$ or $n=1$ but not both $m$ and $n$ are one), SSD reduces to vanilla \SGD whereas \Muon without momentum reduces to $\ell_2$-normalized \SGD. On the other hand, when $X$ is a scalar ($m=n=1$), SSD again reduces to vanilla \SGD whereas \Muon without momentum reduces to \signSGD. 
        To see this, without loss of generality, we consider the case where the iterate $x_k\in\RR^{1\times n}$ is a column vector with $n\ge1$. Then the gradient $g_k$ is a nonzero rank one matrix with the SVD $g_k = \sigma_k u_k v_k^\top$, where $\sigma_k = \euclidnorm{g_k}$, $u_k = g_k / \euclidnorm{g_k}$ and $v_k = 1$. Since $\rank(g_k) = 1$, we have $\nucnorm{g_k} = \sigma_k = \euclidnorm{g_k}$. Hence, SSD is equivalent to \SGD: $x_{k+1} = x_k - \gamma_k \euclidnorm{g_k}\cdot g_k / \euclidnorm{g_k} = x_k - \gamma_k g_k$, while \Muon without momentum takes the form of $\ell_2$-normalized \SGD: $x_{k+1} = x_k - \gamma_k g_k / \euclidnorm{g_k}$. Now, if we further set $n=1$, then we have $\sigma_k = \euclidnorm{g_k} = |g_k|$ and $u_k = g_k / |g_k| = \sgn(g_k)$, so that $\ell_2$-normalized \SGD reduces to \signSGD. Similar arguments remain valid when momentum is used. 
        
        As a result, we can see that SSD and \Muon both reduce to vanilla stochastic gradient methods without preconditioning when the parameter is a vector. This suggests that vector preconditioned gradient methods like \textsc{Adam}(W) or \Lion are more favored for vector parameters to accelerate convergence.

        \section{Proofs}
        \label{sec:proofs}
        We provide proofs of the results in \Cref{sec:polar_grad} in this section. 
        
        \subsection{Proofs for \Cref{subsec:conv}}
        
        \begin{proof}[Proof of \Cref{thm:polargrad_strcvx}]
           	By the $L$-Lipschitz smoothness of $f$ (\Cref{def:Lipschitz_smooth}), we have         
           	\begin{align}
           		f(X_{k+1}) &\le f(X_k) + \dotpF{G_k}{X_{k+1} - X_k} + \frac{L}{2}\fronorm{X_{k+1} - X_k}^2 \nonumber\\
           		&= f(X_k) - \gamma_k\nu_k\dotpF{G_k}{U_k} + \frac{L}{2}\gamma_k^2\nu_k^2\fronorm{U_k}^2. \label{eqn:proof_thm1_1}
           	\end{align}
           	Now let us show that $r_k\coloneqq\rank(G_k) = \fronorm{U_k}^2$. If $G_k = \sum_{i=1}^{r_k} \sigma_iu_iv_i^\top$ is the SVD of $G_k$, then the orthogonal polar factor $U_k = \sum_{i=1}^{r_k} u_iv_i^\top$. Thus, its squared Frobenius norm is 
           	\begin{equation*}
           		\fronorm{U_k}^2 = \tr(U_k^\top U_k) = \tr\left(\sum_{i,j=1}^{r_k} v_iu_iu_j^\top v_j^\top\right) = \sum_{i=1}^{r_k}v_i v_i^\top = \tr(I_{r_k}) = r_k\coloneqq\rank(G_k). 
           	\end{equation*}
           	We also recall that $\dotpF{G_k}{U_k} = \nucnorm{G_k}\eqqcolon \nu_k$ since $U_k = \argmax_{U:\specnorm{U}\le1}\dotpF{G_k}{U}$. 
           	Therefore, plugging into \eqref{eqn:proof_thm1_1}, we have 
           	\begin{equation*}
           		f(X_{k+1}) \le f(X_k) - \nu_k^2\left(\gamma_k - \frac{L}{2}\gamma_k^2 r_k\right). 
           	\end{equation*}
           	To ensure descent, we choose $\gamma_k = 1/(Lr_k)$ so that we have 
           	\begin{equation}\label{eqn:proof_thm1_2}
           		f(X_{k+1}) \le f(X_k) - \frac{1}{2Lr_k}\nu_k^2.
           	\end{equation}
           	Using the inequality $\fronorm{\cdot} \le \nucnorm{\cdot}$, we have 
           	\begin{equation}\label{eqn:proof_thm1_3}
           		\nucnorm{G_k} \ge \fronorm{G_k},
           	\end{equation} 
           	which implies $\nu_k^2 \ge \fronorm{G_k}^2$.         
           	By the $\mu$-\PL condition of $f$ \eqref{eqn:strong_cvx}, we also have $\fronorm{G_k}^2 \ge 2\mu(f(X_k) - \fstar)$. Plugging into \eqref{eqn:proof_thm1_2}, we obtain 
           	\[
           	f(X_{k+1}) - \fstar \le \left(1 - \frac{1}{r_k\kappa_H}\right)\left(f(X_k) - \fstar\right). 
           	\]
           	If we choose $r_{\max} \ge r_k$ for all $k\in\NN^*$, then we have 
           	\[
           	f(X_{k+1}) - \fstar \le \left(1 - \frac{1}{r_{\max}\kappa_H}\right)\left(f(X_k) - \fstar\right), 
           	\]
           	which implies
           	\begin{align*}
           		f(X_k) - \fstar &\le \left(1 - \frac{1}{r_{\max}\kappa_H}\right)^{\negthickspace k}\left(f(X_0) - \fstar\right) \\
           		&\le \exp\left(-\frac{k}{r_{\max}\kappa_H}\right)\left(f(X_0) - \fstar\right)  = \scrO(\exp(-k/r_{\max}\kappa_H)). 
           	\end{align*}
           	For the second bound in terms of the gradient condition number $\kappa_{G_k}$, notice that 
           	\[\fronorm{G_k}^2 = \sum_{i=1}^{r_k} \sigma_i^2 \le r_k \sigma_1^2 = r_k\kappa_{G_k}^2\sigma_{r_k}^2 \Rightarrow \sigma_{r_k}^2 \ge \frac1{r_k\cdot\kappa_{G_k}^2}\fronorm{G_k}^2. \]
           	Applying the Cauchy--Schwarz's inequality \eqref{eqn:proof_thm1_3} again, we have 
           	\[\nu_k^2 = \nucnorm{G_k}^2 \ge r_k^2\sigma_{r_k}^2 \ge r_k^2\cdot\frac1{r_k\cdot\kappa_{G_k}^2}\fronorm{G_k}^2 = \frac{r_k}{\kappa_{G_k}^2}\fronorm{G_k}^2. \]
           	Then, we can deduce from \eqref{eqn:proof_thm1_2} and the $\mu$-\PL condition of $f$ that 
           	\begin{align*}
           		f(X_{k+1}) &\le f(X_k) - \frac{1}{2Lr_k}\cdot \frac{r_k}{\kappa_{G_k}^2}\fronorm{G_k}^2 \\
           		&= f(X_k) - \frac{1}{2L\kappa_{G_k}^2}\fronorm{G_k}^2 \\
           		&\le f(X_k) - \frac{\mu}{L\kappa_{G_k}^2}(f(X_k) - \fstar). 
           	\end{align*}
           	Thus, we can conclude that 
           	\[f(X_{k+1}) - \fstar \le \left(1 - \frac{1}{\kappa_H\cdot\kappa_{G_k}^2}\right)\left(f(X_k) - \fstar\right). \]
        \end{proof}

        \begin{proof}[Proof of \Cref{thm:polarsgd_strcvx}]
           	By the $L$-Lipschitz smoothness of $f$, we have 
           	\begin{align*}
           		f(X_{k+1}) &\le f(X_k) + \dotpF{\nabla f(X_k)}{X_{k+1} - X_k} + \frac{L}{2}\fronorm{X_{k+1} - X_k}^2 \\
           		&= f(X_k) - \gamma\hnu_k\dotpF{G_k}{\hU_k} + \frac{L}{2}\gamma^2\hnu_k^2\fronorm{\hU_k}^2,
           	\end{align*}
           	where $\hnu_k \coloneqq \nucnorm{\hG_k}$ and $\hU_k\hH_k = \polar(\hG_k)$. Taking expectation on both sides, we obtain 
           	\begin{equation}\label{eqn:proof_thm2_1}
           		\Ex[f(X_{k+1})] \le f(X_k) - \gamma\Ex\left[ \hnu_k\dotpF{G_k}{\hU_k}\right]  + \frac{L}{2}\gamma^2\Ex\left[ \hnu_k^2\fronorm{\hU_k}^2\right]. 
           	\end{equation}
           	By \Cref{assum:grad_noise}, we can write $\hG_k = G_k + Z_k$, where $\Ex Z_k = 0$ and $\Ex \fronorm{Z_k}^2\le\varsigma^2$. Therefore, we have
           	\begin{equation}\label{eqn:proof_thm2_2}
           		\Ex\left[ \hnu_k\dotpF{G_k}{\hU_k}\right] = \Ex\left[\hnu_k\left( \dotpF{\hG_k}{\hU_k} - \dotpF{Z_k}{\hU_k} \right) \right] = \Ex \hnu_k^2 - \Ex\left[ \hnu_k\dotpF{Z_k}{\hU_k}\right]. 
           	\end{equation}        
           	Using $\fronorm{\cdot} \le \nucnorm{\cdot}$ and Jensen's inequality, we have 
           	\begin{equation}\label{eqn:proof_thm2_3}
           		\Ex\hnu_k^2 \ge \Ex\fronorm{\hG_k}^2 \ge \fronorm{\Ex\hG_k}^2 = \fronorm{G_k}^2. 
           	\end{equation}
           	On the other hand, by Cauchy--Schwarz's inequality, we also have
           	\begin{equation}\label{eqn:proof_thm2_4}
           		\Ex\left[\hnu_k\dotpF{Z_k}{\hU_k}\right] \le \sqrt{\Ex\hnu_k^2\cdot\Ex\dotpF{Z_k}{\hU_k}^2 }. 
           	\end{equation}
           	The first term on the right hand side can be upper bounded by the inequality $\nucnorm{G} \le \sqrt{\rank(G)}\fronorm{G}$ for any $G\in\RR^{m\times n}$:
           	\begin{equation}\label{eqn:proof_thm2_5}
           		\Ex\hnu_k^2 \le r_k\Ex\fronorm{\hG_k}^2 \le r_{\max}\Ex\fronorm{\hG_k}^2 \le r_{\max}\left(\varsigma^2 + \fronorm{G_k}^2\right), 
           	\end{equation}
           	where the last inequality is by \Cref{assum:grad_noise}. The second term can be upper bounded by Cauchy--Schwarz's inequality again: 
           	\begin{equation}\label{eqn:proof_thm2_6}
           		\Ex\dotpF{Z_k}{\hU_k}^2 \le \Ex\left[\fronorm{Z_k}^2 \fronorm{\hU_k}^2\right] = \Ex\left[\fronorm{Z_k}^2 \cdot\rank(\hG_k)\right] \le r_{\max} \Ex\fronorm{Z_k}^2\le \varsigma^2 r_{\max}. 
           	\end{equation}
           	Now, plugging \eqref{eqn:proof_thm2_3}, \eqref{eqn:proof_thm2_4}, \eqref{eqn:proof_thm2_5} and \eqref{eqn:proof_thm2_6} into \eqref{eqn:proof_thm2_2}, we obtain 
           	\begin{equation}\label{eqn:proof_thm2_7}
           		\Ex\left[ \hnu_k\dotpF{G_k}{\hU_k}\right]  \ge \fronorm{G_k}^2 -\varsigma r_{\max}\sqrt{\varsigma^2 + \fronorm{G_k}^2}.
           	\end{equation}
           	Furthermore, we can also use \eqref{eqn:proof_thm2_5} to bound 
           	\begin{equation}\label{eqn:proof_thm2_8}
           		\Ex\left[ \hnu_k^2\fronorm{\hU_k}^2\right] = \Ex\left[\hnu_k^2 \cdot\rank(\hG_k)\right] \le r_{\max}\Ex\hnu_k^2 \le r_{\max}^2\left(\varsigma^2 + \fronorm{G_k}^2\right). 
           	\end{equation}
           	Hence, putting \eqref{eqn:proof_thm2_7} and \eqref{eqn:proof_thm2_8} into \eqref{eqn:proof_thm2_1}, we obtain 
           	\[\Ex[f(X_{k+1})] \le f(X_k) - \left(\gamma - \frac{L}{2}\gamma^2r_{\max}^2\right) \fronorm{G_k}^2 + \varsigma\gamma r_{\max}\sqrt{\varsigma^2 + \fronorm{G_k}^2} + \frac{L}{2}\gamma^2\varsigma^2 r_{\max}^2. \]
           	Now, let us define $\Delta_k\coloneqq f(X_k) - \fstar$. Then, by the $\mu$-\PL condition of $f$, we have $\fronorm{G_k}^2 \ge 2\mu\Delta_k$, so the above bound can be rewritten as 
           	\begin{align*}
           		\Ex[\Delta_{k+1}] &\le \left(1 -2\mu \left(\gamma - \frac{L}{2}\gamma^2r_{\max}^2\right)\right)  \Delta_k + \varsigma\gamma r_{\max}\sqrt{\varsigma^2 + \fronorm{G_k}^2} + \frac{L}{2}\gamma^2\varsigma^2 r_{\max}^2 \\
           		&\le \left(1 -2\mu \left(\gamma - \frac{L}{2}\gamma^2r_{\max}^2\right)\right)  \Delta_k + \varsigma\gamma r_{\max}\left( \varsigma + \fronorm{G_k}\right)  + \frac{L}{2}\gamma^2\varsigma^2 r_{\max}^2, 
           	\end{align*}
           	since $\sqrt{a^2 + b^2} \le |a| + |b|$ for any $a,b\in\RR$. 
           	Furthermore, by the $L$-Lipschitz smoothness of $f$, we have $\fronorm{G_k}^2\le 2L\Delta_k$, implying that
           	\begin{equation*}
           		\Ex[\Delta_{k+1}] \le \left(1 -2\mu \left(\gamma - \frac{L}{2}\gamma^2r_{\max}^2\right)\right)  \Delta_k + \varsigma\gamma r_{\max}\left( \varsigma + \sqrt{2L\Delta_k}\right)  + \frac{L}{2}\gamma^2\varsigma^2 r_{\max}^2. 
           	\end{equation*}
           	Now, we invoke the A.M.-G.M.~inequality $ab\le \frac{a^2}{2\varepsilon} + \frac{\varepsilon b^2}2$ for any $a,b\in\Rp$ and $\varepsilon>0$, with $a=\sqrt{\Delta_k} $ and $b=\varsigma\gamma r_{\max}\sqrt{2L}$. Then we have $\varsigma\gamma r_{\max}\sqrt{2L\Delta_k} \le \Delta_k/(2\varepsilon) + \varepsilon L\gamma^2\varsigma^2r_{\max}^2$. Combining this inequality implies 
           	\begin{equation*}
           		\Ex[\Delta_{k+1}] \le \left(1 -2\mu \left(\gamma - \frac{L}{2}\gamma^2r_{\max}^2\right) + \frac{1}{2\varepsilon}\right)  \Delta_k + \varsigma^2\gamma r_{\max}\left( 1 + L\gamma r_{\max}\left(\varepsilon+\frac12 \right) \right).
           	\end{equation*}
           	Now, let $C_1\coloneqq 2\mu (\gamma - \frac{L}{2}\gamma^2r_{\max}^2) - 1/(2\varepsilon) > 0$, then we have the recursion 
           	\begin{equation}\label{eqn:proof_thm2_9}
           		\Ex[\Delta_{k+1}] \le (1- C_1)\Delta_k + \varsigma^2\gamma r_{\max}\left( 1 + L\gamma r_{\max}\left(\varepsilon+\frac12 \right) \right). 
           	\end{equation}
           	Note that we need $\gamma>0$ and $0<1-C_1<1$. With $\kappa_H \coloneqq L/\mu$, solving these inequalities yields an upper bound of the constant learning rate $\gamma$, given by 
           	\[\gamma < \gamma_{\max} \coloneqq\frac{1 + \sqrt{1-r_{\max}^2\kappa_H/(2\varepsilon)}}{Lr_{\max}^2}, \]
           	which is valid only if we choose $\varepsilon > r_{\max}^2\kappa_H/2$. Since $\gamma_{\max} > 1/(Lr_{\max}^2)$ if we choose $\varepsilon > r_{\max}^2\kappa_H/2$, we can choose a more conservative constant learning rate $\gamma \le 1/(Lr_{\max}^2)$ for simplicity.         	
           	Then, defining $C(\varepsilon) \coloneqq \gamma r_{\max}( 1 + L\gamma r_{\max}(\varepsilon+1/2 ) )$, the recursion \eqref{eqn:proof_thm2_9} becomes 
           	\begin{equation*}
           		\Ex[\Delta_{k+1}] \le (1- C_1)\Delta_k + C(\varepsilon) \varsigma^2. 
           	\end{equation*}
           	By a simple induction argument, we obtain that 
           	\begin{align*}
           		\Ex[\Delta_k] &\le (1- C_1)^k\left( \Delta_0 - \frac{C(\varepsilon)\varsigma^2}{C_1}\right) + \frac{C(\varepsilon)\varsigma^2}{C_1} \\
           		&\le \left( \Delta_0 - \frac{C(\varepsilon)\varsigma^2}{C_1}\right) \exp(-C_1k) + \frac{C(\varepsilon)\varsigma^2}{C_1} \\
           		&= \scrO\left(\exp(-C_1k) + C_2\varsigma^2\right), 
           	\end{align*}
           	where $C_2 \coloneqq C(\varepsilon) / C_1$. 
        \end{proof}

        \begin{proof}[Proof of \Cref{thm:mat_signSGD_strcvx}]
           	We first prove the convergence rate of matrix sign descent.         	
           	By the $L$-Lipschitz smoothness of $f$, we have         
           	\begin{align}
           		f(X_{k+1}) &\le f(X_k) + \dotpF{\nabla f(X_k)}{X_{k+1} - X_k} + \frac{L}{2}\fronorm{X_{k+1} - X_k}^2 \nonumber\\
           		&= f(X_k) - \gamma\dotpF{\nabla f(X_k)}{U_k} + \frac{L}{2}\gamma^2\fronorm{U_k}^2 \nonumber\\
           		&= f(X_k) - \gamma\nucnorm{\nabla f(X_k)} + \frac{L}{2}\gamma^2r_k \nonumber\\
           		&\le f(X_k) - \gamma\fronorm{\nabla f(X_k)} + \frac{L}{2}\gamma^2r_{\max}, \label{eqn:proof_thm3_0}
           	\end{align}
           	since $\fronorm{\cdot} \le \nucnorm{\cdot}$ and $r_k\le r_{\max}$ for all $k\in\setK$. 
            
            Now, let us define $\Delta_k\coloneqq f(X_k) - \fstar$. Then, by the $\mu$-\PL condition of $f$, we have $\fronorm{\nabla f(X_k)}^2 \ge 2\mu\Delta_k$, leading to the following nonlinear recursion: 
            \begin{equation*}
                \Delta_{k+1} \le \Delta_k - \gamma\sqrt{2\mu\Delta_k} + \frac{L}{2}\gamma^2r_{\max},  
            \end{equation*}
    		which converges at most sublinearly. 
    		
           	On the other hand, rearranging terms in \eqref{eqn:proof_thm3_0} gives 
           	\[ \gamma\fronorm{\nabla f(X_k)} \le f(X_k) - f(X_{k+1})+ \frac{L}{2}\gamma^2r_{\max}. \]
           	Summing $k$ from $1$ to $K$ yields 
           	\begin{align*}
           		\min_{k\in\setK}\fronorm{\nabla f(X_k)}  \le \frac1K\sumK \fronorm{\nabla f(X_k)} &\le \frac{1}{\gamma K}(f(X_1) - f(X_{K+1})) + \frac{L\gamma r_{\max}}2 \\
           		&\le \frac{1}{\gamma K}(f(X_1) - \fstar) + \frac{L\gamma r_{\max}}2 \\
           		&\le \scrO\left(\frac1{\gamma K} + \frac{L\gamma r_{\max}}2\right). 
           	\end{align*}
           	
           	Next, we prove the convergence rate of matrix \signSGD.         
           	Again, by the $L$-Lipschitz smoothness of $f$, we have 
           	\begin{align*}
           		f(X_{k+1}) &\le f(X_k) + \dotpF{\nabla f(X_k)}{X_{k+1} - X_k} + \frac{L}{2}\fronorm{X_{k+1} - X_k}^2 \\
           		&= f(X_K) - \gamma\dotpF{G_k}{\hU_k} + \frac{L}{2}\gamma^2\fronorm{\hU_k}^2,
           	\end{align*}
           	where $\hU_k\hH_k = \polar(\hG_k)$. 
           	Taking expectation on both sides, we have 
           	\begin{equation}\label{eqn:proof_thm3_1}
           		\Ex[f(X_{k+1})] \le f(X_k) - \gamma\Ex\dotpF{G_k}{\hU_k} + \frac{L}{2}\gamma^2\Ex[\fronorm{\hU_k}^2].
           	\end{equation}
           	By \Cref{assum:grad_noise}, we can write $\hG_k = G_k + Z_k$, where $\Ex Z_k = 0$ and $\Ex \fronorm{Z_k}^2\le\varsigma^2$. 
           	Then we have 
           	\begin{equation}\label{eqn:proof_thm3_2}
           		\dotpF{G_k}{\hU_k} = \dotpF{\hG_k}{\hU_k} - \dotpF{Z_k}{\hU_k} = \nucnorm{\hG_k} - \dotpF{Z_k}{\hU_k}. 
           	\end{equation}
           	By Cauchy--Schwarz's inequality, we have 
           	\begin{align}
           		\Ex \dotpF{Z_k}{\hU_k} &\le \Ex\left[\fronorm{Z_k}\fronorm{\hU_k}\right] \nonumber\\
           		&\le \sqrt{r_{\max}}\,\Ex\fronorm{Z_k} & \text{since }\fronorm{\hU_k}^2 = \rank(\hG_k) \le r_{\max}  \nonumber\\
           		&\le \sqrt{r_{\max}}\sqrt{\Ex\fronorm{Z_k}^2} & \text{by Jensen's inequality}  \nonumber\\
           		&\le \varsigma\sqrt{r_{\max}}. \label{eqn:proof_thm3_3}
           	\end{align}
           	On the other hand, by Jensen's inequality, we also have 
           	\begin{equation}\label{eqn:proof_thm3_4}
           		\Ex\nucnorm{\hG_k} \ge \Ex\fronorm{\hG_k} \ge \fronorm{\Ex \hG_k} = \fronorm{G_k}. 
           	\end{equation}
           	Consequently, taking expectation on both sides of \eqref{eqn:proof_thm3_2} and plugging in \eqref{eqn:proof_thm3_3} and \eqref{eqn:proof_thm3_4} gives 
           	\[\Ex\dotpF{G_k}{\hU_k} \ge \fronorm{G_k} - \varsigma\sqrt{r_{\max}}. \]
           	Again, since $\fronorm{\hU_k}^2 = \rank(\hG_k) \le r_{\max}$, we can derive from \eqref{eqn:proof_thm3_1} that 
           	\begin{equation*}
           		\Ex[f(X_{k+1})] \le f(X_k) - \gamma\fronorm{G_k} +\gamma \varsigma\sqrt{r_{\max}} + \frac{L}{2}\gamma^2r_{\max}.
           	\end{equation*}         
           	Rearranging terms yields
           	\[\fronorm{G_k} \le \frac1\gamma\Ex[f(X_k) - f(X_{k+1})] + \frac{L\gamma r_{\max}}2 + \varsigma\sqrt{r_{\max}}. \]
           	Summing $k$ from $1$ to $K$ yields 
           	\begin{align*}
           		\min_{k\in\setK}\fronorm{\nabla f(X_k)}  \le \frac1K\sumK \fronorm{\nabla f(X_k)} &\le \frac{1}{\gamma K}\Ex[f(X_1) - f(X_{K+1})] + \frac{L\gamma r_{\max}}2 + \varsigma\sqrt{r_{\max}}\\
           		&\le \frac{1}{\gamma K}\Ex[f(X_1) - \fstar] + \frac{L\gamma r_{\max}}2 + \varsigma\sqrt{r_{\max}} \\
           		&\le \scrO\left(\frac1{\gamma K} + \frac{L\gamma r_{\max}}2 + \varsigma\sqrt{r_{\max}}\right). 
           	\end{align*}
        \end{proof}
        
        \subsection{Proofs for \Cref{subsec:conv_inexact}}
        \label{subsec:proofs_inexact}
        
        \begin{proof}[Proof of \Cref{thm:polargrad_strcvx_inexact}]
            From \Cref{assum:inexact_polar_oracle}(i), we can characterize an \emph{alignment defect}. By H\"{o}lder's inequality, we have 
            \begin{align*}
                \dotpF{G_k}{\tU_k} & = \dotpF{G_k}{U_k} + \dotpF{G_k}{\tU_k - U_k} \\
                &\ge \nucnorm{G_k} - \nucnorm{G_k}\specnorm{\tU_k - U_k} \\
                &\ge (1-\varepsilon_k)\nucnorm{G_k}.  
            \end{align*}
            Let us recall that $\tnu_k \coloneqq\dotpF{\tU_k}{G_k}$ and $\nu_k \coloneqq \nucnorm{G_k}$. The above inequality is equivalent to 
            \begin{equation}\label{eqn:inexact_polar_oracle_1}
                \tnu_k \ge (1-\varepsilon_k) \nu_k. 
            \end{equation}
            From \Cref{assum:inexact_polar_oracle}(ii), we can characterize an \emph{orthogonality defect}: 
            \begin{align*}
                \fronorm{\tU_k}^2 &= \tr\left(\tU_k^\top\tU_k\right) = \tr\left(I_{r_k} + (\tU_k^\top\tU_k - I_{r_k})\right) \nonumber \\
                &= r_k + \tr\left(\tU_k^\top\tU_k - I_{r_k}\right). 
            \end{align*}
            Since $\tU_k^\top\tU_k - I_{r_k}$ is symmetric, we have $|\tr(\tU_k^\top\tU_k - I_{r_k})| \le r_k\specnorm{\tU_k^\top\tU_k - I_{r_k}}$. This implies that 
            \begin{equation}\label{eqn:inexact_polar_oracle_2}
                \fronorm{\tU_k}^2 \le r_k(1+\delta_k). 
            \end{equation}
           We obtain a new descent lemma for \eqref{eqn:polargrad_practical} via the $L$-Lipschitz smoothness of $f$: 
           \begin{align*}
                f(X_{k+1}) &\le f(X_k) - \gamma_k\tnu_k^2 + \frac{L}{2}\gamma_k^2\tnu_k^2\fronorm{\tU_k}^2 \\
                &= f(X_k) + \tnu_k^2\left(-\gamma_k + \frac{L}{2}\gamma_k^2\fronorm{\tU_k}^2\right). 
           \end{align*}
            Now, if we choose $\gamma_k\le \tfrac{c}{L\fronorm{\tU_k}^2}$ for some $c\in\left(0,1\right]$, we have $-\gamma_k + \frac{L}{2}\gamma_k^2\fronorm{\tU_k}^2 \le -\left(1-c/2\right)\gamma_k$. 
            Then, using \eqref{eqn:inexact_polar_oracle_1}, we have 
            \[f(X_{k+1}) \le f(X_k) - \left(1-\frac{c}{2}\right)\gamma_k(1-\varepsilon_k)^2\nu_k^2. \]
            Since $\nu_k \ge \fronorm{G_k}$ and $f$ is $\mu$-\PL, i.e., $\fronorm{G_k}^2 \ge 2\mu\left(f(X_k) - f^\star\right)$, we have 
            \begin{align*}
                f(X_{k+1}) &\le f(X_k) - \left(1 - \frac{c}{2}\right)\gamma_k(1-\varepsilon_k)^2\fronorm{G_k}^2 \\
                &\le f(X_k) - \left(1 - \frac{c}{2}\right)\frac{c}{L\fronorm{\tU_k}^2}(1-\varepsilon_k)^2 \cdot2\mu\left(f(X_k) - f^\star\right) & \text{since $f$ is $\mu$-\PL}\\
                &\le f(X_k) - 2\left(1 - \frac{c}{2}\right)c\frac{(1-\varepsilon_k)^2}{Lr_k(1+\delta_k)}\cdot\mu\left(f(X_k) - f^\star\right) & \text{by \eqref{eqn:inexact_polar_oracle_2}.}
            \end{align*}
            Therefore, the above inequality yields
            \[f(X_{k+1}) - f^\star \le \left(1 - \frac{2c}{\kappa_Hr_k}\left(1 - \frac{c}{2}\right)\frac{(1-\varepsilon_k)^2}{1+\delta_k}\right)\left(f(X_k) - f^\star\right). \]
            Furthermore, since $\fronorm{\tU_k}^2 \le r_k(1+\delta_k) \le r_{\max}(1+\delta_{\max})$ and $\varepsilon_k\le\varepsilon_{\max}$ for all $k\in\NN$, if we apply a constant learning rate $\gamma \coloneqq c/(Lr_{\max}(1 + \delta_{\max}))$ for some $c\in\left(0,1\right]$, we obtain the desired uniform bound:
            \[
                f(X_{k+1}) - \fstar \le \left(1 - \frac{2c}{r_{\max}\kappa_H}\left(1-\frac{c}2\right) \frac{(1-\varepsilon_{\max})^2}{1+\delta_{\max}} \right)(f(X_k) - \fstar).
            \]
        \end{proof}

        \begin{proof}[Proof of \Cref{thm:polarsgd_strcvx_inexact}]
            This proof largely resembles that of \Cref{thm:polarsgd_strcvx}. 
            By the $L$-Lipschitz smoothness of $f$, we have 
               	\begin{align*}
               		f(X_{k+1}) &\le f(X_k) + \dotpF{\nabla f(X_k)}{X_{k+1} - X_k} + \frac{L}{2}\fronorm{X_{k+1} - X_k}^2 \\
               		&= f(X_k) - \gamma\tnu_k\dotpF{G_k}{\tU_k} + \frac{L}{2}\gamma^2\tnu_k^2\fronorm{\tU_k}^2,
               	\end{align*}
               	where $\tnu_k \coloneqq \dotpF{\hG_k}{\tU_k}$ and $\tU_k\tH_k = \polarhat(\hG_k)$. Taking expectation on both sides, we obtain 
               	\begin{equation}\label{eqn:proof_thm5_1}
               		\Ex[f(X_{k+1})] \le f(X_k) - \gamma\Ex\left[\tnu_k\dotpF{G_k}{\tU_k}\right]  + \frac{L}{2}\gamma^2\Ex\left[ \tnu_k^2\fronorm{\tU_k}^2\right]. 
               	\end{equation}
               Similar to the derivation of \eqref{eqn:inexact_polar_oracle_1}, from \Cref{assum:inexact_polar_oracle}(i), we have 
               \begin{equation}\label{eqn:inexact_polar_oracle_3}
                   \tnu_k \ge (1-\hat\varepsilon_k) \hnu_k. 
               \end{equation}
               Also, from \Cref{assum:inexact_polar_oracle}(ii), we can also deduce that 
               \begin{equation}\label{eqn:inexact_polar_oracle_4}
                   \fronorm{\tU_k}^2 \le \hat{r}_k(1+\hat\delta_k). 
               \end{equation}               
               	By \Cref{assum:grad_noise}, we have
               	\begin{equation}\label{eqn:proof_thm5_2}
               		\Ex\left[ \tnu_k\dotpF{G_k}{\tU_k}\right] = \Ex\left[\tnu_k\left( \dotpF{\hG_k}{\tU_k} - \dotpF{Z_k}{\tU_k} \right) \right] = \Ex \tnu_k^2 - \Ex\left[ \tnu_k\dotpF{Z_k}{\tU_k}\right]. 
               	\end{equation}        
               	Using \eqref{eqn:inexact_polar_oracle_3}, $\fronorm{\cdot} \le \nucnorm{\cdot}$ and Jensen's inequality, we have 
               	\begin{equation}\label{eqn:proof_thm5_3}
               		\Ex\tnu_k^2 \ge (1-\hat\varepsilon_k)^2\Ex\hnu_k^2 \ge (1-\hat\varepsilon_k)^2\Ex\fronorm{\hG_k}^2 \ge (1-\hat\varepsilon_k)^2\fronorm{\Ex\hG_k}^2 = (1-\hat\varepsilon_k)^2\fronorm{G_k}^2. 
               	\end{equation}
               	On the other hand, by Cauchy--Schwarz's inequality, we also have
               	\begin{equation}\label{eqn:proof_thm5_4}
               		\Ex\left[\tnu_k\dotpF{Z_k}{\tU_k}\right] \le \sqrt{\Ex\tnu_k^2\cdot\Ex\dotpF{Z_k}{\tU_k}^2 }. 
               	\end{equation}
               	The first term on the right hand side can be upper bounded by Cauchy--Schwarz's inequality and \eqref{eqn:inexact_polar_oracle_4}:
               	\begin{equation}\label{eqn:proof_thm5_5}
               		\Ex\tnu_k^2 = \Ex \dotpF{\hG_k}{\tU_k}^2 \le\Ex\left[ \fronorm{\hG_k}^2\fronorm{\hU_k}^2\right]  \le \hat{r}_k(1+\hat\delta_k)\Ex\fronorm{\hG_k}^2 \le \hat{r}_k(1+\hat\delta_k)\left(\varsigma^2 + \fronorm{G_k}^2\right), 
               	\end{equation}
               	where the last inequality is by \Cref{assum:grad_noise}. The second term can be upper bounded by Cauchy--Schwarz's inequality again: 
               	\begin{equation}\label{eqn:proof_thm5_6}
               		\Ex\dotpF{Z_k}{\tU_k}^2 \le \Ex\left[\fronorm{Z_k}^2 \fronorm{\tU_k}^2\right] = \hat{r}_k(1+\hat\delta_k)\Ex\fronorm{Z_k}^2 \le \hat{r}_k(1+\hat\delta_k)\varsigma^2. 
               	\end{equation}
               	Now, plugging \eqref{eqn:proof_thm5_3}, \eqref{eqn:proof_thm5_4}, \eqref{eqn:proof_thm5_5} and \eqref{eqn:proof_thm5_6} into \eqref{eqn:proof_thm5_2}, we obtain 
               	\begin{equation}\label{eqn:proof_thm5_7}
               		\Ex\left[ \tnu_k\dotpF{G_k}{\tU_k}\right]  \ge (1-\hat\varepsilon_k)^2\fronorm{G_k}^2 -\varsigma \hat{r}_k(1+\hat\delta_k)\sqrt{\varsigma^2 + \fronorm{G_k}^2}.
               	\end{equation}
               	Furthermore, we can also use \eqref{eqn:proof_thm5_5} to bound 
               	\begin{equation}\label{eqn:proof_thm5_8}
               		\Ex\left[ \tnu_k^2\fronorm{\tU_k}^2\right] \le \hat{r}_k(1+\hat\delta_k)\Ex\tnu_k^2 \le \hat{r}_k^2(1+\hat\delta_k)^2\left(\varsigma^2 + \fronorm{G_k}^2\right). 
               	\end{equation}
               	Hence, putting \eqref{eqn:proof_thm5_7} and \eqref{eqn:proof_thm5_8} into \eqref{eqn:proof_thm5_1}, we obtain 
                   \begin{multline*}
                    \Ex[f(X_{k+1})] \le f(X_k) - \left(\gamma(1-\hat\varepsilon_k)^2 - \frac{L}{2}\gamma^2\hat{r}_k^2(1+\hat\delta_k)^2\right) \fronorm{G_k}^2 \\
                    + \varsigma\gamma \hat{r}_k(1+\hat\delta_k)\sqrt{\varsigma^2 + \fronorm{G_k}^2} + \frac{L}{2}\gamma^2\varsigma^2 \hat{r}_k^2(1+\hat\delta_k)^2.
                   \end{multline*}
               	Now, let us define $\Delta_k\coloneqq f(X_k) - \fstar$. Then, by the $\mu$-\PL condition of $f$, we have $\fronorm{G_k}^2 \ge 2\mu\Delta_k$, so the above bound can be rewritten as 
               	\begin{align*}
               		\Ex[\Delta_{k+1}] &\le \left(1 -2\mu \left(\gamma(1-\hat\varepsilon_k)^2 - \frac{L}{2}\gamma^2\hat{r}_k^2(1+\hat\delta_k)^2\right)\right)  \Delta_k \\
                    &\qquad+ \varsigma\gamma \hat{r}_k(1+\hat\delta_k)\sqrt{\varsigma^2 + \fronorm{G_k}^2} + \frac{L}{2}\gamma^2\varsigma^2 \hat{r}_k^2(1+\hat\delta_k)^2 \\
               		&\le \left(1 -2\mu \left(\gamma(1-\hat\varepsilon_k)^2 - \frac{L}{2}\gamma^2\hat{r}_k^2(1+\hat\delta_k)^2\right)\right)  \Delta_k \\
                    &\qquad+ \varsigma\gamma \hat{r}_k(1+\hat\delta_k)\left(\varsigma + \fronorm{G_k}\right) + \frac{L}{2}\gamma^2\varsigma^2 \hat{r}_k^2(1+\hat\delta_k)^2, 
               	\end{align*}
               	since $\sqrt{a^2 + b^2} \le |a| + |b|$ for any $a,b\in\RR$. 
               	Furthermore, by the $L$-Lipschitz smoothness of $f$, we have $\fronorm{G_k}^2\le 2L\Delta_k$, implying that
               	\begin{multline*}
               		\Ex[\Delta_{k+1}] \le \left(1 -2\mu \left(\gamma(1-\hat\varepsilon_k)^2 - \frac{L}{2}\gamma^2\hat{r}_k^2(1+\hat\delta_k)^2\right)\right) \Delta_k \\
                       + \varsigma\gamma \hat{r}_k(1+\hat\delta_k)\left( \varsigma + \sqrt{2L\Delta_k}\right) + \frac{L}{2}\gamma^2\varsigma^2 \hat{r}_k^2(1+\hat\delta_k)^2. 
               	\end{multline*}
               	Now, we invoke the A.M.-G.M.~inequality $ab\le \frac{a^2}{2\omega} + \frac{\omega b^2}2$ for any $a,b\in\Rp$ and $\omega>0$, with $a=\sqrt{\Delta_k}$ and $b=\varsigma\gamma \hat{r}_k(1+\hat\delta_k)\sqrt{2L}$. Then we have $\varsigma\gamma  \hat{r}_k(1+\hat\delta_k)\sqrt{2L\Delta_k} \le \Delta_k/(2\omega) + \omega L\gamma^2\varsigma^2\hat{r}_k^2(1+\hat\delta_k)^2$. Combining this inequality implies 
               	\begin{multline*}
               		\Ex[\Delta_{k+1}] \le \left(1 -2\mu \left(\gamma(1-\hat\varepsilon_k)^2 - \frac{L}{2}\gamma^2\hat{r}_k^2(1+\hat\delta_k)^2\right)+ \frac{1}{2\omega}\right)  \Delta_k \\
                       + \varsigma^2\gamma \hat{r}_k(1+\hat\delta_k)\left( 1 + L\gamma \hat{r}_k(1+\hat\delta_k)\left(\omega+\frac12 \right) \right).
               	\end{multline*}
               	Now, let $\tC_1\coloneqq 2\mu (\gamma(1-\hat\varepsilon_{\max})^2 - \frac{L}{2}\gamma^2\hat{r}_{\max}^2(1+\hat\delta_{\max})^2) - 1/(2\omega) > 0$, then we have the recursion 
               	\begin{equation}\label{eqn:proof_thm5_9}
               		\Ex[\Delta_{k+1}] \le (1- \tC_1)\Delta_k + \varsigma^2\gamma \hat{r}_{\max}(1+\hat\delta_{\max})\left( 1 + L\gamma \hat{r}_{\max}(1+\hat\delta_{\max})\left(\omega+\frac12 \right) \right). 
               	\end{equation}
               	Note that we need $\gamma>0$ and $0<1-\tC_1<1$. With $\kappa_H \coloneqq L/\mu$, solving these inequalities yields an upper bound of the constant learning rate $\gamma$, given by 
               	\[\gamma < \gamma_{\max} \coloneqq\frac{(1-\hat\varepsilon_{\max})^2 + \sqrt{(1-\hat\varepsilon_{\max})^4-\hat{r}_{\max}^2(1+\hat\delta_{\max})^2\kappa_H/(2\omega)}}{L\hat{r}_{\max}^2(1+\hat\delta_{\max})^2}, \]
               	which is valid only if we choose $\omega > \hat{r}_{\max}^2(1+\hat\delta_{\max})^2\kappa_H/(2(1-\hat\varepsilon_{\max})^4)$. Since $\gamma_{\max} > (1-\hat\varepsilon_{\max})^2/(L\hat{r}_{\max}^2(1+\hat\delta_{\max})^2)$ if we choose $\omega > \hat{r}_{\max}^2(1+\hat\delta_{\max})^2\kappa_H/(2(1-\hat\varepsilon_{\max})^4)$, we can choose a more conservative constant learning rate $\gamma \le (1-\hat\varepsilon_{\max})^2/(L\hat{r}_{\max}^2(1+\hat\delta_{\max})^2)$ for simplicity.         	
               	Then, defining $\tC(\omega) \coloneqq \gamma \hat{r}_{\max}(1+\hat\delta_{\max})( 1 + L\gamma \hat{r}_{\max}(1+\hat\delta_{\max})(\omega+1/2 ) )$, the recursion \eqref{eqn:proof_thm5_9} becomes 
               	\begin{equation*}
               		\Ex[\Delta_{k+1}] \le (1- \tC_1)\Delta_k + \tC(\omega) \varsigma^2. 
               	\end{equation*}
               	By a simple induction argument, we obtain that 
               	\begin{align*}
               		\Ex[\Delta_k] &\le (1- \tC_1)^k\left( \Delta_0 - \frac{\tC(\omega)\varsigma^2}{\tC_1}\right) + \frac{\tC(\omega)\varsigma^2}{\tC_1} \\
               		&\le \left( \Delta_0 - \frac{\tC(\omega)\varsigma^2}{\tC_1}\right) \exp(-\tC_1k) + \frac{\tC(\omega)\varsigma^2}{\tC_1} \\
               		&= \scrO\left(\exp(-\tC_1k) + \tC_2\varsigma^2\right), 
               	\end{align*}
               	where $\tC_2 \coloneqq \tC(\omega) / \tC_1$. 
        \end{proof}

        \begin{proof}[Proof of \Cref{thm:mat_signSGD_strcvx_inexact}]
            This proof is also similar to that of \Cref{thm:mat_signSGD_strcvx}. 
            We first prove the convergence rate of matrix sign descent.         	
           	By the $L$-Lipschitz smoothness of $f$, \eqref{eqn:inexact_polar_oracle_1} and \eqref{eqn:inexact_polar_oracle_2}, we have         
           	\begin{align}
           		f(X_{k+1}) &\le f(X_k) - \gamma\dotpF{\nabla f(X_k)}{\tU_k} + \frac{L}{2}\gamma^2\fronorm{\tU_k}^2 \nonumber\\
           		&= f(X_k) - \gamma(1-\varepsilon_k)\nu_k+ \frac{L}{2}\gamma^2r_k(1+\delta_k) \nonumber\\
           		&\le f(X_k) - \gamma(1-\varepsilon_k)\fronorm{\nabla f(X_k)} + \frac{L}{2}\gamma^2r_k(1+\delta_k), \label{eqn:proof_thm6_0}
           	\end{align}
           	since $\fronorm{\cdot} \le \nucnorm{\cdot}$ and $r_k\le r_{\max}$ for all $k\in\setK$. 
            
            Now, let us define $\Delta_k\coloneqq f(X_k) - \fstar$. Then, by the $\mu$-\PL condition of $f$, we have $\fronorm{\nabla f(X_k)}^2 \ge 2\mu\Delta_k$, leading to the following nonlinear recursion: 
            \begin{equation*}
                \Delta_{k+1} \le \Delta_k - \gamma(1-\varepsilon_k)\sqrt{2\mu\Delta_k} + \frac{L}{2}\gamma^2r_k(1+\delta_k),  
            \end{equation*}
    		which converges at most sublinearly. 
    		
           	On the other hand, rearranging terms in \eqref{eqn:proof_thm6_0} gives 
           	\[ \gamma(1-\varepsilon_{\max})\fronorm{\nabla f(X_k)} \le f(X_k) - f(X_{k+1})+ \frac{L}{2}\gamma^2r_{\max}(1+\delta_{\max}). \]
           	Summing $k$ from $1$ to $K$ yields 
           	\begin{align*}
           		\min_{k\in\setK}\fronorm{\nabla f(X_k)}  \le \frac1K\sumK \fronorm{\nabla f(X_k)} &\le \frac{1}{\gamma (1-\varepsilon_{\max})K}(f(X_1) - f(X_{K+1})) + \frac{L\gamma r_{\max}(1+\delta_{\max})}{2(1-\varepsilon_{\max})} \\
           		&\le \frac{1}{\gamma (1-\varepsilon_{\max})K}(f(X_1) - \fstar) + \frac{L\gamma r_{\max}(1+\delta_{\max})}{2(1-\varepsilon_{\max})} \\
           		&\le \scrO\left(\frac1{\gamma (1-\varepsilon_{\max})K} + \frac{L\gamma r_{\max}(1+\delta_{\max})}{2(1-\varepsilon_{\max})} \right). 
           	\end{align*}
           	
           	Next, we prove the convergence rate of matrix \signSGD.         
           	Again, by the $L$-Lipschitz smoothness of $f$, we have 
           	\begin{equation*}
           		f(X_{k+1}) \le f(X_K) - \gamma\dotpF{G_k}{\tU_k} + \frac{L}{2}\gamma^2\fronorm{\tU_k}^2,
           	\end{equation*}
           	where $\tU_k\tH_k = \polarhat(\hG_k)$. 
           	Taking expectation on both sides, we have 
           	\begin{equation}\label{eqn:proof_thm6_1}
           		\Ex[f(X_{k+1})] \le f(X_k) - \gamma\Ex\dotpF{G_k}{\tU_k} + \frac{L}{2}\gamma^2\Ex[\fronorm{\tU_k}^2].
           	\end{equation}
           	By \Cref{assum:grad_noise}, we can write $\hG_k = G_k + Z_k$, where $\Ex Z_k = 0$ and $\Ex \fronorm{Z_k}^2\le\varsigma^2$. 
           	Then, by \eqref{eqn:inexact_polar_oracle_3}, we have 
           	\begin{equation}\label{eqn:proof_thm6_2}
           		\dotpF{G_k}{\tU_k} = \dotpF{\hG_k}{\tU_k} - \dotpF{Z_k}{\tU_k} \ge (1-\hat\varepsilon_k) \hnu_k - \dotpF{Z_k}{\tU_k}. 
           	\end{equation}
           	Applying Cauchy--Schwarz's inequality twice and \eqref{eqn:inexact_polar_oracle_4}, we have 
           	\begin{align}
           		\Ex \dotpF{Z_k}{\tU_k} &\le \Ex\left[\fronorm{Z_k}\fronorm{\tU_k}\right] \nonumber\\
           		&\le \sqrt{\Ex\fronorm{Z_k}^2\cdot\Ex\fronorm{\tU_k}^2} \nonumber\\
           		&\le \sqrt{\hat{r}_k(1+\hat\delta_k)}\sqrt{\Ex\fronorm{Z_k}^2} \nonumber\\
           		&\le \varsigma\sqrt{\hat{r}_k(1+\hat\delta_k)}. \label{eqn:proof_thm6_3}
           	\end{align}
           	On the other hand, by Jensen's inequality, we also have 
           	\begin{equation}\label{eqn:proof_thm6_4}
           		\Ex\nucnorm{\hG_k} \ge \Ex\fronorm{\hG_k} \ge \fronorm{\Ex \hG_k} = \fronorm{G_k}. 
           	\end{equation}
           	Consequently, taking expectation on both sides of \eqref{eqn:proof_thm6_2} and plugging in \eqref{eqn:proof_thm6_3} and \eqref{eqn:proof_thm6_4} gives 
           	\[\Ex\dotpF{G_k}{\tU_k} \ge (1-\hat\varepsilon_k)\fronorm{G_k} - \varsigma\sqrt{\hat{r}_k(1+\hat\delta_k)}. \]
           	Again, by \eqref{eqn:inexact_polar_oracle_4}, we can derive from \eqref{eqn:proof_thm6_1} that 
           	\begin{equation*}
           		\Ex[f(X_{k+1})] \le f(X_k) - \gamma(1-\hat\varepsilon_k)\fronorm{G_k} +\gamma \varsigma\sqrt{\hat{r}_k(1+\hat\delta_k)} + \frac{L}{2}\gamma^2\hat{r}_k(1+\hat\delta_k).
           	\end{equation*}         
           	Rearranging terms yields
           \begin{align*}
            \fronorm{G_k} &\le \frac1{\gamma(1-\hat\varepsilon_k)}\Ex[f(X_k) - f(X_{k+1})] + \frac{L\gamma \hat{r}_k(1+\hat\delta_k)}{2(1-\hat\varepsilon_k)} + \frac{\varsigma\sqrt{\hat{r}_k(1+\hat\delta_k)}}{1-\hat\varepsilon_k} \\
            &\le \frac1{\gamma(1-\hat\varepsilon_{\max})}\Ex[f(X_k) - f(X_{k+1})] + \frac{L\gamma \hat{r}_{\max}(1+\hat\delta_{\max})}{2(1-\hat\varepsilon_{\max})} + \frac{\varsigma\sqrt{\hat{r}_{\max}(1+\hat\delta_{\max})}}{1-\hat\varepsilon_{\max}}. 
           \end{align*}
           	Summing $k$ from $1$ to $K$ yields 
           	\begin{align*}
           		&\quad\,\min_{k\in\setK}\fronorm{\nabla f(X_k)} \le \frac1K\sumK \fronorm{\nabla f(X_k)} \\
                &\le \frac{1}{\gamma (1-\hat\varepsilon_{\max})K}\Ex[f(X_1) - f(X_{K+1})] + \frac{L\gamma \hat{r}_{\max}(1+\hat\delta_{\max})}{2(1-\hat\varepsilon_{\max})} + \frac{\varsigma\sqrt{\hat{r}_{\max}(1+\hat\delta_{\max})}}{1-\hat\varepsilon_{\max}} \\
           		&\le \frac{1}{\gamma (1-\hat\varepsilon_{\max})K}\Ex[f(X_1) - \fstar] + \frac{L\gamma \hat{r}_{\max}(1+\hat\delta_{\max})}{2(1-\hat\varepsilon_{\max})} + \frac{\varsigma\sqrt{\hat{r}_{\max}(1+\hat\delta_{\max})}}{1-\hat\varepsilon_{\max}} \\
           		&\le \scrO\left(\frac1{\gamma (1-\hat\varepsilon_{\max})K} + \frac{L\gamma \hat{r}_{\max}(1+\hat\delta_{\max})}{2(1-\hat\varepsilon_{\max})} + \frac{\varsigma\sqrt{\hat{r}_{\max}(1+\hat\delta_{\max})}}{1-\hat\varepsilon_{\max}}\right). 
           	\end{align*}
            
        \end{proof}

        \begin{proof}[Proof of \Cref{thm:NS}]
            We first introduce some additional notation. Let us recall that $\tU_{k,0} = G_k/\fronorm{G_k}$ and the Newton--Schulz iteration with quintic polynomials is given by        
            \begin{equation}\label{eqn:NS}
                (\forall j\in\{0,\ldots,T\})\quad \tU_{k,j+1} = a\tU_{k,j} + b\tU_{k,j}M_{k,j} + c\tU_{k,j}M_{k,j}^2, \quad M_{k,j}\coloneqq \tU_{k,j}^\top \tU_{k,j}.
            \end{equation}
            In the following, we drop the dependence on $k$ for notation simplicity. Then, \eqref{eqn:NS} can be rewritten as 
            \[\tU_{j+1} = \tU_j p(M_j), \quad M_j\coloneqq \tU_j^\top\tU_j, \quad p(t) \coloneqq a+bt+ct^2. \]
            Hence, if $M_j$ has eigenvalue $t$, the corresponding singular value of $\tU_j$ is $\sqrt{t}$. After one iteration of \eqref{eqn:NS}, the new squared singular value is $t_+ = \varphi\coloneqq tp(t)^2$. Let us define $e\coloneqq t-1$. We are interested in the behavior of $\varphi$ near the orthogonal point $t=1$ and from that we can determine $(a,b,c)$. 
            
            We expand $\varphi$ at $1+e$ using the Taylor expansion:
            \[\varphi(1+e) = (a+b+c)^2 + (a+3b+5c)e + \alpha_2e^2 + \alpha_3e^3 + \scrO(e^4), \]
            where 
            \[\alpha_2\coloneqq b^2+4bc+2b+4c^2 +2(bc+c^2)+2abc(b+2c) ,\] 
            and 
            \[\alpha_3\coloneqq b^2+6bc+8c^2+2c(a+b+c). \]
            Now, solving the fixed-point condition $\varphi(1)=1$, $\varphi'(1) = 0$ (no linear term) and $\varphi''(1) = 0$ (no quadratic term), we have $(a,b,c)=(15/8,-5/4,3/8)$. Putting these back, we have
            \[\varphi(1+e) = 1 + \frac58 e^3 - \frac{15}{64}e^4 + \scrO(e^5). \]
            Hence, for $|e|$ small enough, there exists a constant $\zeta>0$ such that $|e_+| \le \zeta|e|^3$. 
            
            Let $\{\lambda_j^{(i)}\}_{i\in\{1,\ldots,n\}}$ be the eigenvalues of $M_j$. We define 
            \[e_j\coloneqq\specnorm{M_j - I} = \max_{i\in\{1,\ldots,n\}} |\lambda_j^{(i)} - 1|, \]
            and we also have $\lambda_{j+1}^{(i)} = \varphi(\lambda_j^{(i)})$. Consequently, $e_j^{(i)} \coloneqq \lambda_j^{(i)} - 1$ satisfies 
            \[|e_{j+1}^{(i)}| = |\varphi(1+e_j^{(i)}) - 1| \le\zeta|e_j^{(i)}|^3\]
            for sufficiently small $|e_j^{(i)}|\le\bar{e}$. Thus, we conclude that 
            \[e_{j+1} \coloneqq \max_{i\in\{1,\ldots,n\}} |e_{j+1}^{(i)}| \le \zeta\max_{i\in\{1,\ldots,n\}} |e_j^{(i)}|^3 \le \zeta e_j^3.\]
            Recursively, we have $e_T\le\zeta^{1+3+3^2+\cdots+3^{T-1}} e_0^{3^T} = C_T e_0^{3^T}$, with a moderate constant $C_T = \zeta^{(3^T-1)/2} = \scrO(1)$ for fixed small $T$. 
            
            Now, we determine the value of $e_0$. Adding back the index on $k$, let us recall that $e_{k,0} \coloneqq \specnorm{\tU_{k,0}^\top\tU_{k,0} - I} = \max_{i\in\{1,\ldots,n\}} |\sigma_i(\tU_{k,0})^2-1|$, where $\tU_{k,0} \coloneqq G_k/\fronorm{G_k}$. Therefore, we have 
            \[0 < \sigma_i(\tU_{k,0})^2 = \frac{\sigma_i(G_k)^2}{\fronorm{G_k}^2} \le \frac{\sigma_{\max}(G_k)^2}{\fronorm{G_k}^2} \le 1, \]
            i.e., we always have $\sigma_i(\tU_{k,0})^2 \in\left(0, 1\right]$ for each $i\in\{1,\ldots,n\}$. 
            The worst deviation is at the minimum singular value $e_{k,0} = 1 - \sigma_{\min}(G_k) / \fronorm{G_k}^2$, which gives 
            \[e_0 \coloneqq \max_{k\in\{0,\ldots,K\}} e_{k,0} = 1 - \min_{k\in\{0,\ldots,K\}} \frac{\sigma_{\min}(G_k)^2}{\fronorm{G_k}^2}. \] 
            Hence, $e_0$ depends on the ``Frobenius condition number'' 
            \[\kappa_{\mathrm{F}}(G_k)\coloneqq \frac{\fronorm{G_k}}{\sigma_{\min}(G_k)}, \quad e_{k,0}^{\mathrm{F}} \coloneqq 1 - \frac{1}{\kappa_{\mathrm{F}}(G_k)^2}. \]
            If we use the spectral norm for normalization, then the standard condition number determines the error bound 
            \[e_{k,0}^{\mathrm{S}} \coloneqq 1 - \frac{1}{\kappa_2(G_k)^2}. \]
            Since we always have $\kappa_{\mathrm{F}} \ge \kappa_2$, the Frobenius norm normalization gives a larger $e_0$ in the worst case. 
            
            Recall that we have $e_{k,T} \le C_\delta e_0^{3^T}$ for some constant $C_\delta\ge 1$ depending on $\zeta$, $\bar{e}$, $T$ but not $k$. This implies that 
            \[\delta_{\max}(T) \coloneqq \max_{k\in\{0,\ldots,K\}} \specnorm{\tU_{k,T}^\top\tU_{k,T} - I} \le C_\delta e_0^{3^T}. \]
    		If $\specnorm{\tU_{k,T}^\top\tU_{k,T} - I} \le \delta_{\max}(T)$ for all $k$ and the singular values of $\tU_{k,0}$ lie in $[\ell, 1]$ with $\ell>0$, then $\specnorm{\tU_{k,T} - U_k} \le  C_{\textrm{pol}}\cdot\delta_{\max}(T)$ for some constant $C_{\textrm{pol}}\in\left(0,1\right]$. 
    		To see this, we use the following perturbation argument. 
    		
    		We write $\tU_{k,T} = U_k + E_{k,T}$ for some small error matrix $E_{k,T}$. Then we have 
    		\[\tU_{k,T}^\top\tU_{k,T} - I = (U_k + E_{k,T})^\top(U_k+E_{k,T}) - I = U_k^\top E_{k,T} + E_{k,T}^\top U_k + E_{k,T}^\top E_{k,T}\]
    		since $U_k^\top U_k = I$. Thus we can deduce that
    		\[\specnorm{\tU_{k,T}^\top\tU_{k,T} - I} \le 2\specnorm{E_{k,T}} + \specnorm{E_{k,T}}^2, \]
    		that is, if $\fronorm{E_{k,T}} = \specnorm{\tU_{k,T} - U_k} \le \varepsilon_k$, then $\delta_k \coloneqq \specnorm{\tU_{k,T}^\top\tU_{k,T} - I} \le 2\varepsilon_k + \varepsilon_k^2$. 
            
            Now, for brevity, we define $C_\varepsilon \coloneqq C_{\textrm{pol}}\cdot C_\delta$ so that $\varepsilon_{\max}(T) = \specnorm{\tU_{k,T} - U_k} \le C_\varepsilon e_0^{3^T}$. Since the oracle factor is $(1 - \varepsilon_{\max}(T))^2/(1+\delta_{\max}(T))$, its first-order approximation is 
            \[\frac{(1-\varepsilon_{\max}(T))^2}{1+\delta_{\max}(T)} \approx 1 - (2\varepsilon_{\max}(T) + \delta_{\max}(T)). \]
            Note that $\varepsilon_{\max}(T)$ reduces the strength of descent, while $\delta_{\max}(T)$ weakens the resulting orthogonality. Both of them must be $\ll1$ to preserve fast convergence. To stay within $1-\eta$ ($\eta\in\left[0,1\right)$) of the exact rate, we need 
            \[2\varepsilon_{\max}(T) + \delta_{\max}(T) = (2C_\varepsilon + C_\delta)e_0^{3^T} \le \eta. \]
            Solving for $T$ yields 
            \[T \ge \left\lceil \frac{1}{\log 3}\log\left(\frac{\log((2C_\varepsilon + C_\delta)/\eta)}{\log(1/e_0)}\right) \right\rceil, \]
            since $\log e_0 < 0$, leading to the required number of inner steps. 
        \end{proof}

        \begin{proof}[Proof of \Cref{thm:QDWH}]
            Let us recall that the QDWH algorithm (\Cref{alg:QWDH}) has an equivalent update (the DWH iteration (3.3) in \citep{nakatsukasa2010optimizing}) as follows:
            \[\tU_{k, j+1} = \tU_{k,j}R_{k,j}, \quad R_{k,j} \coloneqq (a_j I + b_j M_{k,j})(I + c_jM_{k,j})^{-1}, \quad M_{k,j} \coloneqq \tU_{k,j}^\top \tU_{k,j}, \]
            with $\tU_{k,0}= G_k/\specnorm{G_k}$ and scalars $a_j$, $b_j$, $c_j>0$ chosen dynamically from a lower bound $\ell_j$ on the smallest singular values to optimize convergence. 
            
            For now, we drop the dependence on $k$ for notational simplicity. Let us define the orthogonal defect of $\tU_j$ by $E_j \coloneqq M_j - I$ and $e_j \coloneqq \specnorm{E_j}= \max_{i\in\setn} |\lambda_j^{(i)}-1|$, where $\lambda_j^{(i)}$ is the $i$th eigenvalue of $M_j$ (in descending order). Since $M_j = \tU_j^\top \tU_j$ is symmetric positive semidefinite and $R_j$ is a rational function of $M_j$, $R_j$ commutes with $M_j$ and is symmetric. Therefore, we have 
            \begin{equation}\label{eqn:thm_qdwh_1}
            	M_{j+1} = \tU_{j+1}^\top \tU_{j+1} = R_j^\top M_j R_j = R_jM_jR_j = M_jR_j^2. 
            \end{equation}
            Next, we have the eigendecomposition of $M_j$ as $M_j = Q_j\Lambda_j Q_j^\top$ with $\Lambda_j = \Diag((\lambda_j^{(i)})_{1\le i\le n})$ and $Q_j\in\OO^{n\times n}$, where $\lambda_j^{(i)} \coloneqq \sigma_i(\tU_j)^2 \in[\ell^2, 1]$. By the definition of $R_j$, we have         
            \[R_j = Q_j r_j(\Lambda_j)Q_j^\top, \quad r_j(t) \coloneqq \frac{a_j + b_jt}{1+c_j}, \]
            where $r_j(\Lambda_j)$ is understood as $\Diag((r_j(\lambda_j^{(i)}))_{1\le i\le n})$. 
            Then, by \eqref{eqn:thm_qdwh_1}, we have 
            \[M_{j+1} = Q_j \varphi_j(\Lambda_j)Q_j^\top, \quad \varphi_j(t) \coloneqq t r_j(t)^2 = t \left(\frac{a_j + b_jt}{1+c_jt}\right)^{\negthickspace2}, \]
            where $\varphi_j(\Lambda_j)$ is understood as $\Diag((\varphi_j(\lambda_j^{(i)}))_{1\le i\le n})$. Thus, we have the following recursive relation of the eigenvalues of $M_j$'s:         
            $\lambda_{j+1}^{(i)} = \varphi_j(\lambda_j^{(i)})$ for all $i\in\setn$. 
            The orthogonality defect of $\tU_{j+1}$ is therefore
            \[e_{j+1} = \max_{i\in\setn} |\varphi_j(\lambda_j^{(i)}) - 1|\quad\text{since}\quad e_j = \max_{i\in\setn} |\lambda_j^{(i)}-1|.  \]
            
            QDWH chooses positive weighting parameters $a_j$, $b_j$ and $c_j$ dynamically such that the fixed point is preserved, i.e., $\varphi_j(1) =1$ for every $a_j>0$, as well as $b_j = (a_j - 1)^2/4$ and $c_j = a_j + b_j -1$. Note that we can derive fixed weights $(a_j, b_j, c_j) = (3,1,3)$ by further imposing $\varphi_j'(1)=\varphi_j''(1)=0$.         
            For any eigenvalue $t\in[\ell_j^2,1]$ and write $t = 1+\Delta$ with $|\Delta|\le e_j$. The Taylor expansion of $\varphi_j$ at $1$ is given by
            \[\varphi_j(1+\Delta) - 1 = \varphi_j'(1) \Delta + \frac12\varphi_j''(1)\Delta^2 + \frac16\varphi_j^{(3)}(\xi)\Delta^3, \]
            for some $\xi$ between $1$ and $1+\Delta$, which yields        
            \[|\varphi_j(t) - 1| \le |\varphi_j'(1)||\Delta| + \frac12|\varphi_j''(1)| |\Delta|^2 +  \frac16\sup_{s\in[\ell_j^2,1]} |\varphi_j^{(3)}(s)||\Delta|^3. \]
            Taking maximum over $|\Delta|\le e_j$ gives 
            \begin{equation}\label{eqn:thm_qdwh_2}
            	e_{j+1} \le |\varphi_j'(1)|e_j + \frac12|\varphi_j''(1)| e_j^2 + C_{3,j} |\Delta|^3, 
            \end{equation}
            where $C_{3,j} \coloneqq \frac16\sup_{s\in[\ell_j^2,1]} |\varphi_j^{(3)}(s)|$ is finite as long as $1+c_js$ is bounded away from $0$, which indeed holds in QDWH since $a_j, b_j, c_j > 0$ and $s\ge0$. 
            
            Now, we show that the cubic term dominates the linear and quadratic terms. Under the dynamic-weighting constraints $b_j = (a_j - 1)^2/4$ and $c_j = a_j + b_j -1$, we can derive that
            \[\varphi_j'(1) = \left(\frac{a_j-3}{a_j+1}\right)^{\negthickspace2}, \quad \varphi_j''(1) = \frac{32(a_j-3)(a_j-1)}{(a_j+1)^4}. \]
            Also let us recall from the definition of $a_j$ that 
            \[a_j = h(\ell_j), \quad h(\ell) = \sqrt{1+\gamma} + \frac12\sqrt{8- 4\gamma + \frac{8(2-\ell^2)}{\ell^2\sqrt{1+\gamma}}}, \quad \gamma = \sqrt[3]{\frac{4(1-\ell^2)}{\ell^4}}, \]
            which is a smooth function of the current lower bound $\ell_j$ and $\ell_j\to1$. As $\ell\to1$, $a=h(\ell)=3+\scrO(1-\ell)$. Since the spectrum $\sigma(M_j)\subseteq[\ell_j^2,1]$, we have $e_j = \specnorm{M_j - I} = \max\{1-\ell_j^2, 1-1\} = 1-\ell_j^2$, which implies $1-\ell_j = (1-\ell_j^2)/(1+\ell_j^2) \le e_j$. Hence, for $j$ large enough (once $\ell_j$ is close to $1$), we have $|a_j-3| \le C_a(1-\ell_j) \le C_a e_j$ for some bounded constant $C_a>0$.                 
            Then, the linear coefficient is upper bounded by 
            \[|\varphi_j'(1)| = \left(\frac{a_j-3}{a_j+1}\right)^{\negthickspace2} \le C(a_j - 3)^2 \le Ce_j^2,\]
            while the quadratic coefficient is upper bounded by 
            \[|\varphi_j''(1)| \le C|a_j-3| \le Ce_j, \]
            for some bounded constant $C>0$.         
            Plugging these into \eqref{eqn:thm_qdwh_2} yields 
            \[e_{j+1} \le (Ce_j^2)e_j + (Ce_j)e_j^2 + C_{3,j}e_j^3 \le \zeta e_j^3, \]
            for sufficiently large $j$ and some bounded constant $\zeta>0$. 
            The remaining part of this proof is similar to that of \Cref{thm:NS} with $e_0 = 1 - \ell_0^2$ and $\ell_0 = \sigma_{\min}(G_k)/\sigma_{\max}(G_k) = 1/\kappa_2(G_k)$, and is thus omitted. 
        \end{proof}

        \section{Numerical Experiments}    
        \label{sec:expt}
        We compare various \PolarGrad optimizers with \textsc{Adam}(W) and \Muon, and study the effect of different numerical polar decomposition algorithms. 
        For more comprehensive understanding of \Muon and more generally \PolarGrad optimizers for different types of matrix optimization problems, including the use of deterministic and stochastic gradients, convexity of the problem (strongly convex, convex and nonconvex problems), and different applications including traditional statistical learning problems and language model pre-training, we include a number of numerical experiments in this section.     
        We start with (strongly) convex problems including a matrix quadratic regression and a matrix logistic regression, followed by a nonconvex low-rank matrix completion problem with simulated data. We then perform Qwen2.5 and GPT-2 Small pre-training experiment. Details of the experiments are given in \Cref{sec:details}. Additional numerical experiments on GPT-2 Medium pre-training are given in \Cref{sec:add_expt}. Open-source implementation of \PolarGrad is available at \url{https://github.com/timlautk/polargrad}.
        
        \subsection{Matrix Quadratic Regression}
        \label{subsec:mat_quad_reg}    
        To better understand the similarities and differences between curvature- and gradient-anisotropy preconditioning, we revisit \Cref{example:quad} numerically, i.e., the quadratic regression objective $\sff(X) = \frac12\fronorm{AXB-C}^2$, where $X\in\RR^{m\times n}$, $A\in\RR^{p\times m}$, $B\in\RR^{n\times q}$ and $C\in\RR^{p\times q}$. We set $(m, n, p, q) = (500, 100, 1000, 250)$ so that $\sff$ is strongly convex. 
        \begin{figure}[h]
              	\centering
              	\includegraphics[width=\textwidth]{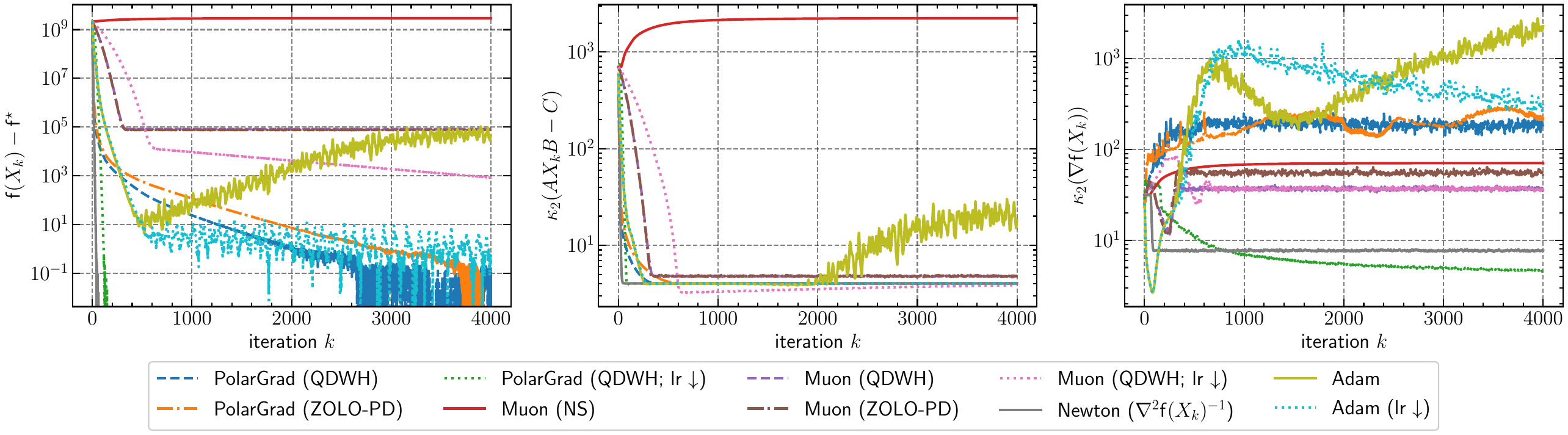}
              	\caption{Losses, residuals and gradient condition numbers of matrix quadratic regression. }
              	\label{fig:mat_quad_reg}
        \end{figure}    
        
        From \Cref{fig:mat_quad_reg}, we make several important observations: (i) The use of better numerical polar decomposition algorithms improves \Muon; (ii) \PolarGrad enjoys much faster early convergence than \Muon (with QDWH and ZOLO-PD) and \Adam, even comparable to Newton's method which enjoys local quadratic convergence for strongly convex functions with Lipschitz Hessian, and empirically verifying the difference between the convergence rates of \PolarGrad (linear; cf.~\Cref{thm:polargrad_strcvx}) and \Muon (sublinear; cf.~\Cref{thm:mat_signSGD_strcvx}); (iii) Learning rate decay of \Muon with deterministic gradients is necessary for convergence to the global minimum even for strongly convex problems (cf.~\Cref{thm:mat_signSGD_strcvx}); (iv) The condition number of the residual $\kappa_2(E_k)$ is indicative of the convergence behavior of optimizers as mentioned in \Cref{example:quad}; (v) Unlike other optimizers, the gradient condition number $\kappa_2(\nabla\sff(X_k))$ in \Adam grows rapidly throughout training, which could be a potential cause for training instabilities. We remark that the intrinsic reason that the optimality gap of \Muon ceases to descend and plateaus at a floor is its failure to satisfy \emph{null-gradient consistency} (\Cref{def:null-grad_consistency}). 
        
        We also plot the gradient nuclear norms to evaluate the difference between \PolarGrad and \Muon. 
        \begin{figure}[h!]
           	\centering
           	\includegraphics[width=\textwidth]{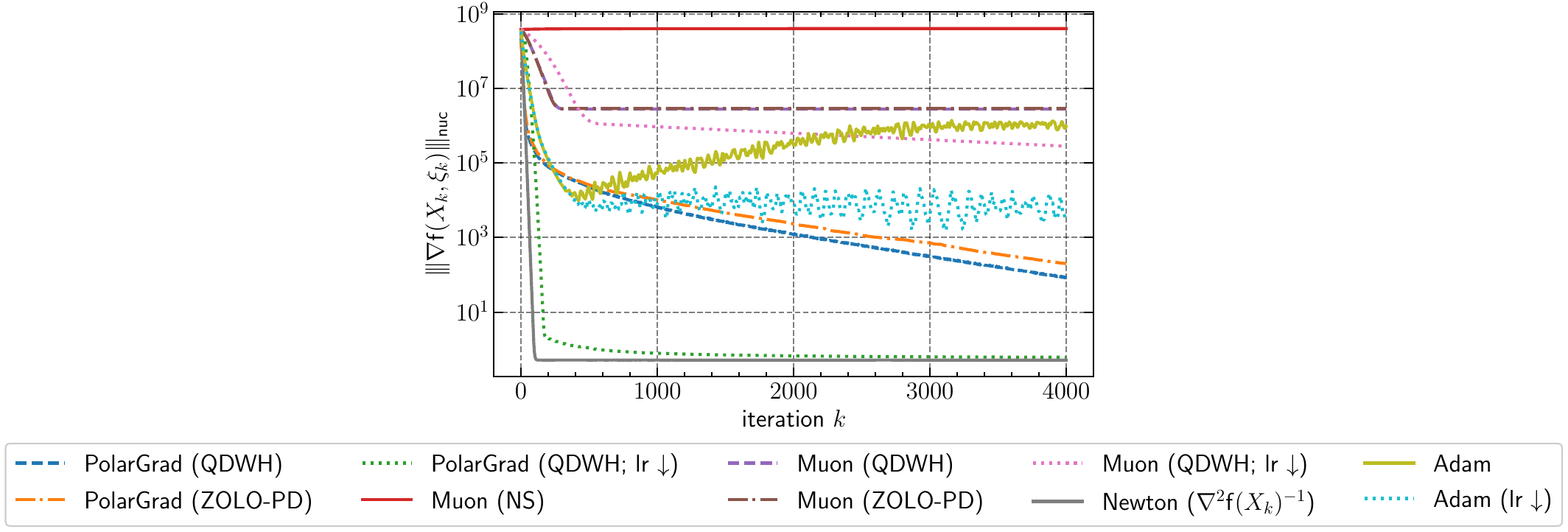}
           	\caption{Gradient nuclear norms of matrix quadratic regression (1st seed).}
           	\label{fig:mat_quad_reg_nuc_1}
        \end{figure}
        
        For this strongly convex problem, \Cref{fig:mat_quad_reg_nuc_1} reveals that the evolution of the gradient nuclear norms is also indicative of the loss convergence of different optimizers. 
        
        \subsection{Matrix Logistic Regression}
        We study a matrix logistic regression problem with the objective $\sff(X) = \sum_{i=1}^N\log(1+\exp(-c_i\odot (a_iXB)))$, where $X\in\RR^{m\times n}$, $A\in\RR^{N\times m}$, $B\in\RR^{n\times q}$ and $C\in\RR^{N\times q}$, and $a_i\in\RR^{1\times m}$ and $c_i\in\RR^{1\times q}$ are the row vectors of $A$ and $C$, respectively.  We set $(m, n, N, q) = (1000, 100, 10000, 400)$. We use minibatch gradients with a batch size of $1000$, sampling with replacement. 
        
        \begin{figure}[h!]
           	\centering
           	\includegraphics[width=\textwidth]{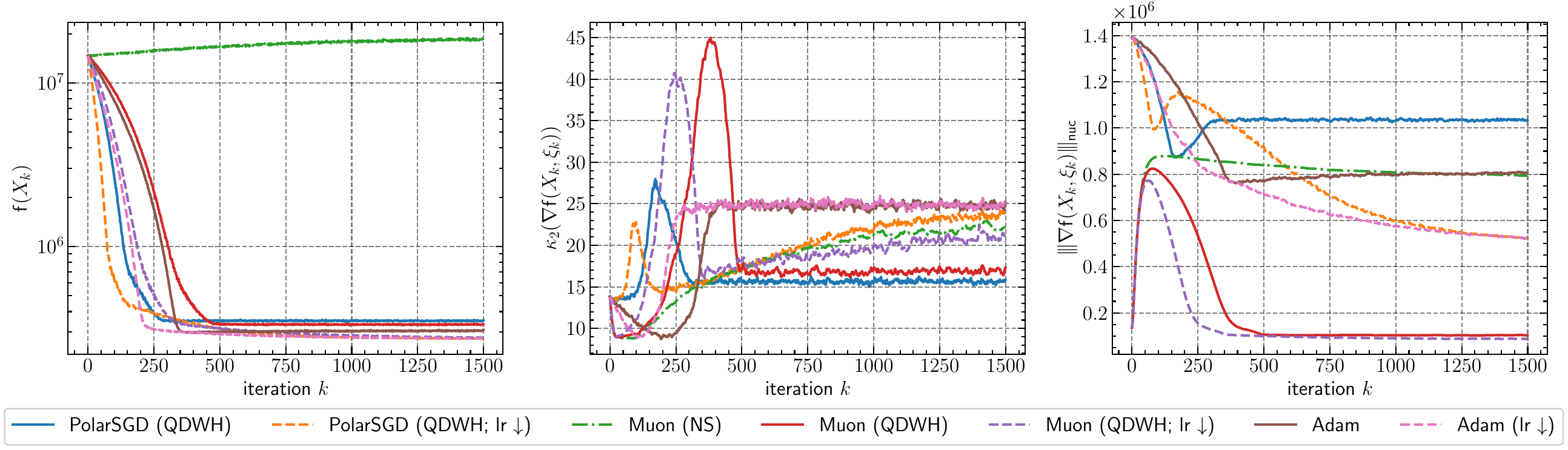}
           	\caption{Losses, gradient condition numbers and nuclear norms of matrix logistic regression.}
           	\label{fig:mat_logistic_reg}
        \end{figure}
        
        From \Cref{fig:mat_logistic_reg}, we also make the following observations: (i) \PolarSGD again enjoys faster early convergence than \Muon and \Adam with constant learning rates; (ii) Learning rate decay is also necessary for all considered optimizers with stochastic gradients even for (strongly) convex problems; (iii) Early loss convergence corresponds to early gradient condition number convergence; (iv) Recall that the nuclear norm of the stochastic gradient $\nucnorm{\nabla\sff(X_k, \xi_k)}$ is the main difference between \PolarSGD and \Muon as a scaling factor of the learning rate---a warmup-then-decay variation can be seen for \Muon with constant learning rate, suggesting that the popular warmup-then-decay learning rate schedule could be used to compensate for the omission of the dual norm scaling factor (cf.~gradient $\ell_1$-norm for \Adam or \signSGD). 
        
        \subsection{Low-Rank Matrix Completion}
        \label{subsec:mat_fac}    
\begin{figure}[h]
        	\centering
        	\includegraphics[width=\textwidth]{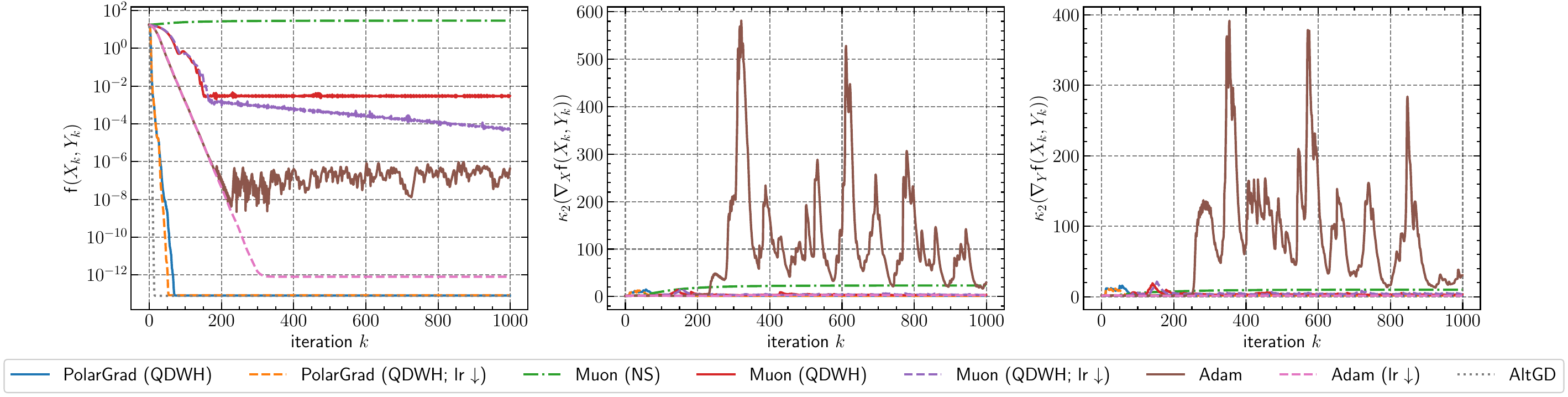}
        	\caption{Losses and gradient condition numbers of low-rank matrix completion.}
        	\label{fig:low_rank_mat_comp}
        \end{figure}

        We first study a simple nonconvex low-rank matrix completion problem with a mask $\calA =(a_{i,j})_{1\le i\le m, 1\le j\le n}\in\RR^{m\times n}$ to mimic missing entries (see e.g., Section IV.C of \citep{chi2019nonconvex}). This model can be viewed as a very simplified neural network. The objective function is $\sff(X, Y) =\fronorm{\calA\odot(XY^\top - M_\star)}^2/\fronorm{\calA}^2$, where $X\in\RR^{m\times r}$, $Y\in\RR^{n\times r}$. We choose $(m, n, r)=(500, 250, 5)$. We also consider alternating gradient descent (\AltGD) method which alternates between solving the two subproblems of $X$ and $Y$ \citep{chi2019nonconvex}. From \Cref{fig:low_rank_mat_comp}, we observe that (i) \PolarGrad has fast and stable convergence compared to \Muon and \Adam, despite converging slower than \AltGD; (ii) the convergence of \Muon plateaus even with learning rate decay, likely due to the omission of the nuclear norm scaling term; (iii) the gradient condition numbers of \Adam is highly unstable unless learning rate decay is used, which is another piece of empirical evidence of the training instabilities of \Adam for nonconvex problems due to poor gradient-anisotropy preconditioning.
        
        \subsection{Qwen2.5 Pre-Training}
        \label{subsec:qwen}
        Our goals here are to understand the general applicability of polar gradient methods including \Muon and \PolarSGDM as matrix optimizers for all matrix parameters in language models including the embedding and head weight matrices in place of \textsc{Adam}(W) and its potential benefits, as well as the potential improvement of \Muon with better numerical polar decomposition algorithms. We keep the use of \AdamW for vector and scalar parameters.     
        We pre-train a modified version of Qwen2.5 \citep{qwen2025qwen_full} with 12 hidden layers and untied embedding on the OpenWebText-100k dataset for one epoch. We plot the training losses and gradient condition numbers of the embedding and head weight matrices. 
        \begin{figure}[h]
        	\centering
        	\includegraphics[width=\textwidth]{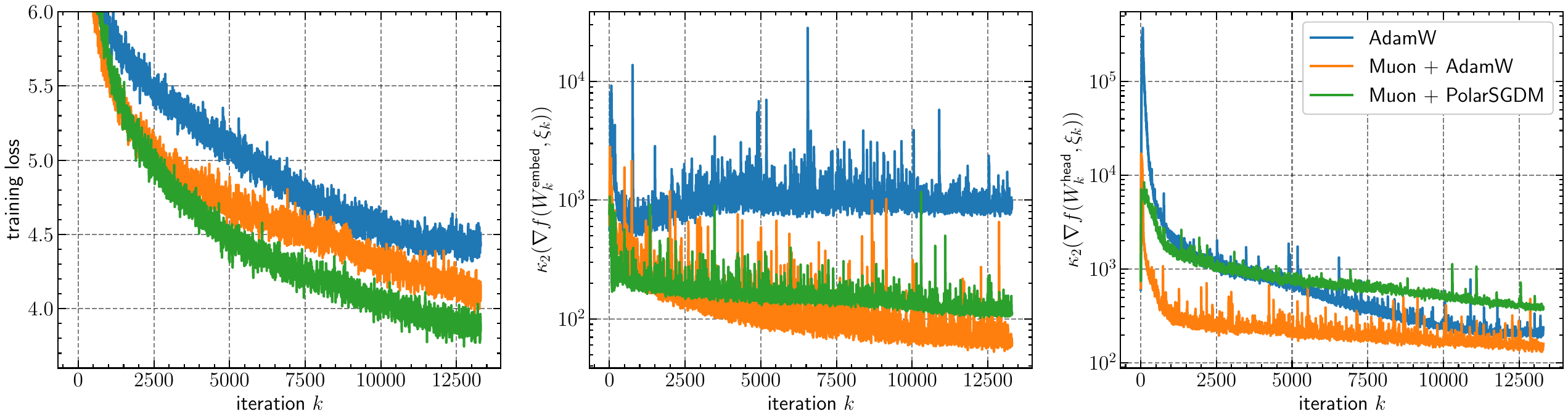}
        	\caption{Training losses and gradient condition numbers of Qwen2.5 pre-training: \AdamW---\AdamW for all parameters; \Muon $+$ \AdamW (\PolarSGDM)---\Muon for hidden layers and \AdamW (\PolarSGDM) for embedding and head layers. }
        	\label{fig:qwen_1}
        \end{figure}
        
        While it is widely agreed that \Muon converges faster when it is applied to matrix parameters in the hidden layers, \AdamW is still used for the embedding and head layers. 
        In \Cref{remark:optim_embed}, we provide an explanation for this choice and that \PolarGrad can still be used for such layers with proper numerical polar decomposition algorithms. 
        From \Cref{fig:qwen_1}, we observe that using \PolarSGDM for these two layers is able to further accelerate convergence. We also observe that there are various large spikes in the gradient condition number of the embedding layer for \AdamW, which could indicate training instability using \AdamW for such ``fat'' matrices. Besides, the current implementation of \Muon relies on the NS iteration, which might not be numerically stable for ill-conditioned matrices, thus hindering its applicability for such matrices; see \Cref{remark:NS_QDWH,subsubsec:stability} for details. 
        
        \subsection{GPT-2 Small 124M Pre-Training}
        With a primary purpose of speedrunning, the \texttt{modded-nanogpt} repository \citep{modded_nanogpt_2024} focuses on pre-training GPT-2 models \citep{radford2019language} on the FineWeb dataset \citep{penedo2024fineweb}, achieving a validation loss at $3.28$ using the least amount of training time. As a result, there are a lot of different aspects of optimization involved in the codebase including implementation and architecture. 
        
\begin{figure}[h!]
           	\centering
           	\includegraphics[width=\textwidth]{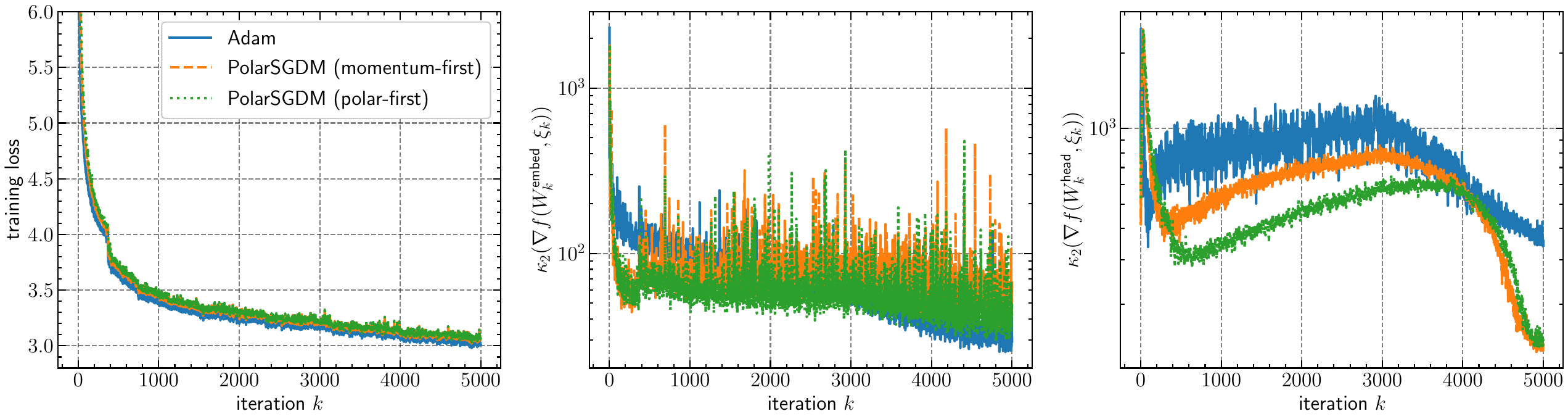}
           	\caption{Training losses and gradient condition numbers of GPT-2 Small 124M pre-training.}
           	\label{fig:gpt2_small}
        \end{figure}

\begin{figure}[h!]
           	\centering
           	\includegraphics[width=\textwidth]{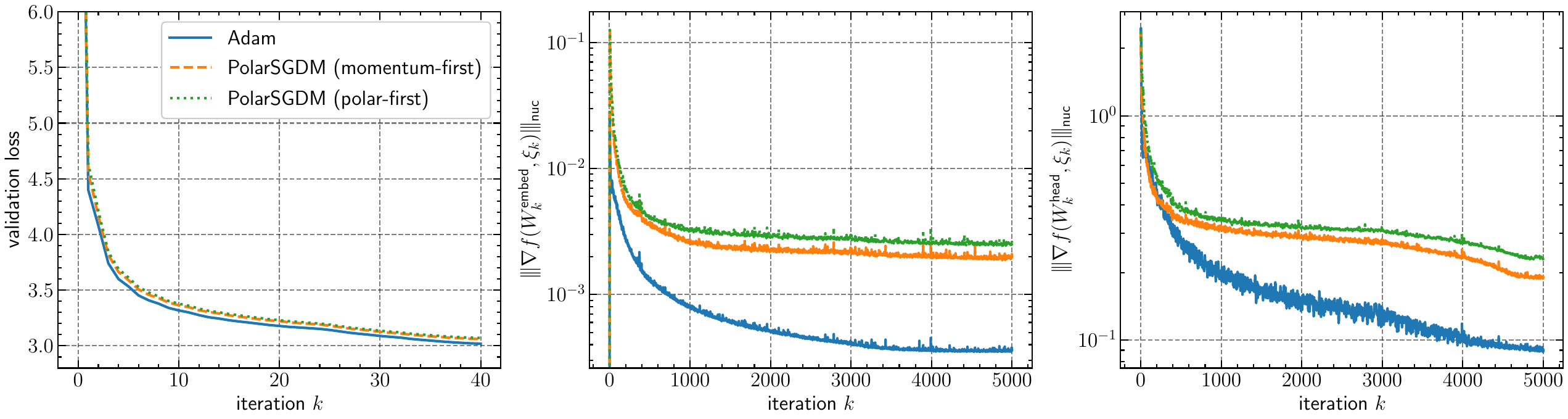}
           	\caption{Validation losses and gradient nuclear norms of GPT-2 Small 124M pre-training.}
           	\label{fig:gpt2_small_2}
        \end{figure}

        Instead, the goal of the experiments on GPT-2 Small and Medium here is to only explore the effect of optimizers without over-optimizing other components of language model development. We hence make use of the setting of the 01/04/25 record. 
        Since this implementation is quite optimized for \Muon for the hidden layers, we keep the use of \Muon for them and just vary the use of \Adam or \PolarSGDM for the embedding and head layers. We also use \Adam for scalar and vector parameters. The implementation of \PolarSGDM is based on the QDWH algorithm.     
        We also compare both types of EMA momentum and plot the results in \Cref{fig:gpt2_small,fig:gpt2_small_2}, with the goal of understanding the different behavior of \PolarSGDM and \Adam for training the embedding and head layers. 
                
        We observe from the plots that while both \Adam and \PolarSGDM have similar training loss curves, it is not the case for the gradient condition numbers and the gradient nuclear norms for the embedding and head layers. The gradient condition numbers of the embedding layer for both \Adam and \PolarSGDM appear to be both noisy despite being convergent, but the gradient condition number of the head layer looks much noisier when \Adam is used. The gradient condition number of the head layer with \PolarSGDM drops rapidly at the last 1000 iterations, which aligns with the interpretation of (better) gradient-anisotropy preconditioning for polar gradient methods. The distinction between momentum-first and polar-first \PolarSGDM is not obvious in this set of experiments.

        \section{Discussion}
        In this work, we establish a unifying preconditioning view for interpreting most deep learning optimizers. Through all these viewpoints and arguments, we compare, contrast and connect most of them under the same umbrella. These optimizers consist of three notable distinctions: (i) types of preconditioning---addressing curvature vs.~gradient anisotropy; (ii) algebraic structures---vectors vs.~matrices, leading to different norms for steepest descent; and (iii) forms of preconditioners---explicit (memory-bound) vs.~implicit (compute-bound). 
        We emphasize the importance of these principles when developing deep learning optimizers and their connections through a unifying preconditioing lens, which is currently missing in the literature. These enhance our understanding of the similarities and differences of these optimizers in a more principled way and pave the road for more efficient and scalable optimizers for large-scale training.     
        Motivated by these principles, we introduce the class of polar gradient methods, as both deep learning optimizers and standalone matrix preconditioned optimization methods which could be of independent interest. Despite their similarities to \Muon in terms of algorithmic designs, our proposed optimizers possess two striking differences---the nuclear norm term and the application of better numerical polar decomposition algorithms. We expect that our proposed optimizers applied to matrix parameters in neural networks are able to mitigate training instability issues arising from \textsc{Adam}(W), hence avoiding the need for instability mitigation tricks such as learning rate warmup. 
        Regarding future work, we plan to develop a more efficient distributed implementation of \PolarGrad similar to that of \citep{shi2023distributed} for \Shampoo and polar decomposition algorithms \citep{lewis2022large}, hence enabling model training on an even larger scale. We also aim to perform in-depth studies on hyperparameter scaling and transfer of \PolarGrad \citep{yang2021tuning,essentialai2025practical_full}, as well as more numerical experiments on other families of models including multi-modal models and MoEs.

        \section*{Acknowledgments}
        The authors would like to thank Damek Davis and Antonio Silveti-Falls for helpful discussion. This work was supported in part by NIH grant U01CA274576, ARPA-H Award D24AC00253, NSF grant DMS-2310679, a Meta Faculty Research Award, and Wharton AI for Business. This work was also supported in part through the computational resources provided by Prime Intellect.

        {\small
        \addcontentsline{toc}{section}{\protect\textbf{References}}
        \bibliographystyle{abbrvnat}
        \begin{CJK*}{UTF8}{gbsn}
        \bibliography{ref}

@article{robbins1951,
	title = {A stochastic approximation method},
	author = {Herbert Robbins and Sutton Monro},
	journal = {The Annals of Mathematical Statistics},
	volume = {22},
	number = {3},
	pages = {400--407},
	year = {1951}
}

@book{beck2017,
	author = {Beck, Amir},
	title = {First-Order Methods in Optimization},
	publisher = {Society for Industrial and Applied Mathematics (SIAM)},
	year = {2017},
	address = {Philadelphia, PA},
}

@article{bottou2018optimization,
  title={Optimization methods for large-scale machine learning},
  author={Bottou, L{\'e}on and Curtis, Frank E. and Nocedal, Jorge},
  journal={SIAM Review},
  volume={60},
  number={2},
  pages={223--311},
  year={2018},
  publisher={SIAM}
}

@book{bauschke2017,
	author="Bauschke, Heinz H.
	and Combettes, Patrick L.",
	title="Convex Analysis and Monotone Operator Theory in Hilbert Spaces",
	year="2017",
	publisher="Springer",
	edition={2nd}
}

@book{rockafellar1970,
	title={Convex Analysis},
	author={Rockafellar, Ralph Tyrell},
	year={1970},
	publisher={Princeton University Press},
	address = {Princeton, NJ},
}

@Book{rockafellar1998,
	Title  = {Variational Analysis},
	Author = {Rockafellar, Ralph Tyrell and Wets, Roger J.-B.},
	Publisher = {Springer},
	Year = {1998},
}

@inproceedings{kingma2015,
	title = {Adam: a method for stochastic optimization},
	author = {Kingma, Diederik P. and Ba, Jimmy Lei},
	booktitle={International Conference on Learning Representations (ICLR)},
	year={2015},
}

@article{duchi2011adagrad,
	author = {Duchi, John and Hazan, Elad and Singer, Yoram},
	title = {Adaptive Subgradient Methods for Online Learning and Stochastic Optimization},
	journal = {Journal of Machine Learning Research},
	volume = {12},
	year = {2011},
	pages = {2121--2159},
}

@article{zeiler2012adadelta,
	author    = {Matthew D. Zeiler},
	title     = {{ADADELTA}: An Adaptive Learning Rate Method},
	journal   = {arXiv preprint arXiv:1212.5701},
	year      = {2012},
}

@misc{tieleman2012,
	title={{Lecture 6.5---RMSProp: Divide the gradient by a running average of its recent magnitude}},
	author={Tieleman, Tijmen and Hinton, Geoffrey},
	howpublished={Coursera: Neural Networks for Machine Learning},
	year={2012}
}

@article{combettes2021fixed,
	title={Fixed point strategies in data science},
	author={Combettes, Patrick L. and Pesquet, Jean-Christophe},
	journal={IEEE Transactions on Signal Processing},
	volume={69},
    pages={3878--3905},
    year={2021},
}

@article{teboulle2018simplified,
	title={A simplified view of first order methods for optimization},
	author={Teboulle, Marc},
	journal={Mathematical Programming},
	volume={170},
	number={1},
	pages={67--96},
	year={2018},
	publisher={Springer}
}

@article{beck2003mirror,
  title={Mirror descent and nonlinear projected subgradient methods for convex optimization},
  author={Beck, Amir and Teboulle, Marc},
  journal={Operations Research Letters},
  volume={31},
  number={3},
  pages={167--175},
  year={2003},
  publisher={Elsevier}
}

@inproceedings{chen2023symbolic,
	title={Symbolic discovery of optimization algorithms},
	author={Chen, Xiangning and Liang, Chen and Huang, Da and Real, Esteban and Wang, Kaiyuan and Liu, Yao and Pham, Hieu and Dong, Xuanyi and Luong, Thang and Hsieh, Cho-Jui and Lu, Yifeng and Le, Quoc V.},
	booktitle={Advances in Neural Information Processing Systems (NeurIPS)},
	year={2023}
}

@inproceedings{liu2024sophia,
	title={Sophia: A Scalable Stochastic Second-order Optimizer for Language Model Pre-training},
	author={Hong Liu and Zhiyuan Li and David Hall and Percy Liang and Tengyu Ma},
	booktitle={International Conference on Learning Representations (ICLR)},
	year={2024}
}

@inproceedings{hoffmann2022training,
	title={Training Compute-Optimal Large Language Models},
	author={Jordan Hoffmann and Sebastian Borgeaud and Arthur Mensch and Elena Buchatskaya and Trevor Cai and Eliza Rutherford and Diego de las Casas and Lisa Anne Hendricks and Johannes Welbl and Aidan Clark and Tom Hennigan and Eric Noland and Katherine Millican and George van den Driessche and Bogdan Damoc and Aurelia Guy and Simon Osindero and Karen Simonyan and Erich Elsen and Oriol Vinyals and Jack William Rae and Laurent Sifre},
	booktitle={Advances in Neural Information Processing Systems (NeurIPS)},
	year={2022}
}

@article{radford2019language,
	title={Language models are unsupervised multitask learners},
	author={Radford, Alec and Wu, Jeffrey and Child, Rewon and Luan, David and Amodei, Dario and Sutskever, Ilya},
	journal={OpenAI blog},
	year={2019}
}

@article{kaplan2020scaling,
	title={Scaling laws for neural language models},
	author={Kaplan, Jared and McCandlish, Sam and Henighan, Tom and Brown, Tom B. and Chess, Benjamin and Child, Rewon and Gray, Scott and Radford, Alec and Wu, Jeffrey and Amodei, Dario},
	journal={arXiv preprint arXiv:2001.08361},
	year={2020}
}

@article{dahl2023benchmarking,
	title={Benchmarking Neural Network Training Algorithms},
	author={George E. Dahl and Frank Schneider and Zachary Nado and Naman Agarwal and Chandramouli Shama Sastry and Philipp Hennig and Sourabh Medapati and Runa Eschenhagen and Priya Kasimbeg and Daniel Suo and Juhan Bae and Justin Gilmer and Abel L. Peirson and Bilal Khan and Rohan Anil and Mike Rabbat and Shankar Krishnan and Daniel Snider and Ehsan Amid and Kongtao Chen and Chris J. Maddison and Rakshith Vasudev and Michal Badura and Ankush Garg and Peter Mattson},
	journal={arXiv preprint arXiv:2306.07179},
	year={2023}
}

@inproceedings{loshchilov2019decoupled,
	title={Decoupled Weight Decay Regularization},
	author={Loshchilov, Ilya and Hutter, Frank},
	booktitle={International Conference on Learning Representations (ICLR)},
	year={2019}
}

@inproceedings{paszke2019pytorch,
    title={{PyTorch}: An imperative style, high-performance deep learning library},
    author={Adam Paszke and Sam Gross and Francisco Massa and Adam Lerer and James Bradbury and Gregory Chanan and Trevor Killeen and Zeming Lin and Natalia Gimelshein and Luca Antiga and Alban Desmaison and Andreas Köpf and Edward Yang and Zach DeVito and Martin Raison and Alykhan Tejani and Sasank Chilamkurthy and Benoit Steiner and Lu Fang and Junjie Bai and Soumith Chintala},
    booktitle={Advances in Neural Information Processing Systems (NeurIPS)},
    year={2019}
}

@misc{jax2018github,
    author = {James Bradbury and Roy Frostig and Peter Hawkins and Matthew James Johnson and Chris Leary and Dougal Maclaurin and George Necula and Adam Paszke and Jake Vander{P}las and Skye Wanderman-{M}ilne and Qiao Zhang},
    title = {{JAX}: composable transformations of {P}ython+{N}um{P}y programs},
    url = {http://github.com/google/jax},
    version = {0.2.5},
    year = {2018},
}

@inproceedings{sutskever2013importance,
    title={On the importance of initialization and momentum in deep learning},
    author={Sutskever, Ilya and Martens, James and Dahl, George and Hinton, Geoffrey},
    booktitle = {Proceedings of the International Conference on Machine Learning (ICML)},
    year={2013},
}

@article{shi2023distributed,
    title={A Distributed Data-Parallel {PyTorch} Implementation of the Distributed {Shampoo} Optimizer for Training Neural Networks At-Scale},
    author={Shi, Hao-Jun Michael and Lee, Tsung-Hsien and Iwasaki, Shintaro and Gallego-Posada, Jose and Li, Zhijing and Rangadurai, Kaushik and Mudigere, Dheevatsa and Rabbat, Michael},
    journal={arXiv preprint arXiv:2309.06497},
    year={2023}
}

@inproceedings{mcmahan2010adaptive,
    title={Adaptive bound optimization for online convex optimization},
    author={McMahan, H. Brendan and Streeter, Matthew},
    booktitle = {Proceedings of the Conference on Learning Theory (COLT)},
    year={2010}
}

@inproceedings{zhang2024transformers,
  title={Why Transformers Need {A}dam: A {H}essian Perspective},
  author={Zhang, Yushun and Chen, Congliang and Ding, Tian and Li, Ziniu and Sun, Ruoyu and Luo, Zhi-Quan},
  booktitle={Advances in Neural Information Processing Systems (NeurIPS)},
  year={2024}
}

@article{kunstner2024heavy,
  title={Heavy-Tailed Class Imbalance and Why {A}dam Outperforms Gradient Descent on Language Models},
  author={Kunstner, Frederik and Yadav, Robin and Milligan, Alan and Schmidt, Mark and Bietti, Alberto},
  journal={arXiv preprint arXiv:2402.19449},
  year={2024}
}

@inproceedings{yang2021tuning,
  title={Tuning Large Neural Networks via Zero-Shot Hyperparameter Transfer},
  author={Greg Yang and Edward J. Hu and Igor Babuschkin and Szymon Sidor and Xiaodong Liu and David Farhi and Nick Ryder and Jakub Pachocki and Weizhu Chen and Jianfeng Gao},
  booktitle={Advances in Neural Information Processing Systems (NeurIPS)},
  year={2021},
}

@inproceedings{penedo2024fineweb,
  title={The {FineWeb} Datasets: Decanting the Web for the Finest Text Data at Scale},
  author={Penedo, Guilherme and Kydl{\'\i}{\v{c}}ek, Hynek and Lozhkov, Anton and Mitchell, Margaret and Raffel, Colin and Von Werra, Leandro and Wolf, Thomas},
  booktitle={Advances in Neural Information Processing Systems (NeurIPS) Datasets and Benchmarks Track},
  year={2024}
}

@article{condat2023proximal,
	title = {Proximal Splitting Algorithms for Convex Optimization: A Tour of Recent Advances, with New Twists},
	author={Condat, Laurent and Kitahara, Daichi and Contreras, Andr{\'e}s and Hirabayashi, Akira},
	journal={SIAM Review},
	volume = {65},
    number = {2},
    pages = {375--435},
    year = {2023},
}

@inproceedings{bernstein2024modular,
  title={Modular Duality in Deep Learning},
  author={Bernstein, Jeremy and Newhouse, Laker},
  booktitle={Proceedings of the International Conference on Machine Learning (ICML)},
  year={2025}
}

@inproceedings{bernstein2024old,
  title={Old Optimizer, New Norm: An Anthology},
  author={Bernstein, Jeremy and Newhouse, Laker},
  booktitle={OPT 2024: Optimization for Machine Learning},
  year={2024}
}

@inproceedings{large2024scalable,
  title={Scalable Optimization in the Modular Norm},
  author={Large, Tim and Liu, Yang and Huh, Minyoung and Bahng, Hyojin and Isola, Phillip and Bernstein, Jeremy},
  booktitle={Advances in Neural Information Processing Systems (NeurIPS)},
  year={2024}
}

@article{lewis1995convex,
  title={The convex analysis of unitarily invariant matrix functions},
  author={Lewis, Adrian S.},
  journal={Journal of Convex Analysis},
  volume={2},
  number={1},
  pages={173--183},
  year={1995}
}

@article{lewis1996eigenvalue,
  title={Eigenvalue Optimization},
  author={Lewis, Adrian S. and Overton, Michael L.},
  journal={Acta Numerica},
  volume={5},
  pages={149--190},
  year={1996},
  publisher={Cambridge University Press}
}

@article{lewis2003mathematics,
  title={The mathematics of eigenvalue optimization},
  author={Lewis, Adrian S.},
  journal={Mathematical Programming},
  volume={97},
  pages={155--176},
  year={2003},
  publisher={Springer}
}

@misc{jordan2024muon,
  author       = {Keller Jordan and Yuchen Jin and Vlado Boza and You Jiacheng and
                  Franz Cecista and Laker Newhouse and Jeremy Bernstein},
  title        = {Muon: An optimizer for hidden layers in neural networks},
  year         = {2024},
  url          = {https://kellerjordan.github.io/posts/muon/}
}

@inproceedings{gupta2018shampoo,
  title={Shampoo: Preconditioned stochastic tensor optimization},
  author={Gupta, Vineet and Koren, Tomer and Singer, Yoram},
  booktitle={Proceedings of the International Conference on Machine Learning (ICML)},
  year={2018},
}

@article{anil2020scalable,
  title={Scalable second order optimization for deep learning},
  author={Anil, Rohan and Gupta, Vineet and Koren, Tomer and Regan, Kevin and Singer, Yoram},
  journal={arXiv preprint arXiv:2002.09018},
  year={2020}
}

@inproceedings{vyas2024soap,
  title={{SOAP}: Improving and stabilizing {Shampoo} using {Adam}},
  author={Vyas, Nikhil and Morwani, Depen and Zhao, Rosie and Shapira, Itai and Brandfonbrener, David and Janson, Lucas and Kakade, Sham},
  booktitle={International Conference on Learning Representations (ICLR)},
  year={2025}
}

@inproceedings{duvvuri2024combining,
  title={Combining Axes Preconditioners through {K}ronecker Approximation for Deep Learning},
  author={Duvvuri, Sai Surya and Devvrit, Fnu and Anil, Rohan and Hsieh, Cho-Jui and Dhillon, Inderjit S.},
  booktitle={International Conference on Learning Representations (ICLR)},
  year={2024}
}

@misc{modded_nanogpt_2024,
  author       = {Keller Jordan and Jeremy Bernstein and Brendan Rappazzo and
                  @fernbear.bsky.social and Boza Vlado and You Jiacheng and
                  Franz Cesista and Braden Koszarsky and @Grad62304977},
  title        = {\texttt{modded-nanogpt}: Speedrunning the {NanoGPT} baseline},
  year         = {2024},
  url          = {https://github.com/KellerJordan/modded-nanogpt},
}

@article{ziketak1988characterization,
  title={On the characterization of the extremal points of the unit sphere of matrices},
  author={Zi\k{e}tak, Krystyna},
  journal={Linear Algebra and Its Applications},
  volume={106},
  pages={57--75},
  year={1988},
  publisher={Elsevier}
}

@article{ziketak1993subdifferentials,
  title={Subdifferentials, faces, and dual matrices},
  author={Zi\k{e}tak, Krystyna},
  journal={Linear Algebra and Its Applications},
  volume={185},
  pages={125--141},
  year={1993},
  publisher={Elsevier}
}

@article{ma2024swan,
      title={{SWAN}: Preprocessing {SGD} Enables {Adam}-Level Performance On {LLM} Training With Significant Memory Reduction}, 
      author={Chao Ma and Wenbo Gong and Meyer Scetbon and Edward Meeds},
      year={2024},
      journal={arXiv preprint arXiv:2412.13148},
}

@article{li2017preconditioned,
  title={Preconditioned stochastic gradient descent},
  author={Li, Xi-Lin},
  journal={IEEE Transactions on Neural Networks and Learning Systems},
  volume={29},
  number={5},
  pages={1454--1466},
  year={2017},
  publisher={IEEE}
}

@inproceedings{morwani2024new,
  title={A New Perspective on {S}hampoo's Preconditioner},
  author={Morwani, Depen and Shapira, Itai and Vyas, Nikhil and Malach, Eran and Kakade, Sham M. and Janson, Lucas},
  booktitle={International Conference on Learning Representations (ICLR)},
  year={2025}
}

@misc{su2024muon,
    title={Appreciating the {Muon} Optimizer: From Vectors to Matrices, An Essential Leap},
    author={Jianlin Su},
    year={2024},
    month={December},
    url={https://kexue.fm/archives/10592},
}

@misc{su2024hessian,
    title={Adaptive Learning Rate Optimizer from the Perspective of {Hessian} Approximation},
    author={Jianlin Su},
    year={2024},
    month={November},
    url={https://kexue.fm/archives/10588},
}

@misc{su2025muon,
    title={Why We Chose {Muon}: Our Chain of Thought},
    author={Jianlin Su},
    year={2025},
    month={February},
    url={https://x.com/Kimi_Moonshot/status/1897929976948965870},
}

@misc{su2025msign,
	title={{N}ewton--{S}chulz Iteration of the msign Operator},
	author={Jianlin Su},
	year={2025},
	month={May},
	url={https://kexue.fm/archives/10922},
}

@inproceedings{bernstein2018signsgd,
  title={sign{SGD}: Compressed optimisation for non-convex problems},
  author={Bernstein, Jeremy and Wang, Yu-Xiang and Azizzadenesheli, Kamyar and Anandkumar, Animashree},
  booktitle={Proceedings of the International Conference on Machine Learning (ICML)},
  year={2018},
}

@article{watson1992characterization,
  title={Characterization of the subdifferential of some matrix norms},
  author={Watson, G. Alistair},
  journal={Linear Algebra and its Applications},
  volume={170},
  number={1},
  pages={33--45},
  year={1992}
}

@inproceedings{carlson2015stochasticRBM,
  title={Stochastic spectral descent for restricted {B}oltzmann machines},
  author={Carlson, David and Cevher, Volkan and Carin, Lawrence},
  booktitle={Proceedings of the International Conference on Artificial Intelligence and Statistics (AISTATS)},
  year={2015},
}

@article{carlson2016stochastic,
  title={Stochastic spectral descent for discrete graphical models},
  author={Carlson, David and Hsieh, Ya-Ping and Collins, Edo and Carin, Lawrence and Cevher, Volkan},
  journal={IEEE Journal of Selected Topics in Signal Processing},
  volume={10},
  number={2},
  pages={296--311},
  year={2016},
  publisher={IEEE}
}

@inproceedings{carlson2015preconditioned,
  title={Preconditioned spectral descent for deep learning},
  author={Carlson, David and Collins, Edo and Hsieh, Ya-Ping and Carin, Lawrence and Cevher, Volkan},
  booktitle={Advances in Neural Information Processing Systems (NeurIPS)},
  year={2015}
}

@inproceedings{shazeer2018adafactor,
  title={Adafactor: Adaptive learning rates with sublinear memory cost},
  author={Shazeer, Noam and Stern, Mitchell},
  booktitle={Proceedings of the International Conference on Machine Learning (ICML)},
  year={2018},
}

@article{flynn2017duality,
  title={The duality structure gradient descent algorithm: analysis and applications to neural networks},
  author={Flynn, Thomas},
  journal={arXiv preprint arXiv:1708.00523},
  year={2017}
}

@article{hsieh2018non,
  title={A non-{E}uclidean gradient descent framework for non-convex matrix factorization},
  author={Hsieh, Ya-Ping and Kao, Yu-Chun and Mahabadi, Rabeeh Karimi and Yurtsever, Alp and Kyrillidis, Anastasios and Cevher, Volkan},
  journal={IEEE Transactions on Signal Processing},
  volume={66},
  number={22},
  pages={5917--5926},
  year={2018},
  publisher={IEEE}
}

@inproceedings{kelner2014almost,
  title={An almost-linear-time algorithm for approximate max flow in undirected graphs, and its multicommodity generalizations},
  author={Kelner, Jonathan A. and Lee, Yin Tat and Orecchia, Lorenzo and Sidford, Aaron},
  booktitle={Proceedings of the ACM-SIAM Symposium on Discrete Algorithms (SODA)},
  year={2014},
}

@article{vonneumann1937some,
    title={Some Matrix-Inequalities and Metrization of Matrix-Space},
    author={von Neumann, John},
    year={1937},
    journal={Tomskii University Review},
    volume={1},
    pages={286--300},
    note = "In: Collected Works, {\rm (A. H. Taub Editor)}, Pergamon, Oxford, 1962, Volume IV, 205--218.",
}

@article{tuddenham2022orthogonalising,
  title={Orthogonalising gradients to speed up neural network optimisation},
  author={Tuddenham, Mark and Pr{\"u}gel-Bennett, Adam and Hare, Jonathan},
  journal={arXiv preprint arXiv:2202.07052},
  year={2022}
}

@article{mirsky1960symmetric,
  title={Symmetric gauge functions and unitarily invariant norms},
  author={Mirsky, Leon},
  journal={The Quarterly Journal of Mathematics},
  volume={11},
  number={1},
  pages={50--59},
  year={1960},
}

@inproceedings{xie2024adam,
  title={Adam Exploits $\ell_\infty$-geometry of Loss Landscape via Coordinate-wise Adaptivity},
  author={Xie, Shuo and Mohamadi, Mohamad Amin and Li, Zhiyuan},
  booktitle={International Conference on Learning Representations (ICLR)},
  year={2025}
}

@inproceedings{xie2024implicit,
  title={Implicit Bias of {AdamW}: $\ell_\infty$-Norm Constrained Optimization},
  author={Xie, Shuo and Li, Zhiyuan},
  booktitle={Proceedings of the International Conference on Machine Learning (ICML)},
  year={2024}
}

@inproceedings{riedmiller1993direct,
  title={A direct adaptive method for faster backpropagation learning: The {RPROP} algorithm},
  author={Riedmiller, Martin and Braun, Heinrich},
  booktitle={Proceedings of the IEEE International Conference on Neural Networks},
  year={1993},
}

@article{chen2021spectral,
  title={Spectral methods for data science: A statistical perspective},
  author={Chen, Yuxin and Chi, Yuejie and Fan, Jianqing and Ma, Cong},
  journal={Foundations and Trends{\textregistered} in Machine Learning},
  volume={14},
  number={5},
  pages={566--806},
  year={2021},
  publisher={Now Publishers, Inc.}
}

@inproceedings{balles2018dissecting,
  title={Dissecting {A}dam: The sign, magnitude and variance of stochastic gradients},
  author={Balles, Lukas and Hennig, Philipp},
  booktitle={Proceedings of the International Conference on Machine Learning (ICML)},
  year={2018},
}

@article{pooladzandi2024curvature,
  title={Curvature-Informed {SGD} via General Purpose {L}ie-Group Preconditioners},
  author={Pooladzandi, Omead and Li, Xi-Lin},
  journal={arXiv preprint arXiv:2402.04553},
  year={2024}
}

@article{nakatsukasa2013stable,
  title={Stable and efficient spectral divide and conquer algorithms for the symmetric eigenvalue decomposition and the {SVD}},
  author={Nakatsukasa, Yuji and Higham, Nicholas J.},
  journal={SIAM Journal on Scientific Computing},
  volume={35},
  number={3},
  pages={A1325--A1349},
  year={2013},
}

@article{nakatsukasa2016computing,
  title={Computing fundamental matrix decompositions accurately via the matrix sign function in two iterations: The power of {Z}olotarev's functions},
  author={Nakatsukasa, Yuji and Freund, Roland W.},
  journal={SIAM Review},
  volume={58},
  number={3},
  pages={461--493},
  year={2016},
}

@article{ltaief2019massively,
  title={Massively parallel polar decomposition on distributed-memory systems},
  author={Ltaief, Hatem and Sukkari, Dalal and Esposito, Aniello and Nakatsukasa, Yuji and Keyes, David},
  journal={ACM Transactions on Parallel Computing (TOPC)},
  volume={6},
  number={1},
  pages={1--15},
  year={2019},
}

@article{benfenati2020proximal,
  title={Proximal approaches for matrix optimization problems: Application to robust precision matrix estimation},
  author = {Benfenati, Alessandro and Chouzenoux, Emilie and Pesquet, Jean-Christophe},
  journal={Signal Processing},
  volume={169},
  pages={107417},
  year={2020},
}

@article{lewis2022large,
  title={Large-scale distributed linear algebra with tensor processing units},
  author={Lewis, Adam G. M. and Beall, Jackson and Ganahl, Martin and Hauru, Markus and Mallick, Shrestha Basu and Vidal, Guifre},
  journal={Proceedings of the National Academy of Sciences},
  volume={119},
  number={33},
  pages={e2122762119},
  year={2022},
}

@article{nakatsukasa2010optimizing,
  title={Optimizing {H}alley's iteration for computing the matrix polar decomposition},
  author={Nakatsukasa, Yuji and Bai, Zhaojun and Gygi, Fran{\c{c}}ois},
  journal={SIAM Journal on Matrix Analysis and Applications},
  volume={31},
  number={5},
  pages={2700--2720},
  year={2010},
  publisher={SIAM}
}

@inproceedings{kasimbeg2025accelerating,
title={Accelerating neural network training: An analysis of the {AlgoPerf} competition},
author={Priya Kasimbeg and Frank Schneider and Runa Eschenhagen and Juhan Bae and Chandramouli Shama Sastry and Mark Saroufim and Boyuan Feng and Less Wright and Edward Z. Yang and Zachary Nado and Sourabh Medapati and Philipp Hennig and Michael Rabbat and George E. Dahl},
booktitle={International Conference on Learning Representations (ICLR)},
year={2025},
}

@book{higham2008functions,
  title={Functions of Matrices: Theory and Computation},
  author={Higham, Nicholas J.},
  year={2008},
  publisher = {Society for Industrial and Applied Mathematics},
}

@article{qu2025optimal,
  title={Optimal Diagonal Preconditioning},
  author={Qu, Zhaonan and Gao, Wenzhi and Hinder, Oliver and Ye, Yinyu and Zhou, Zhengyuan},
  journal={Operations Research},
  volume={73},
  number={3},
  pages={1479--1495},
  year={2025},
  publisher={INFORMS}
}

@article{jambulapati2020fast,
  title={Fast and near-optimal diagonal preconditioning},
  author={Jambulapati, Arun and Li, Jerry and Musco, Christopher and Sidford, Aaron and Tian, Kevin},
  journal={arXiv preprint arXiv:2008.01722},
  year={2020}
}

@inproceedings{martens2015optimizing,
  title={Optimizing neural networks with {K}ronecker-factored approximate curvature},
  author={Martens, James and Grosse, Roger},
  booktitle={Proceedings of the International Conference on Machine Learning (ICML)},
  year={2015},
}

@book{boyd2004convex,
  title={Convex Optimization},
  author={Boyd, Stephen and Vandenberghe, Lieven},
  year={2004},
  publisher={Cambridge University Press}
}

@article{auslender2006interior,
  title={Interior gradient and proximal methods for convex and conic optimization},
  author={Auslender, Alfred and Teboulle, Marc},
  journal={SIAM Journal on Optimization},
  volume={16},
  number={3},
  pages={697--725},
  year={2006},
  publisher={SIAM}
}

@article{turing1948rounding,
  title={Rounding-off errors in matrix processes},
  author={Turing, Alan M.},
  journal={The Quarterly Journal of Mechanics and Applied Mathematics},
  volume={1},
  number={1},
  pages={287--308},
  year={1948},
  publisher={Oxford University Press}
}

@article{bian2024preconditioned,
  title={A preconditioned {R}iemannian gradient descent algorithm for low-rank matrix recovery},
  author={Bian, Fengmiao and Cai, Jian-Feng and Zhang, Rui},
  journal={SIAM Journal on Matrix Analysis and Applications},
  volume={45},
  number={4},
  pages={2075--2103},
  year={2024},
  publisher={SIAM}
}

@article{higham1986computing,
  title={Computing the polar decomposition---with applications},
  author={Higham, Nicholas J.},
  journal={SIAM Journal on Scientific and Statistical Computing},
  volume={7},
  number={4},
  pages={1160--1174},
  year={1986},
  publisher={SIAM}
}

@article{higham1990fast,
  title={Fast polar decomposition of an arbitrary matrix},
  author={Higham, Nicholas J. and Schreiber, Robert S.},
  journal={SIAM Journal on Scientific and Statistical Computing},
  volume={11},
  number={4},
  pages={648--655},
  year={1990},
  publisher={SIAM}
}

@article{higham1994matrix,
	title={The matrix sign decomposition and its relation to the polar decomposition},
	author={Higham, Nicholas J.},
	journal={Linear Algebra and its Applications},
	volume={212},
	pages={3--20},
	year={1994},
	publisher={Elsevier}
}

@article{watson1993matrix,
  title={On matrix approximation problems with {Ky Fan} $k$ norms},
  author={Watson, G. Alistair},
  journal={Numerical Algorithms},
  volume={5},
  pages={263--272},
  year={1993},
  publisher={Springer}
}

@article{lewis1996group,
  title={Group invariance and convex matrix analysis},
  author={Lewis, Adrian S.},
  journal={SIAM Journal on Matrix Analysis and Applications},
  volume={17},
  number={4},
  pages={927--949},
  year={1996},
  publisher={SIAM}
}

@article{watson1991algorithm,
  title={An algorithm for optimal $\ell_2$ scaling of matrices},
  author={Watson, G. Alistair},
  journal={IMA Journal of Numerical Analysis},
  volume={11},
  number={4},
  pages={481--492},
  year={1991},
  publisher={Oxford University Press}
}

@inproceedings{zhuang2020adabelief,
  title={{AdaBelief} optimizer: Adapting stepsizes by the belief in observed gradients},
  author={Zhuang, Juntang and Tang, Tommy and Ding, Yifan and Tatikonda, Sekhar C. and Dvornek, Nicha and Papademetris, Xenophon and Duncan, James},
  booktitle={Advances in Neural Information Processing Systems (NeurIPS)},
  year={2020}
}

@article{chi2019nonconvex,
  title={Nonconvex optimization meets low-rank matrix factorization: An overview},
  author={Chi, Yuejie and Lu, Yue M. and Chen, Yuxin},
  journal={IEEE Transactions on Signal Processing},
  volume={67},
  number={20},
  pages={5239--5269},
  year={2019},
  publisher={IEEE}
}

@book{horn2012matrix,
  title={Matrix Analysis},
  author={Horn, Roger A. and Johnson, Charles R.},
  year={2012},
  edition={2nd},
  publisher={Cambridge University Press}
}

@book{horn1994topics,
  title={Topics in Matrix Analysis},
  author={Horn, Roger A. and Johnson, Charles R.},
  year={1994},
  publisher={Cambridge University Press}
}

@article{nakatsukasa2012backward,
  title={Backward stability of iterations for computing the polar decomposition},
  author={Nakatsukasa, Yuji and Higham, Nicholas J.},
  journal={SIAM Journal on Matrix Analysis and Applications},
  volume={33},
  number={2},
  pages={460--479},
  year={2012},
  publisher={SIAM}
}

@book{mordukhovich2022convex,
  title={Convex Analysis and Beyond. Volume I: Basic Theory},
  author={Mordukhovich, Boris S. and Nam, Nguyen Mau},
  year={2022},
  publisher={Springer},
  series={Springer Series in Operations Research and Financial Engineering}
}

@article{li2025muon,
      title={A Note on the Convergence of {Muon} and Further}, 
      author={Jiaxiang Li and Mingyi Hong},
      year={2025},
      journal={arXiv preprint arXiv:2502.02900},
}

@article{nesterov1983method,
author="Nesterov, Yurii",
title="A method for solving the convex programming problem with convergence rate $o(1/k^2)$",
journal={Doklady Akademii Nauk},
volume={269},
number={3},
pages={543--547},
year={1983},
organization={Russian Academy of Sciences}
}

@InProceedings{dozat2016incorporating,
  title={Incorporating {N}esterov momentum into {A}dam},
  author={Dozat, Timothy},
  booktitle = {International Conference on Learning Representations (ICLR), Workshop Track},
  year={2016}
}

@inproceedings{pethick2025training,
      title={Training Deep Learning Models with Norm-Constrained {LMOs}}, 
      author={Thomas Pethick and Wanyun Xie and Kimon Antonakopoulos and Zhenyu Zhu and Antonio Silveti-Falls and Volkan Cevher},
      year={2025},
      booktitle={Proceedings of the International Conference on Machine Learning (ICML)},
}

@article{liu2025muon,
  author = {Jingyuan Liu and Jianlin Su and Xingcheng Yao and Zhejun Jiang and Guokun Lai and Yulun Du and Yidao Qin and Weixin Xu and Enzhe Lu and Junjie Yan and Yanru Chen and Huabin Zheng and Yibo Liu and Shaowei Liu and Bohong Yin and Weiran He and Han Zhu and Yuzhi Wang and Jianzhou Wang and Mengnan Dong and Zheng Zhang and Yongsheng Kang and Hao Zhang and Xinran Xu and Yutao Zhang and Yuxin Wu and Xinyu Zhou and Zhilin Yang},
  title = {Muon is Scalable For {LLM} Training},
  journal={arXiv preprint arXiv:2502.16982},
  year = {2025},
}

@article{fan1955some,
  title={Some metric inequalities in the space of matrices},
  author={Fan, Ky and Hoffman, Alan J.},
  journal={Proceedings of the American Mathematical Society},
  volume={6},
  number={1},
  pages={111--116},
  year={1955},
}

@article{liu2025cosmos,
      title={{COSMOS}: A Hybrid Adaptive Optimizer for Memory-Efficient Training of {LLMs}}, 
      author={Liming Liu and Zhenghao Xu and Zixuan Zhang and Hao Kang and Zichong Li and Chen Liang and Weizhu Chen and Tuo Zhao},
      year={2025},
      journal={arXiv preprint arXiv:2502.17410}, 
}

@book{bach2024learning,
  title={Learning Theory from First Principles},
  author={Bach, Francis},
  year={2024},
  publisher={MIT Press}
}

@misc{bernstein2025deriving,
  author = {Jeremy Bernstein},
  title = {Deriving {Muon}},
  url = {https://jeremybernste.in/writing/deriving-muon},
  month = {March},
  year = {2025}
}

@article{kovalev2025understanding,
      title={Understanding Gradient Orthogonalization for Deep Learning via Non-{E}uclidean Trust-Region Optimization}, 
      author={Dmitry Kovalev},
      year={2025},
      journal={arXiv preprint arXiv:2503.12645},
}

@article{medapati2025training,
	title={Training neural networks faster with minimal tuning using pre-computed lists of hyperparameters for {NAdamW}},
	author={Medapati, Sourabh and Kasimbeg, Priya and Krishnan, Shankar and Agarwal, Naman and Dahl, George},
	journal={arXiv preprint arXiv:2503.03986},
	year={2025}
}

@inproceedings{defazio2024the,
	title={The Road Less Scheduled},
	author={Aaron Defazio and Xingyu Alice Yang and Ahmed Khaled and Konstantin Mishchenko and Harsh Mehta and Ashok Cutkosky},
	booktitle={Advances in Neural Information Processing Systems (NeurIPS)},
	year={2024},
}

@inproceedings{xie2025structured,
  title={Structured Preconditioners in Adaptive Optimization: A Unified Analysis},
  author={Xie, Shuo and Wang, Tianhao and Reddi, Sashank and Kumar, Sanjiv and Li, Zhiyuan},
  booktitle={Proceedings of the International Conference on Machine Learning (ICML)},
  year={2025}
}

@article{an2025asgo,
      title={{ASGO}: Adaptive Structured Gradient Optimization}, 
      author={Kang An and Yuxing Liu and Rui Pan and Shiqian Ma and Donald Goldfarb and Tong Zhang},
      year={2025},
      journal={arXiv preprint arXiv:2503.20762},
}

@article{polyak1964some,
  author={Polyak, Boris T.},
  title = {Some methods of speeding up the convergence of iteration methods},
  journal = {USSR Computational Mathematics and Mathematical Physics},
  volume={4},
  number={5},
  pages={1--17},
  year={1964},
  publisher={Elsevier}
}

@article{amari1998natural,
  title={Natural gradient works efficiently in learning},
  author={Amari, Shun-Ichi},
  journal={Neural Computation},
  volume={10},
  number={2},
  pages={251--276},
  year={1998},
  publisher={MIT Press}
}

@inproceedings{karimi2016linear,
  title={Linear convergence of gradient and proximal-gradient methods under the {P}olyak-{{\L}}ojasiewicz condition},
  author={Karimi, Hamed and Nutini, Julie and Schmidt, Mark},
  booktitle={Joint European Conference on Machine Learning and Knowledge Discovery in Databases (ECML-PKDD)},
  year={2016},
}

@article{essentialai2025practical_full,
      title={Practical Efficiency of {Muon} for Pretraining}, 
      author={{Essential AI} and Ishaan Shah and Anthony M. Polloreno and Karl Stratos and Philip Monk and Adarsh Chaluvaraju and Andrew Hojel and Andrew Ma and Anil Thomas and Ashish Tanwer and Darsh J. Shah and Khoi Nguyen and Kurt Smith and Michael Callahan and Michael Pust and Mohit Parmar and Peter Rushton and Platon Mazarakis and Ritvik Kapila and Saurabh Srivastava and Somanshu Singla and Tim Romanski and Yash Vanjani and Ashish Vaswani},
      year={2025},
      journal={arXiv preprint arXiv:2505.02222},
}

@inproceedings{yao2021adahessian,
  title={{AdaHessian}: An adaptive second order optimizer for machine learning},
  author={Yao, Zhewei and Gholami, Amir and Shen, Sheng and Mustafa, Mustafa and Keutzer, Kurt and Mahoney, Michael},
  booktitle={Proceedings of the AAAI Conference on Artificial Intelligence},
  year={2021}
}

@article{autonne1902sur,
     author = {Autonne, Léon},
     title = {Sur les groupes lin\'eaires, r\'eels et orthogonaux},
     journal = {Bulletin de la Soci\'et\'e Math\'ematique de France},
     pages = {121--134},
     publisher = {Soci\'et\'e math\'ematique de France},
     volume = {30},
     year = {1902},
}

@article{qwen2025qwen_full,
      title={Qwen2.5 Technical Report}, 
      author={Qwen and An Yang and Baosong Yang and Beichen Zhang and Binyuan Hui and Bo Zheng and Bowen Yu and Chengyuan Li and Dayiheng Liu and Fei Huang and Haoran Wei and Huan Lin and Jian Yang and Jianhong Tu and Jianwei Zhang and Jianxin Yang and Jiaxi Yang and Jingren Zhou and Junyang Lin and Kai Dang and Keming Lu and Keqin Bao and Kexin Yang and Le Yu and Mei Li and Mingfeng Xue and Pei Zhang and Qin Zhu and Rui Men and Runji Lin and Tianhao Li and Tianyi Tang and Tingyu Xia and Xingzhang Ren and Xuancheng Ren and Yang Fan and Yang Su and Yichang Zhang and Yu Wan and Yuqiong Liu and Zeyu Cui and Zhenru Zhang and Zihan Qiu},
      year={2025},
      journal={arXiv preprint arXiv:2412.15115},
}

@article{amsel2025polar,
  title={The {Polar Express}: Optimal Matrix Sign Methods and Their Application to the {Muon} Algorithm},
  author={Amsel, Noah and Persson, David and Musco, Christopher and Gower, Robert},
  journal={arXiv preprint arXiv:2505.16932},
  year={2025}
}

@article{dong2025towards,
  title={Towards quantifying the {H}essian structure of neural networks},
  author={Dong, Zhaorui and Zhang, Yushun and Luo, Zhi-Quan and Yao, Jianfeng and Sun, Ruoyu},
  journal={arXiv preprint arXiv:2505.02809},
  year={2025}
}

@article{shen2025convergence,
  title={On the Convergence Analysis of {M}uon},
  author={Shen, Wei and Huang, Ruichuan and Huang, Minhui and Shen, Cong and Zhang, Jiawei},
  journal={arXiv preprint arXiv:2505.23737},
  year={2025}
}

@article{chen2025muon,
  title={Muon Optimizes Under Spectral Norm Constraints},
  author={Chen, Lizhang and Li, Jonathan and Liu, Qiang},
  journal={arXiv preprint arXiv:2506.15054},
  year={2025}
}

@article{chen2014stable,
  title={A stable scaling of {N}ewton-{S}chulz for improving the sign function computation of a {H}ermitian matrix},
  author={Chen, Jie and Chow, Edmond},
  journal={Preprint ANL/MCS-P5059-0114},
  year={2014}
}

@inproceedings{goldfarb2020practical,
  title={Practical quasi-{N}ewton methods for training deep neural networks},
  author={Goldfarb, Donald and Ren, Yi and Bahamou, Achraf},
  booktitle={Advances in Neural Information Processing Systems (NeurIPS)},
  year={2020}
}

@article{riabinin2025gluon,
  title={{Gluon}: Making {Muon} \& {Scion} Great Again! (Bridging Theory and Practice of {LMO}-based Optimizers for {LLMs})},
  author={Riabinin, Artem and Shulgin, Egor and Gruntkowska, Kaja and Richt{\'a}rik, Peter},
  journal={arXiv preprint arXiv:2505.13416},
  year={2025}
}

@article{ahn2025dion,
  title={Dion: Distributed Orthonormalized Updates},
  author={Kwangjun Ahn and Byron Xu and Natalie Abreu and Ying Fan and Gagik Magakyan and Pratyusha Sharma and Zheng Zhan and John Langford},
  journal={arXiv preprint arXiv:2504.05295},
  year={2025}
}

@article{grishina2025accelerating,
  title={Accelerating {N}ewton-{S}chulz Iteration for Orthogonalization via {C}hebyshev-type Polynomials},
  author={Grishina, Ekaterina and Smirnov, Matvey and Rakhuba, Maxim},
  journal={arXiv preprint arXiv:2506.10935},
  year={2025}
}

@article{shulgin2025beyond,
  title={Beyond the Ideal: Analyzing the Inexact {Muon} Update},
  author={Shulgin, Egor and AlRashed, Sultan and Orabona, Francesco and Richt{\'a}rik, Peter},
  journal={arXiv preprint arXiv:2510.19933},
  year={2025}
}

@article{armijo1966minimization,
  title={Minimization of functions having {L}ipschitz continuous first partial derivatives},
  author={Armijo, Larry},
  journal={Pacific Journal of Mathematics},
  volume={16},
  number={1},
  pages={1--3},
  year={1966},
  publisher={Mathematical Sciences Publishers}
}

@article{su2025isotropic,
  title={Isotropic Curvature Model for Understanding Deep Learning Optimization: Is Gradient Orthogonalization Optimal?},
  author={Su, Weijie},
  journal={arXiv preprint arXiv:2511.00674},
  year={2025}
}

@article{page2025muonall,
      title={{MuonAll}: {Muon} Variant for Efficient Finetuning of Large Language Models}, 
      author={Saurabh Page and Advait Joshi and S. S. Sonawane},
      year={2025},
      journal={arXiv preprint arXiv:2511.06086}
}

@misc{pytorch_muon2025,
	author = {{PyTorch Contributors}},
	title = {Muon},
	howpublished = {\url{https://docs.pytorch.org/docs/stable/generated/torch.optim.Muon.html}},
	note = {Accessed: 2025-12-17},
	year={2025}
}

@misc{optax_muon2025,
	author = {{Optax Contributors}},
	title = {Muon},
	howpublished = {\url{https://optax.readthedocs.io/en/stable/api/contrib.html#optax.contrib.muon}},
	note = {Accessed: 2025-12-17},
	year={2025}
}

@InProceedings{zeng2019global,
  title = "Global convergence of block coordinate descent in deep learning",
  author = 	 {Jinshan Zeng and Tim Tsz-Kit Lau and Shaobo Lin and Yuan Yao},
  booktitle = {Proceedings of the International Conference on Machine Learning (ICML)},
  year = 	 {2019},
}

@InProceedings{lau2018proximal,
	title = "A Proximal Block Coordinate Descent Algorithm for Deep Neural Network Training",
	author = 	 {Tim Tsz-Kit Lau and Jinshan Zeng and Baoyuan Wu and Yuan Yao},
	booktitle = {International Conference on Learning Representations (ICLR), Workshop Track},
	year = 	 {2018},
}
    	\end{CJK*}
    	}
        
        \newpage
        \appendix
            \addcontentsline{toc}{section}{\protect\textbf{Appendix}}
            \numberwithin{equation}{section}
            \numberwithin{theorem}{section}
            \numberwithin{algorithm}{section}
            \numberwithin{figure}{section}
            \numberwithin{table}{section}

            \begin{center}
                {\LARGE \textbf{Appendix}}
            \end{center}
        
            \tableofcontents
        
            \newpage
            In the appendix, we provide discussion on supplementary technical background materials, omitted proofs from the main text, as well as details of polar gradient methods. We also provide further details of the numerical experiments and additional numerical experiments.

            \section{Supplementary Technical Background}
            In this section, we provide supplementary technical background omitted from the main text due to space constraint. 
            
            \subsection{Convex Analysis}
            In the following, we introduce various notions from convex analysis which will be useful in the later parts of this paper, mostly taken from \citep{bauschke2017,beck2017}. For more background on convex analysis, we refer readers to standard texts such as \citep{bauschke2017,rockafellar1970,rockafellar1998,beck2017,mordukhovich2022convex}. 
            For generality, we consider a Euclidean space $\calE$ endowed with an inner product $\dotp{\cdot}{\cdot}$ and an associated norm $\norm{\cdot}$, which subsumes Euclidean spaces $\RR^d$ and $\RR^{m\times n}$. The following notions are indeed also well defined for more general infinite-dimensional spaces (i.e., real Hilbert spaces; see \citep{bauschke2017}), but we stick with finite-dimensional spaces for simplicity. 
            
            \begin{definition}[Subdifferential]
            	\label{def:subdiff}
            	Let $f\colon\calE\to\oRR$ be a proper function. The \emph{subdifferential} of $f$ is the set-valued operator 
            	\[\partial f \colon\calE\to2^{\calE}\colon x\mapsto \left\lbrace y\in\calE : (\forall z\in\calE)\;\; f(x) + \dotp{y}{z-x} \le f(z)  \right\rbrace.  \]
            	Let $x\in\calE$. Then $f$ is \emph{subdifferentiable} at $x$ if $\partial f(x) \ne\varnothing$; the elements of $\partial f(x)$ are the \emph{subgradients} of $f$ at $x$. In particular, if $f$ is convex and G\^{a}teaux differentiable at $x\in\calE$, the subdifferential of $f$ at $x$ is the set of gradients of $f$ at $x$, i.e., $\partial f(x) = \{\nabla f(x)\}$. 
            \end{definition}
    
            \begin{definition}[Fenchel conjugate]\label{def:conjugate}
            	The \emph{Fenchel conjugate} of a proper function $f\colon\calE\to\oRR$ is the function $f^*\colon \calE\to\RR\cup\{\pm\infty\}$ such that 
            	\[(\forall u\in\calE)\quad f^*(u) \coloneqq \sup_{x\in\dom f} \,\left\lbrace \dotp{x}{u} - f(x) \right\rbrace. \]		
            \end{definition}
            
            We now mention the famous Fenchel--Moreau theorem, which relates the biconjugate $f^{**}$ of $f$ and itself. 
            \begin{theorem}[Fenchel--Moreau]
            	\label{thm:Fenchel_Moreau}
            	Let $f\colon\euE\to\oRR$ be a proper function. Then $f$ is lower semi-continuous and convex if and only if $f^{**} = f$. In this case, $f^*$ is also proper. 
            \end{theorem}
                          
            \subsection{Matrix Analysis}
            We also include some notions and results from matrix analysis \citep{horn2012matrix,horn1994topics} which will be useful to understand some of the theoretical results and arguments of this paper. 
            
            Let us denote the vector space of $m\times n$ real matrices by $\calM_{m, n}$, and we only discuss the case over the field of real numbers $\RR$. 
            \begin{definition}[Matrix norm; \S5.6 of \citep{horn2012matrix}]
            	A function $\matsnorm{\cdot}{}\colon\calM_{m,n}\to\RR$ is a \emph{matrix norm} if, for all $A,B\in\calM_{m, n}$, it satisfies the following five axioms:
            	\begin{enumerate}[label=(\roman*)]
            		\item $\matsnorm{A}{}\ge0$ (nonnegativity)
            		\item $\matsnorm{A}{}=0$ if and only if $A=0$ (positivity)
            		\item $\matsnorm{cA}{} = |c|\matsnorm{A}{}$ for all $c\in\RR$ (homogeneity)
            		\item $\matsnorm{A+B}{} \le \matsnorm{A}{}\matsnorm{B}{}$ (triangle inequality)
            		\item $\matsnorm{AB}{} \le \matsnorm{A}{}\matsnorm{B}{}$ for $A\in\RR^{m\times p}$ and $B\in\RR^{p\times n}$ (submultiplicativity)
            	\end{enumerate}
            \end{definition}
            A norm on matrices which does not satisfy (v) submultiplicativity for all $A$ and $B$ is a \emph{vector norm on matrices}, sometimes called \emph{generalized matrix norm}. 
            In particular, the vector $\ell_\infty$-norm defined for $A\in\calM_{m, n}$, which is referred to as the \emph{max norm} in this paper, is not a matrix norm.     
            \begin{example}[Max norm]
            	The max norm $\matsnorm{\cdot}{\max}\colon X\in\RR^{m\times n}\mapsto \max_{1\le i\le m, 1\le j\le n} |x_{i,j}|$ is not a matrix norm. To see this, consider the matrix $J = \begin{pmatrix}
            		1 & 1  \\ 1 & 1 
            	\end{pmatrix}\in\RR^{2\times 2}$. Then $J^2 = 2J$, $\matsnorm{J}{\max} = 1$, $\matsnorm{J^2}{\max} = \matsnorm{2J}{\max} = 2\matsnorm{J}{\max} = 2$. Thus, the max norm is not submultiplicative as $\matsnorm{J^2}{\max} > \matsnorm{J}{\max}^2$. However, a scalar multiple of the max norm on matrices is a matrix norm. Indeed, $\sqrt{mn}\matsnorm{\cdot}{\max}$ on $\calM_{m, n}$ is a matrix norm.            
            \end{example}
            
            Next, we introduce the notion of invariant matrix norms, which is originated from \citep{vonneumann1937some}. Unitarily invariant norms have important implications for the polar decomposition from a matrix approximation perspective. There is a long line of works studying this class of objects in matrix analysis and linear algebra, e.g.,    \citep{mirsky1960symmetric,watson1993matrix,watson1991algorithm}. They are also a central notion to convex matrix analysis and eigenvalue optimization \citep{lewis1995convex,lewis1996group,lewis1996eigenvalue,lewis2003mathematics}.

            In the following, we mainly follow the notation from Chapters 5.6 and 7.4.7 of \citep{horn2012matrix}. Let us denote the set of $d\times d$ real orthogonal matrices by $\OO^d$. 
            \begin{definition}[Unitarily invariant norm]
            	A norm $\matsnorm{\cdot}{}$ on $\calM_{m, n}$ (not necessarily a matrix norm) is said to be \emph{unitarily (or orthogonally) invariant} if $\matsnorm{UXV}{} = \matsnorm{X}{}$ for any $X\in\RR^{m \times n}$, $U\in\OO^m$, $V\in\OO^n$. Unitarily invariant matrix norm is a unitarily invariant norm on $\calM_{m, n}$ that is submultiplicative.         
            \end{definition}
            
            A famous fundamental result of \citet{vonneumann1937some} states that unitarily invariant matrix norms can be characterized as composite functions of the form $\varphi(\sigma(\cdot)) = \varphi \circ\sigma$, where $\varphi$ is a \emph{symmetric gauge function} and $\sigma$ is the singular value function. In what follows, we define $m\wedge n \coloneqq\min\{m,n\}$. 
            \begin{definition}[Symmetric gauge function]
            	A function $\varphi\colon\RR^{m\wedge n}\to\Rp$ is said to be a \emph{symmetric gauge function} if $\varphi$ is an absolute, permutation-invariant norm on the components. 
            \end{definition}
            
            \begin{proposition}
            	Any unitarily invariant norm $\matsnorm{\cdot}{}$ can be written as $\matsnorm{X}{} = \varphi(\sigma(X)) = (\varphi\circ\sigma)(X)$, where $\sigma\colon\RR^{m\times n}\to\RR^{m\wedge n}$ has components $\sigma_1(X)\ge\cdots\ge\sigma_{m\wedge n}(X)\ge 0$ which are the singular values of $X$, and $\varphi\colon\RR^{m\wedge n}\to\RR$ is a \emph{symmetric gauge function}. 
            \end{proposition}
            
            Thus, for any unitarily invariant norm $\matsnorm{\cdot}{}$, we have $\matsnorm{X}{} = \matsnorm{\Diag(\sigma(X))}{}$.     
            The unitarily invariant norms on $\calM_{m, n}$ determined by the $\ell_p$-norm as its symmetric gauge function are known as Schatten $p$-norms. 
            \begin{example}[Schatten $p$-norm]
            	If the symmetric gauge function $\varphi = \norm{\cdot}_p$ is the $\ell_p$-norm, where $1<p\le\infty$, then the Schatten $p$-norm $\matsnorm{\cdot}{p}$ is a unitarily invariant norm.         
            	The nuclear norm $\nucnorm{\cdot}$ is the Schatten $1$-norm. The Frobenius norm $\fronorm{\cdot}$ is the Schatten $2$-norm. The spectral norm $\specnorm{\cdot}$ is the Schatten $\infty$-norm. 
            \end{example}
            
            However, the max norm and the $\ell_p\to\ell_q$ operator norm are not unitarily invariant in general. 
            \begin{example}[Max norm]
            	The max norm $\matsnorm{\cdot}{\max}$ is \emph{not} a unitarily invariant norm. 
            \end{example}
            
            \begin{example}[$\ell_p\to\ell_q$ operator norm]
            	Let $\norm{\cdot}_p$ and $\norm{\cdot}_q$ be the $\ell_p$-norm and $\ell_q$-norm on $\RR^n$ and $\RR^m$, respectively. Then, the $\ell_p\to\ell_q$ operator norm on $\calM_{m, n}$ is defined by the variational characterization
            	\begin{equation*}
            		\matsnorm{X}{p\to q} \coloneqq \max_{u\in\RR^n\setminus\{0_n\}} \frac{\norm{Xu}_q}{\norm{u}_p}.
            	\end{equation*}
            	The $\ell_p\to\ell_q$ operator norm is \emph{not} unitarily invariant in general, except when $p=q=2$ it becomes the spectral norm. 
            \end{example}
            
            To understand the importance of unitarily invariant norms for characterizing best approximation properties of the orthogonal polar factor in polar decomposition, we state the following theorem by \citet{fan1955some}.
            \begin{theorem}\label{thm:polar}
            	Let $A\in\RR^{m\times n}$ ($m\ge n$) have the polar decomposition $A=U_\sfp H$. Then 
            	\[U_\sfp \in \argmin_{Q\in\RR^{m\times n} : Q^\top Q=I_n} \matsnorm{A - Q}{}\] 
            	for any unitarily invariant norm $\matsnorm{\cdot}{}$. The minimizer $U_\sfp$ is unique for the Frobenius norm if $A$ has full rank. 
            \end{theorem}
            Hence, the orthogonal polar factor $U_\sfp$ is the nearest matrix to $A$ with orthonormal columns. This justifies that the polar decomposition offers an optimal way of orthogonalizing a matrix. 
            
            Next, we state a result regarding the subdifferential of the dual norm $\matsnorm{\cdot}{}^*$ of a matrix $X\in\RR^{m\times n}$ and the set of dual matrices of $X$ in the original (primal) norm $\matsnorm{\cdot}{}$. 
            Let $\matsnorm{\cdot}{}$ be a norm on $\calM_{m, n}$. Let us recall from \Cref{def:subdiff} that the \emph{subdifferential} of $\matsnorm{X}{}$ is defined by 
            \begin{equation*}
            	\partial \matsnorm{X}{} = \left\{ Y\in\RR^{m\times n} : (\forall Z\in\RR^{m\times n})\;\; \matsnorm{X}{} + \dotpF{Y}{Z-X} \le \matsnorm{Z}{} \right\}. 
            \end{equation*}

            We now state the following proposition from \citep{watson1992characterization,ziketak1993subdifferentials}, which offers a way of computing the dual norm (i.e., the Fenchel conjugate of the norm) of a matrix if the subdifferential of its dual norm is available. 
            
            \begin{proposition}\label{prop:subdiff_dual}
            	Is is known that the subgradient $G\in\partial \matsnorm{X}{}$ is equivalent to $\matsnorm{X}{} = \dotpF{G}{X}$ and $\matsnorm{G}{}^*\le1$, where $\matsnorm{\cdot}{}^*$ is the \emph{dual norm} of $\matsnorm{\cdot}{}$ defined by 
            	\begin{equation}
            		\matsnorm{G}{}^* = \sup_{Z:\matsnorm{Z}{} \le 1} \dotpF{Z}{G}. 
            	\end{equation}
            	It follows that the subdifferential of $\matsnorm{X}{}^*$ is the set of $\matsnorm{\cdot}{}$-dual matrices of $X$, i.e., 
            	\[\partial\matsnorm{X}{}^* = \{G\in\RR^{m\times n}: \dotpF{X}{G} = \matsnorm{X}{}^*, \matsnorm{G}{} = 1\} \eqqcolon \VV_{\matsnorm{\cdot}{}}(X). \] 
            \end{proposition}
            Consequently, we can compute the set of $\matsnorm{\cdot}{}$-dual matrices of $X$ through the subdifferential $\partial\matsnorm{X}{}^*$.     
            Furthermore, since norms are continuous and convex, by the Fenchel--Moreau theorem (\Cref{thm:Fenchel_Moreau}), the set of $\matsnorm{\cdot}{}^*$-dual matrices of $X$ can also be computed through the subdifferential $\partial\matsnorm{X}{}$. This result is particularly useful for Schatten $p$-norms  (generally unitarily invariant matrix norms) and $\ell_p\to\ell_q$ operator norms since their subdifferentials are generally known in closed forms \citep{ziketak1988characterization,watson1992characterization,ziketak1993subdifferentials}.

            \subsection{Numerical Polar Decomposition Algorithms}
            \label{subsec:NLA}    
            The polar decomposition is an important matrix decomposition in matrix analysis and numerical linear algebra \citep{higham2008functions}. Efficiently computing the polar decomposition of matrices is rudimentary for the practical use of polar gradient methods. In this subsection, we go through various existing numerical polar decomposition algorithms from the numerical linear algebra literature.

            \subsubsection{Details of Numerical Polar Decomposition Algorithms}
            There are numerous numerical algorithms for computing the polar decomposition of a matrix $A\in\RR^{m\times n}$ ($m\ge n$) in the numerical linear algebra literature. We include the pseudocode of these numerical polar decomposition algorithms for readers' convenience. 
            
            The first one is the scaled Newton iteration, which can be found in Chapter 8.6 of \citep{higham2008functions} with different scaling schemes $\mu_k$. 
            
            \begin{algorithm}[H]
            	\caption{Scaled Newton iteration}
            	\label{alg:newton}
            	\begin{algorithmic}[1]
            		\REQUIRE $A\in\RR^{m\times n}$, scaling $(\mu_k)_{1\le k\le K}$
            		\STATE $X_0 = A$
            		\FOR{$k=0, \ldots, K-1$}
            		\STATE $X_{k+1} = \frac12 \left(\mu_k X_k + \mu_k^{-1}X_k^{-\top}\right)$
            		\ENDFOR
            		\ENSURE $U_\sfp = X_K$, $H = \frac12\left(U_\sfp^\top A + (U_\sfp^\top A)^\top \right)$
            	\end{algorithmic}
            \end{algorithm}

            The Newton--Schulz (NS) iteration (\Cref{alg:newton_schulz}), on the other hand, does not involve any matrix inverse. While the original Newton--Schulz iteration in \citep{higham2008functions} makes use of a cubic polynomial, a degree-5 polynomial is used in the implementation of \Muon. The matrix iterative polynomial coefficients in the Newton--Schulz iteration are tuned through gradient descent on heuristic objectives to accelerate convergence, given in \Cref{alg:newton_schulz_5}.   
            \begin{algorithm}[H]
            	\caption{Newton--Schulz iteration (classical)}
            	\label{alg:newton_schulz}
            	\begin{algorithmic}[1]
            		\REQUIRE $A\in\RR^{m\times n}$, small $\delta>0$ 
            		\STATE $X_0 = (\sqrt{3}-\delta)A/\fronorm{A}$
            		\FOR{$k=0, \ldots, K-1$}
            		\STATE $X_{k+1} = \frac12 X_k\left(3I_n + X_k^\top X_k\right)$
            		\ENDFOR
            		\ENSURE $U_\sfp = X_K$, $H = \frac12\left(U_\sfp^\top A + (U_\sfp^\top A)^\top \right)$
            	\end{algorithmic}
            \end{algorithm}

            \begin{algorithm}[H]
                \caption{Newton--Schulz iteration in \Muon \citep{jordan2024muon}}
                \label{alg:newton_schulz_5}
                \begin{algorithmic}[1]
                    \REQUIRE $A\in\RR^{m\times n}$, iterative polynomial coefficients $(a,b,c)=(3.4445, -4.775, 2.0315)$, $\delta=10^{-7}$ 
                    \STATE $X_0 = A/(\fronorm{A} + \delta)$
                    \FOR{$k=0, \ldots, K-1$}
                    \STATE $M_k = X_k^\top X_k$
                    \STATE $X_{k+1} = a X_k + X_k\left(bM_k + cM_k^2\right)$
                    \ENDFOR
                    \ENSURE $U_\sfp = X_K$, $H = \frac12\left(U_\sfp^\top A + (U_\sfp^\top A)^\top \right)$
                \end{algorithmic}
            \end{algorithm}
            
            Unfortunately, the coefficient scheme $(a,b,c)=(3.4445, -4.775, 2.0315)$ in \Cref{alg:newton_schulz_5} does not converge to the desired orthogonal polar factor, as pointed out in \citep{amsel2025polar}. The authors of \citep{amsel2025polar} propose the \textsc{Polar Express}, which dynamically determines the polynomial coefficients at each iteration and converges to the orthogonal polar factor super-exponentially. We refer readers to the paper \citep{amsel2025polar} for the full algorithmic details of the \textsc{Polar Express}. The concurrent work CANS \citep{grishina2025accelerating} is also developed in a similar spirit. 
            
            The QR-based Dynamically Weighted Halley (QDWH) algorithm \citep{nakatsukasa2013stable} is a more recent algorithm based on the QR decomposition and is globally and asymptotically cubically convergent. Its main principle is to derive a dynamic weighting scheme for Halley's iteration, unlike the hand-picked coefficient scheme in the NS iteration in \Muon. It also does not involve explicit matrix inversions, and hence is less likely to suffer from numerical stability issues and minimizes the communication costs by using communication friendly matrix operations such as the QR decomposition (without pivoting).

            \begin{algorithm}[H]
            	\caption{The QR-based Dynamically Weighted Halley (QDWH) algorithm}
            	\label{alg:QWDH}
            	\begin{algorithmic}[1]
            		\REQUIRE $A\in\RR^{m\times n}$
            		\STATE Estimate $\alpha\gtrsim\sigma_{\max}(A)$, $\beta\lesssim\sigma_{\min}(A)$, $X_0=A/\alpha$, $\ell_0=\beta/\alpha$ 
            		\FOR{$k=0, \ldots, K-1$}
            		\STATE $a_k = h(\ell_k)$, $b_k = (a_k - 1)^2/4$, $c_k = a_k + b_k - 1$, where $h(\ell) = \sqrt{1+\gamma} + \frac12\left(8- 4\gamma + 8(2-\ell^2)/(\ell^2\sqrt{1+\gamma})\right)^{\negthickspace\half}$, $\gamma=\gamma(\ell) = \left(4(1-\ell^2)/\ell^4\right)^{\negthickspace\sfrac13}$
            		\STATE Compute QR decomposition $\begin{pmatrix}
            			\sqrt{c_k}X_k \\ I
            		\end{pmatrix} = \begin{pmatrix}
            			Q_1 \\ Q_2
            		\end{pmatrix} R$
            		\STATE $X_{k+1} = (b_k/c_k) X_k + (1/\sqrt{c_k})(a_k - b_k/c_k)Q_1Q_2^\top$
            		\STATE $\ell_{k+1} = \ell_k(a_k + b_k \ell_k^2)/(1+c_k\ell_k^2)$
            		\ENDFOR
            		\ENSURE $U_\sfp = X_K$, $H = \frac12\left(U_\sfp^\top A + (U_\sfp^\top A)^\top \right)$
            	\end{algorithmic}
            \end{algorithm}
            For the practical implementation of the QDWH algorithm, we only need estimates $\hat\alpha$ and $\hat\beta$ of $\alpha$ and $\beta$ satisfying $0 < \hat\beta \le \sigma_{\min}(A) \le \sigma_{\max}(A) \le \hat\alpha$; see Section 4 of \citep{nakatsukasa2010optimizing} for more details. Since the QR decomposition is involved in the QDWH algorithm, its performance relies heavily on the efficiency of the QR decomposition in the deep learning library such as PyTorch \citep{paszke2019pytorch} and JAX \citep{jax2018github}. Notice that the computation of the polar decomposition is available on JAX (\texttt{jax.scipy.linalg.polar}), where the QDWH algorithm is one of the available methods (\texttt{jax.lax.linalg.qdwh}), with the other being the SVD.

            The ZOLO-based Polar Decomposition (ZOLO-PD) algorithm \citep{nakatsukasa2016computing} is a higher-order variant of the QDWH algorithm for the polar decomposition, based on the best rational approximation for the scalar sign function due to Zolotarev in 1877. It converges in just two iterations in double-precision arithmetic with the rate of convergence seventeen. The double-precision requirement might however not be suitable for large-scale applications. It can however be parallelized since the $r$ QR decompositions involved can be performed independently. Therefore, with parallelized implementation, the ZOLO-PD algorithm can be faster than the QDWH algorithm. 
            
            The QDWH and ZOLO-PD algorithms can also be coupled with spectral divide and conquer algorithms for the symmetric eigenvalue problem and computing the singular value decomposition \citep{nakatsukasa2013stable,nakatsukasa2016computing}. 
            
            \begin{algorithm}[h]
            	\caption{The ZOLO-based Polar Decomposition (ZOLO-PD) algorithm}
            	\label{alg:ZOLO-PD}
            	\begin{algorithmic}[1]
            		\REQUIRE $A\in\RR^{m\times n}$, the unit roundoff of IEEE double-precision arithmetic $u=2^{-53}\approx 1.1\times 10^{-16}$
            		\STATE Estimate $\alpha\gtrsim\sigma_{\max}(A)$, $\beta\lesssim\sigma_{\min}(A)$, $X_0=A/\alpha$, $\ell=\beta/\alpha$. 
            		\STATE Choose $r$ based on $\kappa=\ell^{-1}$ from \Cref{tab:required_iterations}. If $\kappa<2$, then $X_1=A$ and skip to (iv). 
            		\STATE Compute $X_1$ and $X_2$:
            		\begin{enumerate}[label=(\roman*)]
            			\item Compute $c_j = \ell^2\mathrm{sn}^2\left(\frac{iK'}{2r+1}; \ell'\right) / \mathrm{cn}^2\left(\frac{iK'}{2r+1}; \ell'\right)$, where $\mathrm{sn}(u; \ell')$ and $\mathrm{cn}(u; \ell')$ are the Jacobi elliptic functions. Also compute $a_j = -\left( \prod_{k=1}^r(c_{2j-1}-c_{2k})\right) \cdot\left(\prod_{k=1,k\ne j}^r(c_{2j-1}-c_{2k-1}) \right)$. 
            			
            			\item Compute $X_1$ by $\hat{M}=\prod_{j=1}^r (1+c_{2j-1})/(1+{2j})$ and $r$ QR decompositions
            			\begin{empheq}[left=\empheqlbrace]{equation}
            				\begin{aligned}                            
            					\begin{pmatrix}
            						X_0 \\ \sqrt{c_{2j-1}I}
            					\end{pmatrix}
            					&= 
            					\begin{pmatrix}
            						Q_{j1} \\ Q_{j2}
            					\end{pmatrix} R_j, \quad j=1,2, \ldots, r, \\
            					X_1 &= \hat{M}\left(X_0 + \sum_{j=1}^r \frac{a_j}{\sqrt{c_{2j-1}}} Q_{j1}Q_{j2}^\top\right). 
            				\end{aligned}
            			\end{empheq}
            			
            			\item Update $\ell\coloneqq\hat{M}\ell\prod_{j=1}^r (\ell^2 + c_{2j})/(\ell^2+c_{2j-1})$ and recompute $c_j$ and $a_j$ as in step (i). 
            			
            			\item Compute $X_2$ by $\hat{M}=\prod_{j=1}^r (1+c_{2j-1})/(1+{2j})$ and 
            			\begin{empheq}[left=\empheqlbrace]{equation}
            				\begin{aligned}
            					Z_{2j-1} &= X_1^\top X_1 + c_{2j-1}I, \quad W_{2j-1} = \mathrm{Chol}(Z_{2j-1}), \\
            					X_2 &= \hat{M}\left(X_1 + \sum_{j=1}^r a_j(X_2W_{2j-1}^{-1})W_{2j-1}^{-\top} \right), 
            				\end{aligned}
            			\end{empheq}
            			where $\mathrm{Chol}$ denotes the Cholesky factor in the Cholesky decomposition of a symmetric positive definite matrix. 
            			
            			Verify that $\fronorm{X_2-X_1}/\fronorm{X_2}\le u^{1/(2r+1)}$ holds. If not, return to step 1 with $A=X_2$. 
            		\end{enumerate}
            		\ENSURE $U_\sfp=X_2$, $H = \frac12\left(U_\sfp^\top A + (U_\sfp^\top A)^\top \right)$
            	\end{algorithmic}
            \end{algorithm}
            
            \begin{table}[h]
            	\centering        
            	\caption{Required number of iterations for varying $\kappa_2(A)$ and $r$, obtained as the smallest $k$ for which $\hat{Z}_{(2r+1)^k}([\ell, 1]) \subseteq [1-\scrO(u), 1]$.}
            	\label{tab:required_iterations}
            	\vspace*{2mm}
            	\begin{tabular}{ccccccccccccc}
            		\toprule
            		$\kappa_2(A)$ & $1.001$ & $1.01$ & $1.1$ & $1.2$ & $1.5$ & $2$ & $10$ & $10^2$ & $10^3$ & $10^5$ & $10^7$ & $10^{16}$ \\
            		\midrule
            		$r=1$ (QDWH) & 2 & 2 & 2 & 3 & 3 & 3 & 4 & 4 & 4 & 5 & 5 & 6 \\
            		$r=2$ & 1 & 2 & 2 & 2 & 2 & 2 & 3 & 3 & 3 & 3 & 4 & 4 \\
            		$r=3$ & 1 & 1 & 2 & 2 & 2 & 2 & 2 & 2 & 3 & 3 & 3 & 3 \\
            		$r=4$ & 1 & 1 & 1 & 2 & 2 & 2 & 2 & 2 & 2 & 3 & 3 & 3 \\
            		$r=5$ & 1 & 1 & 1 & 1 & 2 & 2 & 2 & 2 & 2 & 2 & 3 & 3 \\
            		$r=6$ & 1 & 1 & 1 & 1 & 1 & 2 & 2 & 2 & 2 & 2 & 2 & 3 \\
            		$r=7$ & 1 & 1 & 1 & 1 & 1 & 1 & 2 & 2 & 2 & 2 & 2 & 3 \\
            		$r=8$ & 1 & 1 & 1 & 1 & 1 & 1 & 2 & 2 & 2 & 2 & 2 & 2 \\
            		\bottomrule
            	\end{tabular}        
            \end{table}

            \begin{remark}
            	In the numerical experiments of this paper, the implementation of the QDWH and ZOLO-PD algorithms are based on translation of the MATLAB code of the original paper \citep{nakatsukasa2010optimizing,nakatsukasa2016computing} into PyTorch, which is not optimized for large-scale neural network training. We leave their more optimized implementations for future work. We also believe that the QDWH algorithm implementation in JAX is more efficient but we stick with PyTorch for our experiments. See also \citep{lewis2022large} for large-scale distributed numerical linear algebra algorithms with tensor processing units (TPUs) for the potential of such directions. 
            \end{remark}

            \subsubsection{Backward Stability of Polar Decomposition Algorithms}
            \label{subsubsec:stability}
            In addition to computational efficiency, numerical stability of polar decomposition algorithms is of vital importance to our choice of applications. The notion of \emph{backward stability} of a polar decomposition algorithm \citep{nakatsukasa2012backward} is such one that determines the numerical stability of polar decomposition algorithms. 
            
            \begin{definition}[Backward stability]
            	The polar decomposition of $A$ is said to be computed in a \emph{backward stable} manner if the computed polar factors $\hat{U}_\sfp$ and $\hat{H}$ satisfy that $\hat{H}$ is symmetric, 
            	\[\fronorm{A - \hat{U}_\sfp\hat{H}}/\fronorm{A} = \scrO(u) \quad \text{and}\quad \fronorm{\hat{U}_\sfp^\top\hat{U}_\sfp - I_n}/\sqrt{n} = \scrO(u), \]
            	where $u=2^{-53}\approx1.1\times10^{-16}$ is the unit roundoff for IEEE double precision arithmetic. 
            \end{definition}

            The Newton--Schulz iteration (original form) is only \emph{conditionally stable} \citep{nakatsukasa2012backward}, meaning that it is stable away from, but not very close to, the boundary of its region of convergence. The initialization $X_0$ needs to have norm safely less than $\sqrt{3}$ and to be not too ill-conditioned, i.e., with a small condition number $\kappa_2(X_0)$.
            The QDWH algorithm is backward stable under the assumption that the QR decompositions involved are performed with row sorting (or pivoting) and column sorting      \citep{nakatsukasa2012backward}.    
            Backward stability of the ZOLO-PD algorithm is only demonstrated experimentally \citep{nakatsukasa2016computing}; its proof remains an open problem.

            \section{Details of Polar Gradient Methods}
            \label{sec:details_polargrad}
            We now provide further details of the class of polar gradient methods, as a broad class of matrix optimization algorithms in this section. 
            
            \subsection{\PolarGrad, \PolarGradM, \PolarMuon and \PolarHB}
            We first give the full pseudocode of several optimizers in the \PolarGrad family. Note that the polar decomposition in the following pseudocode are replaced by an inexact polar oracle (i.e., a numerical polar decomposition algorithm) in practice. 
            
            \begin{algorithm}[h!]
            	\caption{\PolarGrad/\PolarSGD (with Decoupled Weight Decay) (\textsc{PolarGrad/PolarSGD(W)})}
            	\label{alg:polar_grad}
            	\begin{algorithmic}
            		\REQUIRE $\{\gamma_k\}_{k=1}^K\subset\Rpp$, $X_0\in\RR^{m\times n}$, $M_0 = 0_{m\times n}$
            		\FOR{$k=0, \ldots, K-1$}
            		\STATE $G_k=\nabla f(X_k)$ or $G_k = \nabla f(X_k,\xi_k)$ with $\xi_k\sim\calD$
            		\STATE $U_kH_k = \polar(G_k)$
            		\STATE $\nu_k = \nucnorm{G_k} \equiv \dotpF{G_k}{U_k} = \tr(H_k)$
            		\STATE $X_{k+1} = (1-\lambda\gamma_k)X_k - \gamma_k \nu_k U_k$
            		\ENDFOR
            	\end{algorithmic}
            \end{algorithm}

            \begin{algorithm}[h!]
            	\caption{\PolarGrad/\PolarSGD with Momentum-First EMA Momentum (and Decoupled Weight Decay) (\textsc{PolarGradM/PolarSGDM(W)}) or \PolarMuon}
            	\label{alg:polar_grad_ema_momentum_first}
            	\begin{algorithmic}
            		\REQUIRE $\{\gamma_k\}_{k=1}^K\subset\Rpp$, $\beta\in(0,1)$, $\lambda\ge0$, $X_0\in\RR^{m\times n}$, $M_0 = 0_{m\times n}$
            		\FOR{$k=0, \ldots, K-1$}
            		\STATE $G_k=\nabla f(X_k)$ or $G_k = \nabla f(X_k,\xi_k)$ with $\xi_k\sim\calD$
            		\STATE $M_k = \beta M_{k-1} + (1-\beta)G_k$    		
            		\STATE $U_kH_k = \polar(M_k)$
            		\STATE $\nu_k = \nucnorm{M_k} \equiv \dotpF{M_k}{U_k} = \tr(H_k)$
            		\STATE $X_{k+1} = (1-\lambda\gamma_k)X_k - \gamma_k \nu_k U_k$
            		\ENDFOR
            	\end{algorithmic}
            \end{algorithm}

            \begin{algorithm}[h!]
            	\caption{\PolarGrad/\PolarSGD with Polar-First EMA Momentum (and Decoupled Weight Decay) (\textsc{PolarGradM/PolarSGDM(W)})}
            	\label{alg:polar_grad_ema}
            	\begin{algorithmic}
            		\REQUIRE $\{\gamma_k\}_{k=1}^K\subset\Rpp$, $\beta\in(0,1)$, $\lambda\ge0$, $X_0\in\RR^{m\times n}$, $M_0 = 0_{m\times n}$
            		\FOR{$k=0, \ldots, K-1$}
            		\STATE $G_k=\nabla f(X_k)$ or $G_k = \nabla f(X_k,\xi_k)$ with $\xi_k\sim\calD$
            		\STATE $U_kH_k = \polar(G_k)$
            		\STATE $\nu_k = \nucnorm{G_k} \equiv \dotpF{G_k}{U_k} = \tr(H_k)$
            		\STATE $M_k = \beta M_{k-1} + (1-\beta)U_k$    		
            		\STATE $X_{k+1} = (1-\lambda\gamma_k)X_k - \gamma_k \nu_k M_k$
            		\ENDFOR
            	\end{algorithmic}
            \end{algorithm}

            \begin{algorithm}[h!]
            	\caption{\PolarGrad or \PolarSGD with (Momentum-First) Polyak's Heavy Ball Momentum (and Decoupled Weight Decay) (\textsc{PolarHB(W)})}
            	\label{alg:polar_grad_hb}
            	\begin{algorithmic}
            		\REQUIRE $\{\gamma_k\}_{k=1}^K\subset\Rpp$, $\beta\in(0,1)$, $\lambda\ge0$, $X_0\in\RR^{m\times n}$, $M_0 = 0_{m\times n}$
            		\FOR{$k=0, \ldots, K-1$}
            		\STATE $G_k=\nabla f(X_k)$ or $G_k = \nabla f(X_k,\xi_k)$ with $\xi_k\sim\calD$
            		\STATE $M_k = \beta M_{k-1} + G_k$    		
            		\STATE $U_kH_k = \polar(M_k)$
            		\STATE $\nu_k = \nucnorm{M_k} \equiv \dotpF{M_k}{U_k} = \tr(H_k)$
            		\STATE $X_{k+1} = (1-\lambda\gamma_k)X_k - \gamma_k \nu_k U_k$
            		\ENDFOR
            	\end{algorithmic}
            \end{algorithm}

            \subsection{Steepest Descent with respect to The $\ell_\infty$-Norm and The Spectral Norm as Preconditioned Gradient Methods with Explicit and Implicit Preconditioners}
            \label{subsec:explicit_implicit_2}
            Following our discussion in \Cref{subsec:explicit_implicit}, we further explain that the (unnormalized) steepest descent w.r.t.~the (squared) $\ell_\infty$-norm and the (squared) spectral norm can both be interpreted as preconditioned gradient methods with either explicit or implicit preconditioners. 
            
            \subsubsection{Unnormalized Sign Descent}
            Let us recall from \eqref{eqn:signSGD} that 
            \[(\forall k\in\NN)\quad x_{k+1} = \argmin_{x\in\RR^d} \,\left\{ \dotp{g_k}{x - x_k} + \frac{1}{2\gamma_k}\infnorm{x - x_k}^2 \right\} = x_k - \gamma_k \onenorm{g_k}\cdot\sgn(g_k). \]
            If we define the explicit preconditioner $P_k\coloneqq\Diag(g_k^2)^{\half} = \Diag(|g_k|)$, then we have 
            \begin{equation*}
            	x_{k+1} = x_k - \gamma_k \onenorm{g_k}\cdot\sgn(g_k) = x_k - \gamma_k \,\tr(P_k)\, P_k^{-1}\, g_k. 
            \end{equation*}
            Consequently, the sign descent method can be viewed as either an explicit preconditioned gradient method with an explicit preconditioner $P_k$ scaled by its trace or an implicit preconditioned gradient method with an implicit preconditioner or preconditioning function $g\mapsto\onenorm{g}\cdot\sgn(g)$.

            \subsubsection{\PolarGrad}
            For \PolarGrad, due to the different definitions of the polar decomposition of a matrix $X\in\RR^{m\times n}$ based on its numbers of rows and columns, we separate our discussion for the cases of $m\ge n$ and $m<n$. 
            
            If $m\ge n$, we recall from \eqref{eqn:polargrad} that 
            \[
            (\forall k\in\NN)\quad U_{\sfp, k}H_k = \polar(G_k) \quad \text{and}\quad X_{k+1} = X_k - \gamma_k \,\tr(H_k)\, U_{\sfp, k}, 
            \]
            where the symmetric polar factor $H_k = (G_k^\top G_k)^{\half} = V_k\Sigma_kV_k^\top$ with $U_k\Sigma_k V_k^\top = \SVD(G_k)$. It turns out that the explicit (right) preconditioner $P_k$ in this case is just the symmetric polar factor $H_k$ itself, since $P_k^{-1} = V_k\Sigma_k^{-1}V_k^\top$ implies $G_kP_k^{-1} = U_kV_k^\top = U_{\sfp, k}$. That is to say, the update of \PolarGrad can be written as 
            \[X_{k+1} = X_k - \gamma_k\, \tr(H_k)\, U_{\sfp, k} = X_k - \gamma_k\,\tr(P_k)\, G_k P_k^{-1}. \]
            
            If $m<n$, the update of \PolarGrad becomes
            \begin{equation*}
            	(\forall k\in\NN)\quad H_kU_{\sfp, k} = \polar(G_k) \quad \text{and}\quad X_{k+1} = X_k - \gamma_k \,\tr(H_k)\, U_{\sfp, k}, 
            \end{equation*}
            where the symmetric polar factor $H_k = (G_k G_k^\top)^{\half} = U_k\Sigma_kU_k^\top$ with $U_k\Sigma_k V_k^\top = \SVD(G_k)$. In this case, the explicit (left) preconditioner $P_k$ is again the symmetric polar factor $H_k$ itself, since $P_k^{-1} = U_k\Sigma_k^{-1}U_k^\top$ implies $P_k^{-1}G_k = U_kV_k^\top = U_{\sfp, k}$. The update of \PolarGrad can then be written as 
            \begin{equation*}
            	X_{k+1} = X_k - \gamma_k \,\tr(H_k)\, U_{\sfp, k} = X_k - \gamma_k\,\tr(P_k) P_k^{-1} G_k. 
            \end{equation*}
            As a result, \PolarGrad can be viewed as either an explicit preconditioned gradient method with an explicit (left or right) preconditioner $P_k$ scaled by its trace or an implicit preconditioned gradient method with an implicit (left or right) preconditioner or preconditioning function $U_\sfp H = G\mapsto\tr(H)\,U_\sfp=\nucnorm{G}\cdot\msgn(G)$.

            \section{Details and Additional Results of Numerical Experiments}
            \label{sec:details}
            The simulated data experiments are performed on a Mac mini with an M4 CPU and 16 GB memory. The language model pre-training experiments for Qwen2.5, GPT-2 Small and Medium are performed on eight NVIDIA H100-SXM5 GPUs. Each set of the simulated data experiments is repeated with three different random seeds. The results of the last two seeds are reported in this section. 
            For the Newton--Schulz iteration in \Muon, we use the default coefficients $(a,b,c) = (3.4445, -4.7750,  2.0315)$ of the matrix iterative polynomial in its original implementation. 
            
        \subsection{Matrix Quadratic Regression}
            The initialization $X_0$ has entries independently drawn from $\Unif(-1, 1)$. The matrices $A$, $B$ and $C$ have independent standard Gaussian entries. 
            No weight decay is used in all optimizers. The learning rate decay schedule is a step scheduler which multiplies the base learning rate $\gamma_0$ by $0.99$ every $25$ steps. The optimizer hyperparameters are given in the table below. Default hyperparameters of \Adam in PyTorch are used ($\varepsilon=10^{-8}$).     
            \begin{table}[h]
            	\centering
            	\caption{Optimizer hyperparameters for matrix quadratic regression}
            	\begin{tabular}{cccc}
            		\toprule
            		Optimizer & $\gamma_0$ & $\beta$ or $(\beta_1, \beta_2)$ & inner steps  \\
            		\midrule
            		\PolarGrad (QDWH) & $4\times10^{-8}$ & N/A & $2$ \\
            		\PolarGrad (ZOLO-PD) & $3\times10^{-8}$ & N/A & N/A  \\
            		\PolarGrad (QDWH; lr $\downarrow$) & $4.75\times10^{-8}$ & N/A & $2$ \\
            		\Muon (NS) & $0.1$ & $0.95$ & $5$ \\
            		\Muon (QDWH) & $0.1$ & $0.95$ & $2$ \\
            		\Muon (ZOLO-PD) & $0.1$ & $0.95$ & N/A  \\
            		\Muon (QDWH; lr $\downarrow$) & $0.05$ & $0.95$ & $2$ \\
            		Newton & $0.25$ &  N/A & N/A \\
            		\Adam & $0.05$ & $(0.9, 0.999)$ & N/A \\
            		\Adam (lr $\downarrow$) & $0.05$ & $(0.9, 0.999)$ & N/A \\
            		\bottomrule
            	\end{tabular}
            \end{table}

\begin{figure}[h]
            	\centering
            	\includegraphics[width=\textwidth]{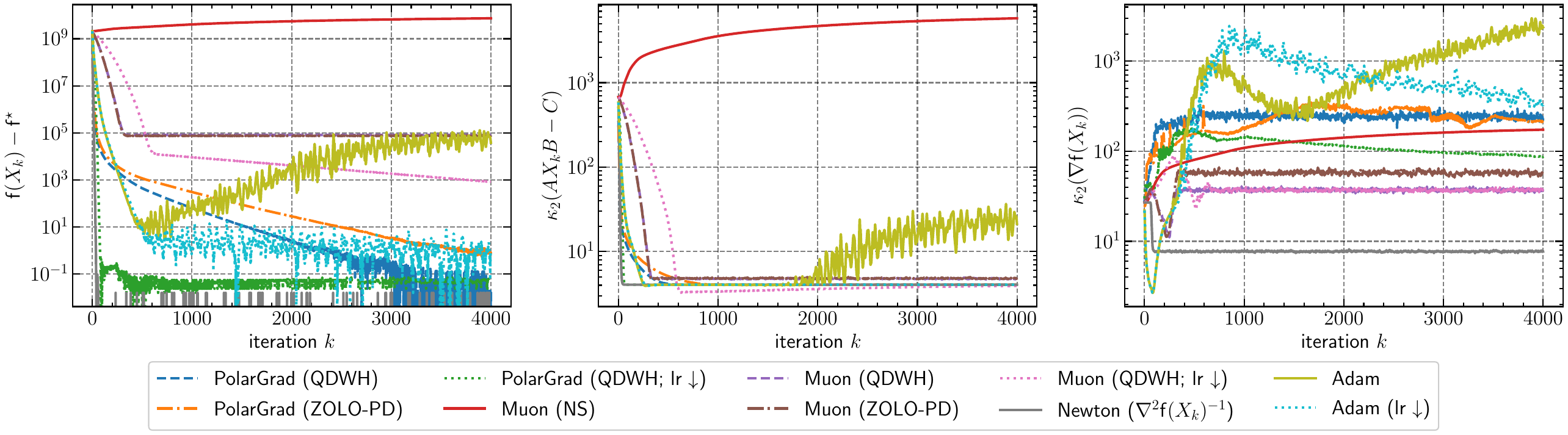}\\[2.5mm]
            	\includegraphics[width=\textwidth]{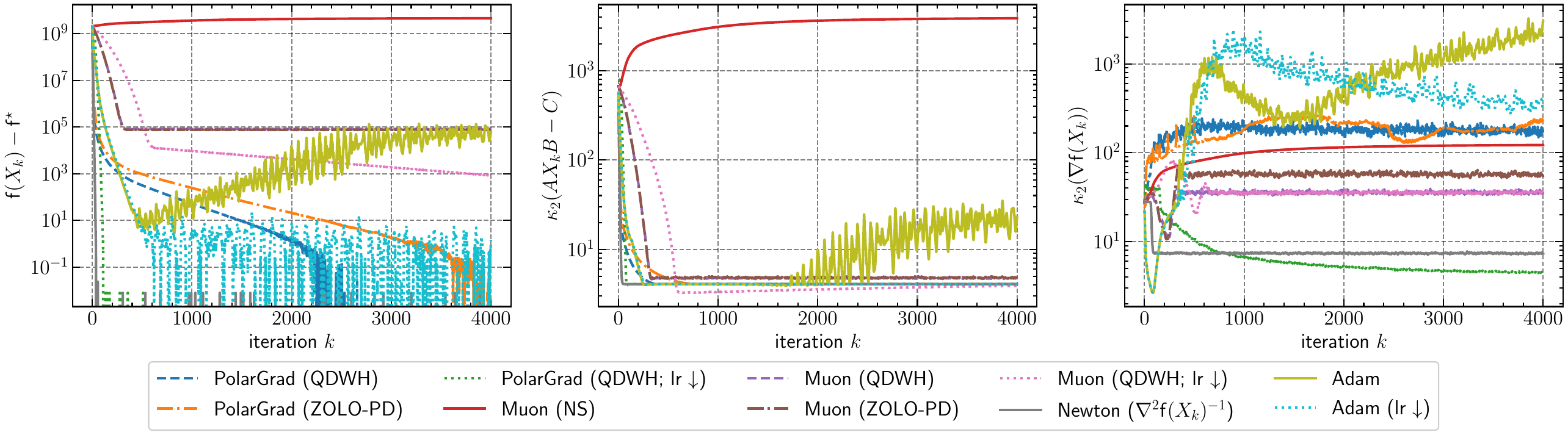}
            	\caption{Losses, residuals and gradient condition numbers of matrix quadratic regression (2nd and 3rd seeds). }
            	\label{fig:mat_quad_reg_2}
            \end{figure}

            We also give the simulation results of the remaining two random seeds in \Cref{fig:mat_quad_reg_2}. 

            \subsubsection{Momentum-First and Polar-First \PolarGradM}
            We are also interested in how the two types of EMA momentum are able to further accelerate convergence, possibly also with learning rate decay. We only use the QDWH algorithm for numerical polar decomposition here, and report the optimizer hyperparameters in \Cref{table:optim_hyperparams_mat_quad_reg_polargradm}.
            \begin{table}[h]
            	\centering
            	\caption{Optimizer hyperparameters for matrix quadratic regression for \PolarGradM}
                \label{table:optim_hyperparams_mat_quad_reg_polargradm}
            	\begin{tabular}{cccc}
            		\toprule
            		Optimizer & $\gamma_0$ & $\beta$ & inner steps  \\
            		\midrule
            		\PolarGradM (polar-first) & $4\times10^{-7}$ & $0.95$ & $2$ \\
            		\PolarGradM (polar-first; lr $\downarrow$) & $5\times10^{-7}$ & $0.95$ & $2$ \\
            		\PolarGradM (momentum-first) & $2\times10^{-7}$ & $0.9$ & $2$ \\
            		\PolarGradM (momentum-first; lr $\downarrow$) & $2.5\times10^{-7}$ & $0.9$ & $2$ \\
            		\bottomrule
            	\end{tabular}
            \end{table}
                    
\begin{figure}[h!]
            	\centering
            	\includegraphics[width=\textwidth]{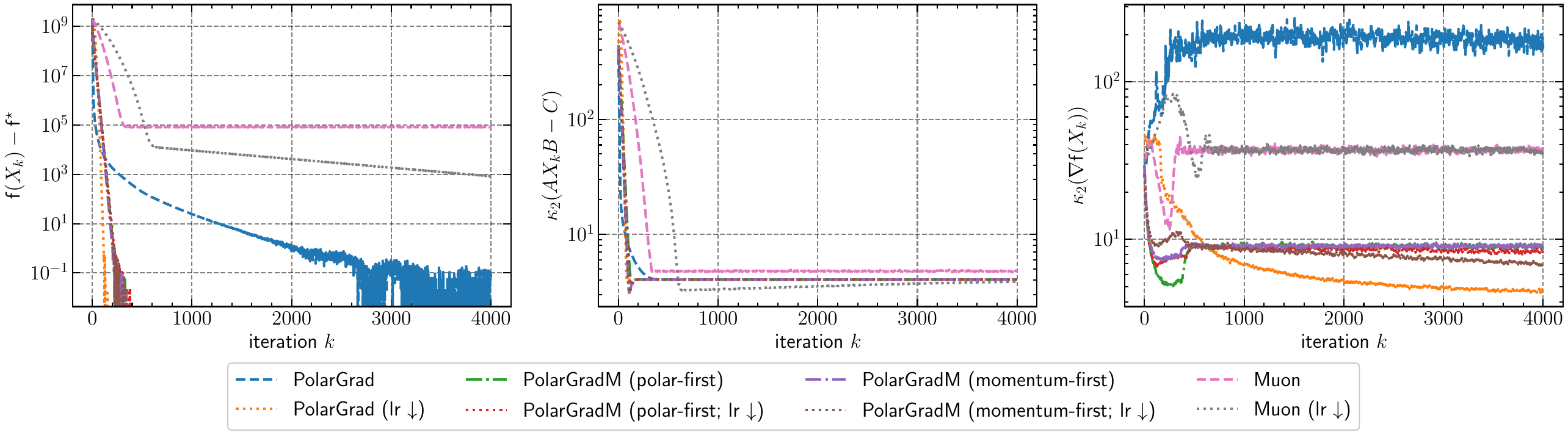}
            	\caption{Losses, residuals, and gradient condition numbers of matrix quadratic regression with momentum-first and polar-first \PolarGradM. }
            	\label{fig:mat_quad_reg_mom_1}
            \end{figure}

            We provide similar plots of losses and condition numbers in \Cref{fig:mat_quad_reg_mom_1}.         
            With either form of EMA momentum (with different values of momentum), we observe that \PolarGradM converges much faster than vanilla \PolarGrad, but slower than \PolarGrad with learning rate decay. Learning rate decay for \PolarGradM does not further accelerate convergence here. 
            
            \begin{figure}[h!]
            	\centering
            	\includegraphics[width=0.8\textwidth]{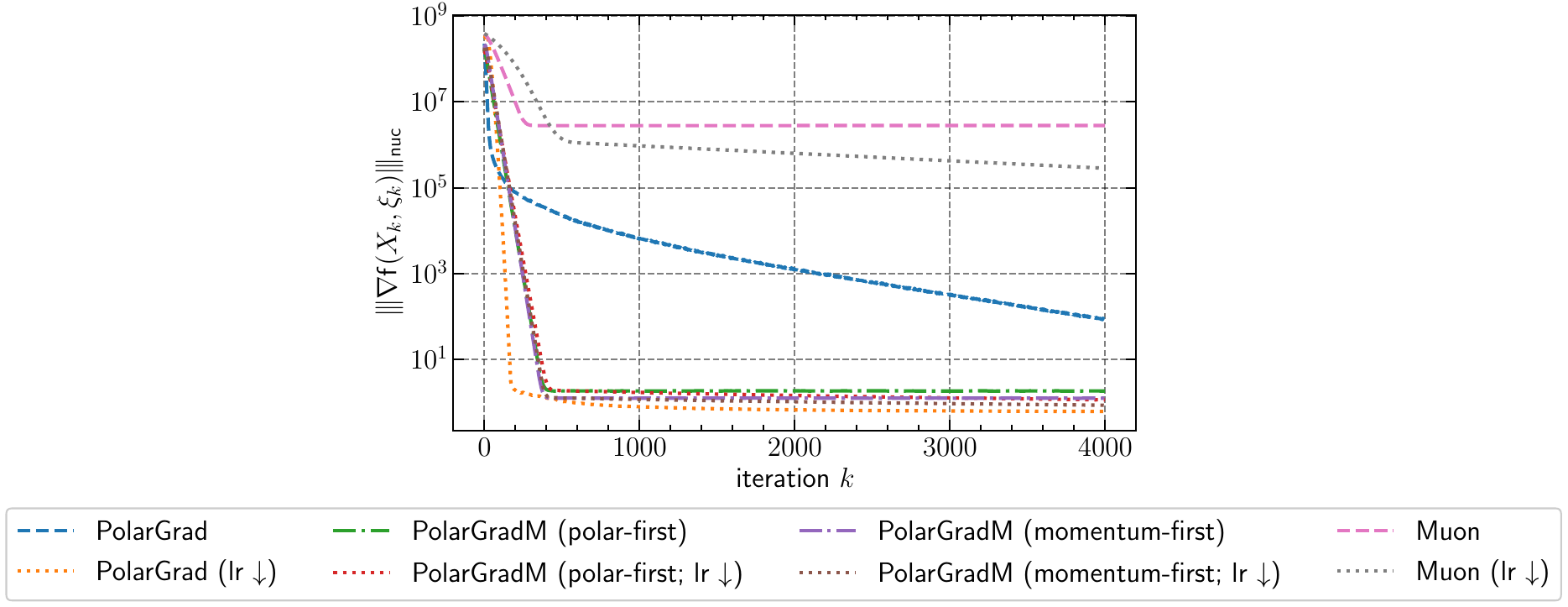}
            	\caption{Gradient nuclear norms of matrix quadratic regression with momentum-first and polar-first \PolarGradM.}
            	\label{fig:mat_quad_reg_mom_nuc_1}
            \end{figure}

            \subsection{Matrix Logistic Regression}
            The initialization $X_0$ has entries independently drawn from $\Unif(-1, 1)$. The matrices $A$ and $B$ have independent standard Gaussian entries. The matrix $C$ is generated by first independently drawing standard Gaussian entries, and setting each entry to be $1$ if it is greater than $0.5$ and $0$ otherwise. No weight decay is used in all optimizers. The learning rate decay schedule is a step scheduler which multiplies the base learning rate $\gamma_0$ by $0.95$ every $25$ steps. The optimizer hyperparameters are given in the table below. Default hyperparameters of \Adam in PyTorch are used.     
            \begin{table}[h!]
            	\centering
            	\caption{Optimizer hyperparameters for matrix logistic regression}
            	\begin{tabular}{cccc}
            		\toprule
            		Optimizer & $\gamma_0$ & $\beta$ or $(\beta_1, \beta_2)$ & inner steps  \\
            		\midrule
            		\PolarSGD (QDWH) & $2.5\times10^{-7}$ & N/A & $2$ \\
            		\PolarSGD (QDWH; lr $\downarrow$) & $5\times10^{-7}$ & N/A & $2$ \\
            		\Muon (NS) & $0.075$ & $0.95$ & $5$ \\
            		\Muon (QDWH) & $0.075$ & $0.95$ & $2$ \\
            		\Muon (QDWH; lr $\downarrow$) & $0.15$ & $0.95$ & $2$ \\
            		\Adam & $0.005$ & $(0.9, 0.999)$ & N/A \\
            		\Adam (lr $\downarrow$) & $0.01$ & $(0.9, 0.999)$ & N/A \\
            		\bottomrule
            	\end{tabular}
            \end{table}

\begin{figure}[h]
            	\centering
            	\includegraphics[width=\textwidth]{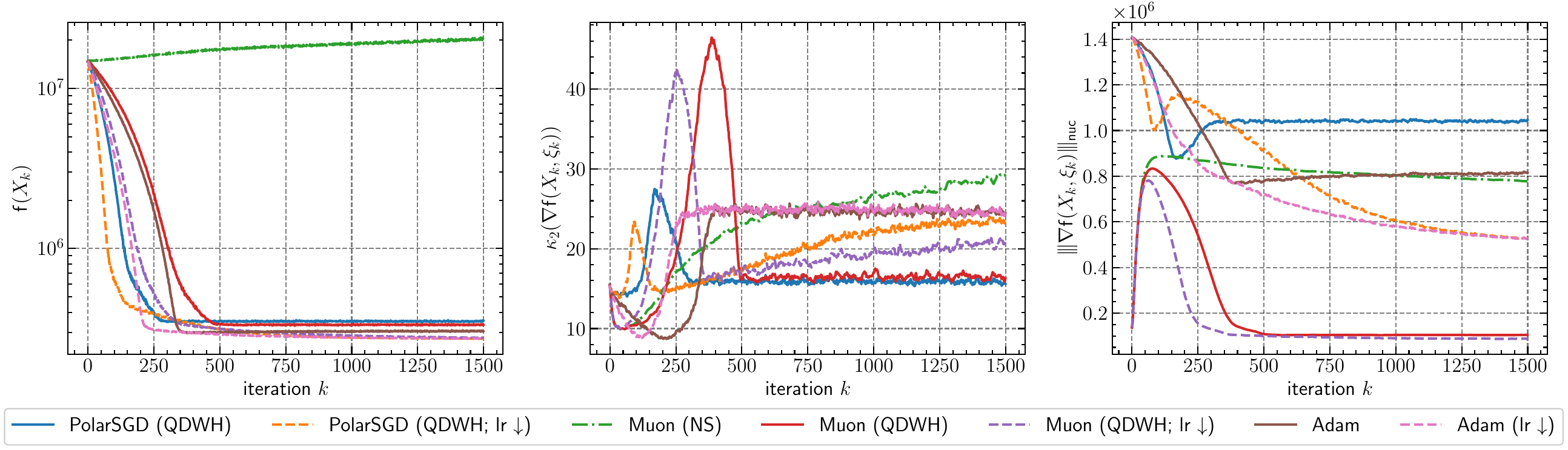}
            	\includegraphics[width=\textwidth]{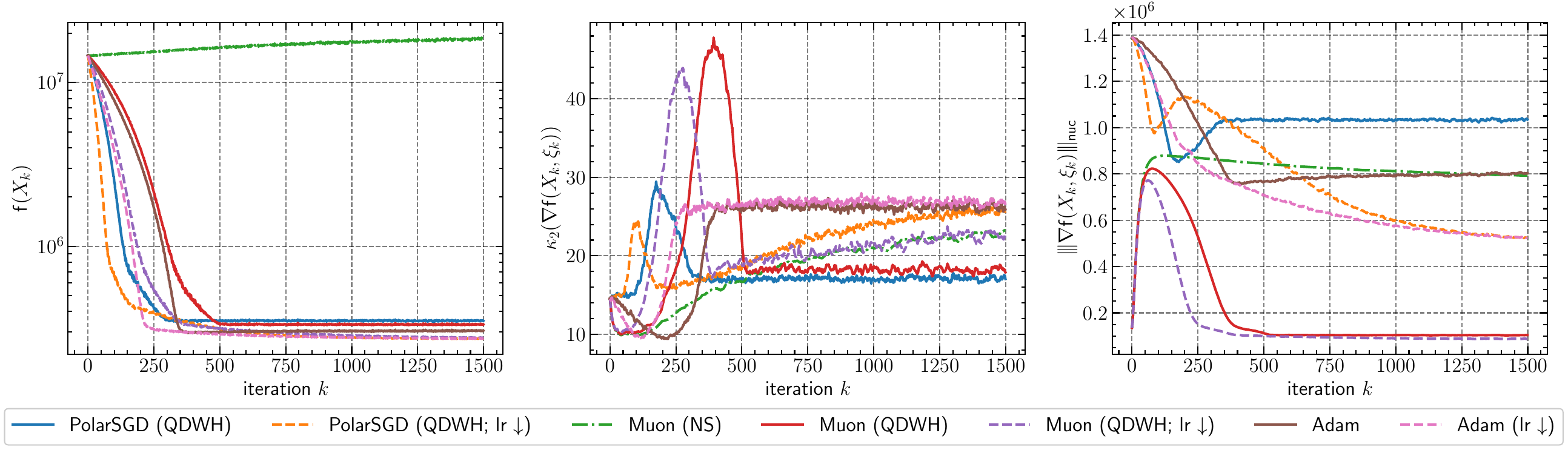}
            	\caption{Losses, gradient condition numbers and nuclear norms of matrix logistic regression (2nd and 3rd seeds).}
            	\label{fig:mat_logistic_reg_2}
            \end{figure}

            We also give the simulation results of the remaining two random seeds in \Cref{fig:mat_logistic_reg_2}. 

            \subsubsection{Momentum-First and Polar-First \PolarSGDM}
            We again study the two types of EMA momentum for this problem. The implementation of \PolarSGDM is based on the QDWH algorithm. 
            
            \begin{table}[h!]
            	\centering    	
            	\caption{Optimizer hyperparameters for matrix logistic regression for \PolarSGDM}
            	\begin{tabular}{cccc}
            		\toprule
            		Optimizer & $\gamma_0$ & $\beta$ & inner steps  \\
            		\midrule
            		\PolarSGDM (polar-first) & $5\times10^{-7}$ & $0.95$ & $2$ \\
            		\PolarSGDM (polar-first; lr $\downarrow$) & $5\times10^{-7}$ & $0.95$ & $2$ \\
            		\PolarSGDM (momentum-first) & $5\times10^{-7}$ & $0.9$ & $2$ \\
            		\PolarSGDM (momentum-first; lr $\downarrow$) & $5\times10^{-7}$ & $0.9$ & $2$ \\
            		\bottomrule
            	\end{tabular}
            \end{table}

            \begin{figure}[h!]
            	\centering
            	\includegraphics[width=\textwidth]{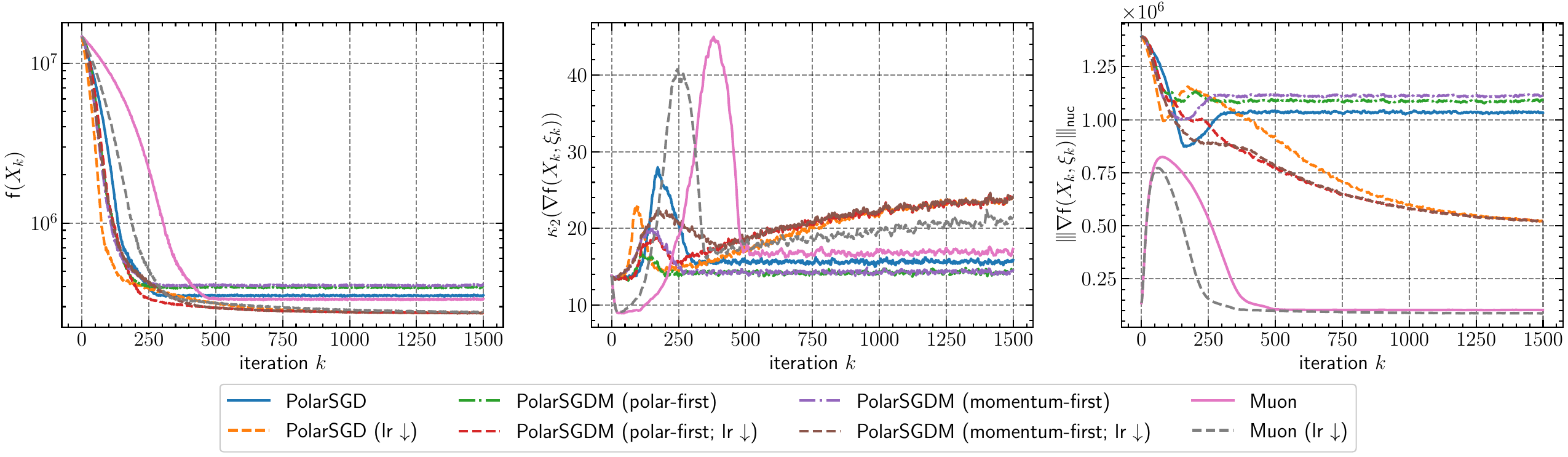}
            	\caption{Losses, gradient condition numbers and nuclear norms of matrix logistic regression with momentum-first and polar-first \PolarSGDM.}
            	\label{fig:mat_logistic_reg_mom}
            \end{figure}
            
            From \Cref{fig:mat_logistic_reg_mom}, there is not a significant distinction of the convergence behavior of these two types of momentum, both not being able to accelerate convergence much compared to vanilla \PolarSGD.

            \subsection{Low-Rank Matrix Completion}    
            The mask is formed by first generating entries from $\Unif(0,1)$ then setting each entry to be $1$ if it is smaller than $0.3$ and $0$ otherwise. 
            The ground truth low-rank matrix $M_\star\in\RR^{m\times n}$ is generated by $M_\star=UV^\top$, where $U\in\RR^{m\times r}$ and $V\in\RR^{n\times r}$ have independent standard Gaussian entries. 
            The initialization $(X_0, Y_0)$ has entries drawn independently from $\Unif(-1, 1)$. 
            No weight decay is used in all optimizers. The learning rate decay schedule is a step scheduler which multiplies the base learning rate $\gamma_0$ by $0.95$ every $25$ steps. The optimizer hyperparameters are given in the table below. Default hyperparameters of \Adam in PyTorch are used. 
            \begin{table}[h]
            	\centering
            	\caption{Optimizer hyperparameters for low-rank matrix completion}
            	\begin{tabular}{cccc}
            		\toprule
            		Optimizer & $\gamma_0$ & $\beta$ or $(\beta_1, \beta_2)$ & inner steps  \\
            		\midrule
            		\PolarGrad (QDWH) & $15$ & N/A & $2$ \\
            		\PolarGrad (QDWH; lr $\downarrow$) & $15$ & N/A & $2$ \\
            		\Muon (NS) & $0.25$ & $0.95$ & $5$ \\
            		\Muon (QDWH) & $0.25$ & $0.95$ & $2$ \\
            		\Muon (QDWH; lr $\downarrow$) & $0.25$ & $0.95$ & $2$ \\
            		\Adam & $0.05$ & $(0.9, 0.999)$ & N/A \\
            		\Adam (lr $\downarrow$) & $0.05$ & $(0.9, 0.999)$ & N/A \\
            		\AltGD & $50$ & N/A & N/A\\
            		\bottomrule
            	\end{tabular} 
            \end{table}
            
\begin{figure}[h!]
            	\centering
            	\includegraphics[width=\textwidth]{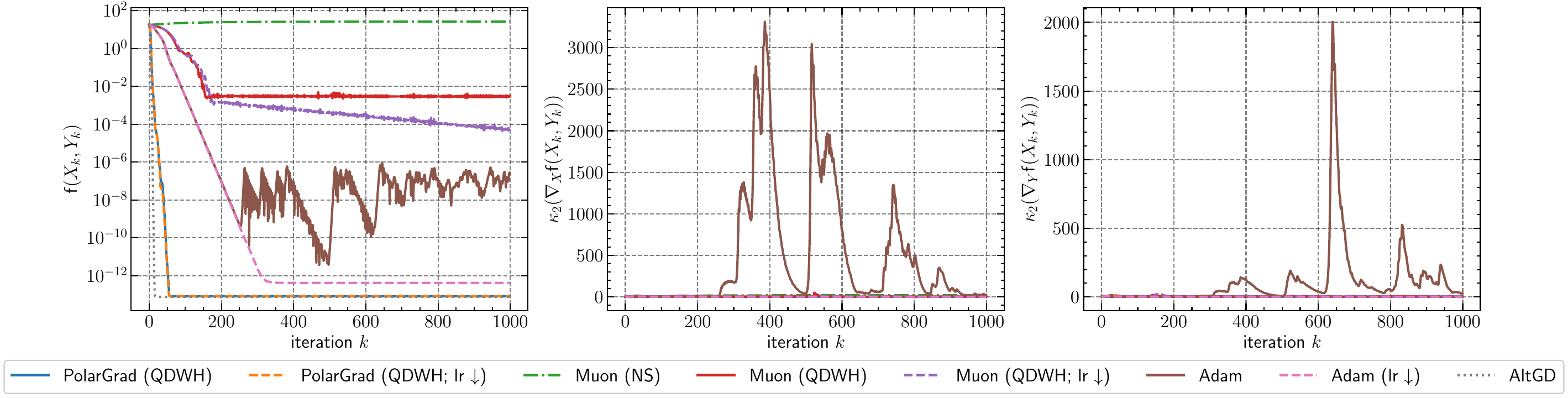}\\[2.5mm]
            	\includegraphics[width=\textwidth]{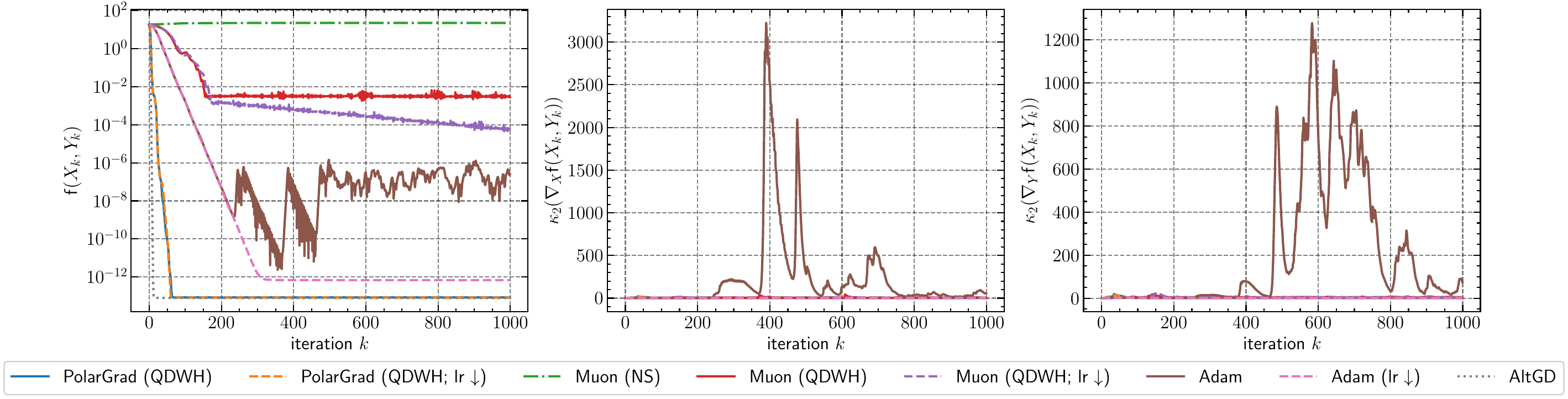}
            	\caption{Losses and gradient condition numbers of low-rank matrix completion (2nd and 3rd seeds).}
            	\label{fig:low_rank_mat_comp_2}
            \end{figure}

            We also give the simulation results of the remaining two random seeds in \Cref{fig:low_rank_mat_comp_2}. 
            Since the gradient condition numbers in \Cref{fig:low_rank_mat_comp_2} are dominated by \Adam, we also plots the figures again without \Adam (and \AltGD). 
            
            \begin{figure}[h!]
            	\centering
            	\includegraphics[width=\textwidth]{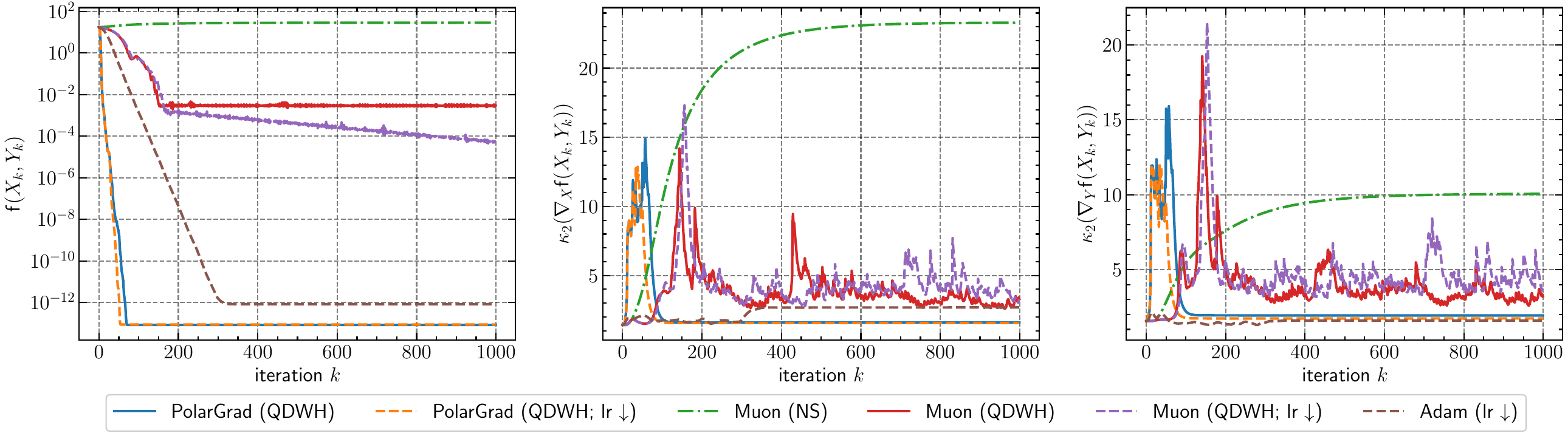}
            	\caption{Losses and gradient condition numbers of low-rank matrix completion.}
            	\label{fig:low_rank_mat_comp_1_2}
            \end{figure}
            
            \begin{figure}[h!]
            	\centering
            	\includegraphics[width=\textwidth]{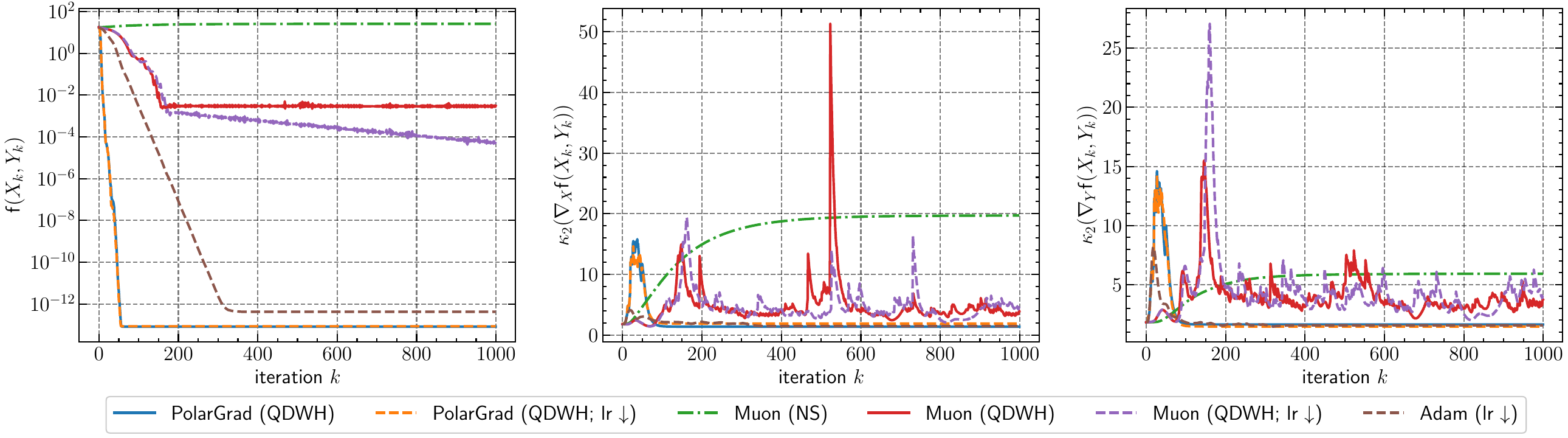}\\[2.5mm]
            	\includegraphics[width=\textwidth]{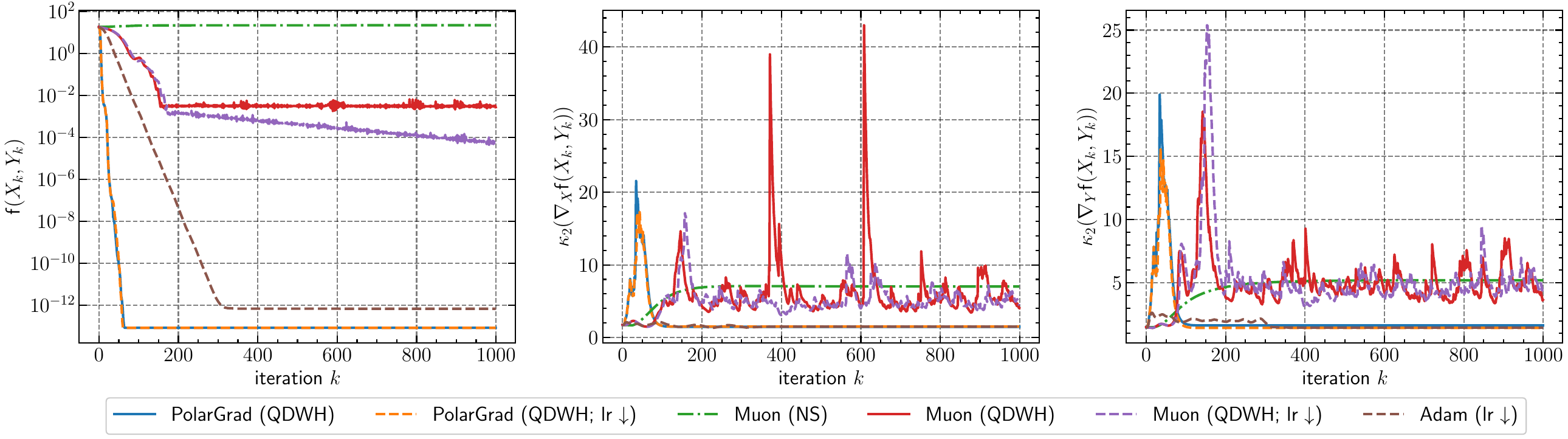}
            	\caption{Losses and gradient condition numbers of low-rank matrix completion (2nd and 3rd seeds).}
            	\label{fig:low_rank_mat_comp_2_2}
            \end{figure}

            We observe further from \Cref{fig:low_rank_mat_comp_1_2,fig:low_rank_mat_comp_2_2} that the gradient condition numbers of \Muon are highly fluctuating even at a later stage of training, whereas \PolarGrad is able to stabilize gradient condition numbers after achieving convergence, again indicating that the gradient condition number is indicative of the convergence of matrix gradient-based optimizers. 
            
\begin{figure}[h!]
            	\centering
            	\includegraphics[width=\textwidth]{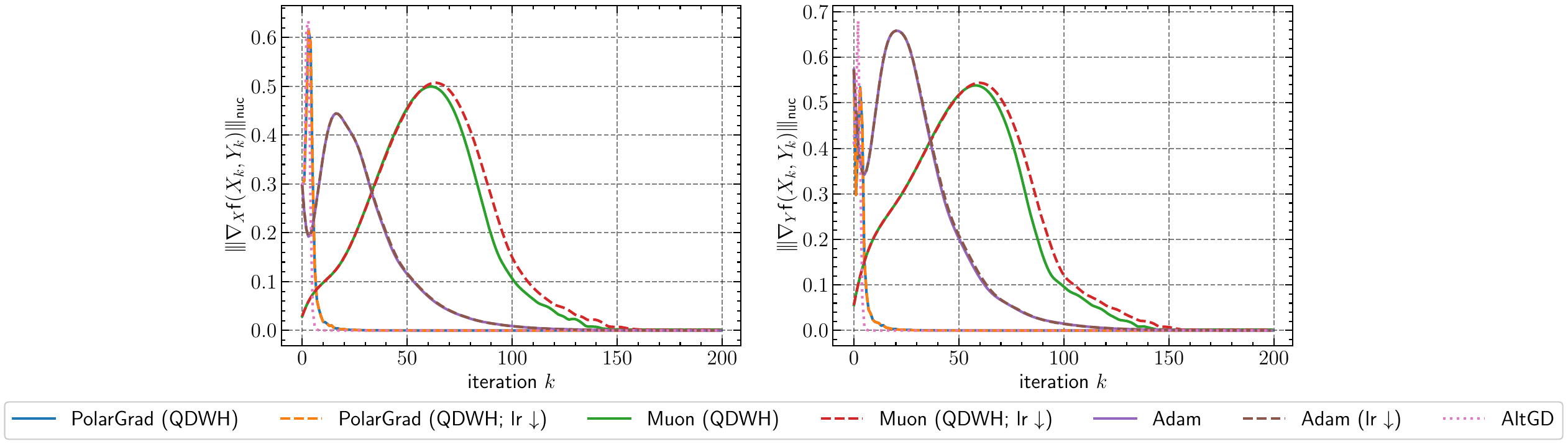}
            	\caption{Gradient nuclear norms of low-rank matrix completion.}
            	\label{fig:low_rank_mat_comp_1_3}
            \end{figure}

            Next, we also plot the gradient nuclear norms to evaluate the difference between \PolarGrad and \Muon in \Cref{fig:low_rank_mat_comp_1_3}. We observe that the gradient nuclear norms of \Muon converge to zero after roughly 150 iterations, but its objective values have not converged. \PolarGrad and \AltGD both converge within 20 iterations in terms of gradient nuclear norms. Again, the bell-shaped curves of the gradient nuclear norms for \Muon and \Adam has led to some potential relationship of a warmup-then-decay learning rate schedule, but we leave a more in-depth study on this for future work. 
            
            \subsubsection{Momentum-First and Polar-First \PolarGradM}
\begin{figure}[h!]
            	\centering
            	\includegraphics[width=\textwidth]{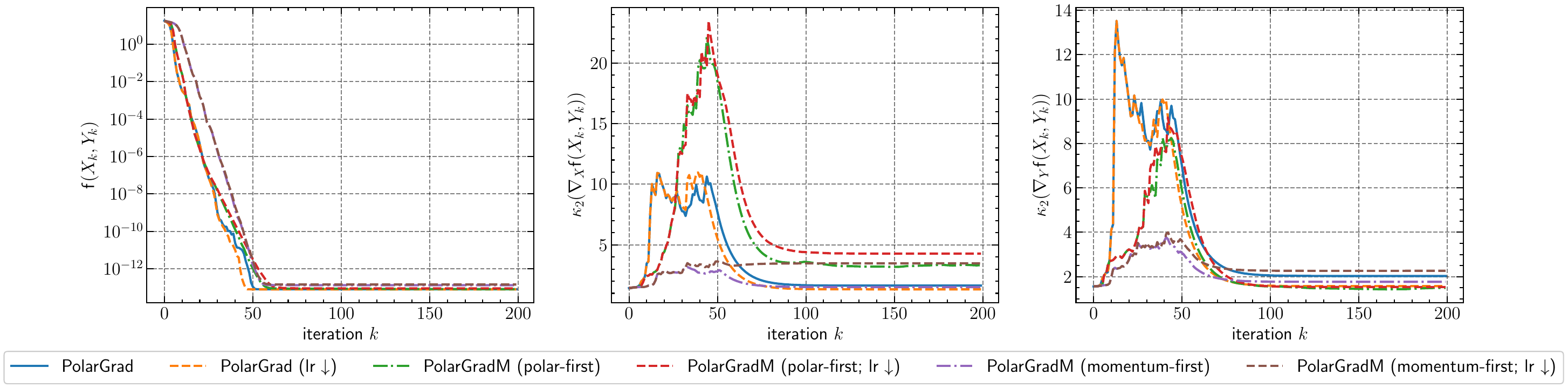}
            	\caption{Losses and gradient condition numbers of low-rank matrix completion with momentum-first and polar-first \PolarGradM.}
            	\label{fig:low_rank_mat_comp_mom_1}
            \end{figure}

\begin{figure}[h!]
            	\centering
            	\includegraphics[width=\textwidth]{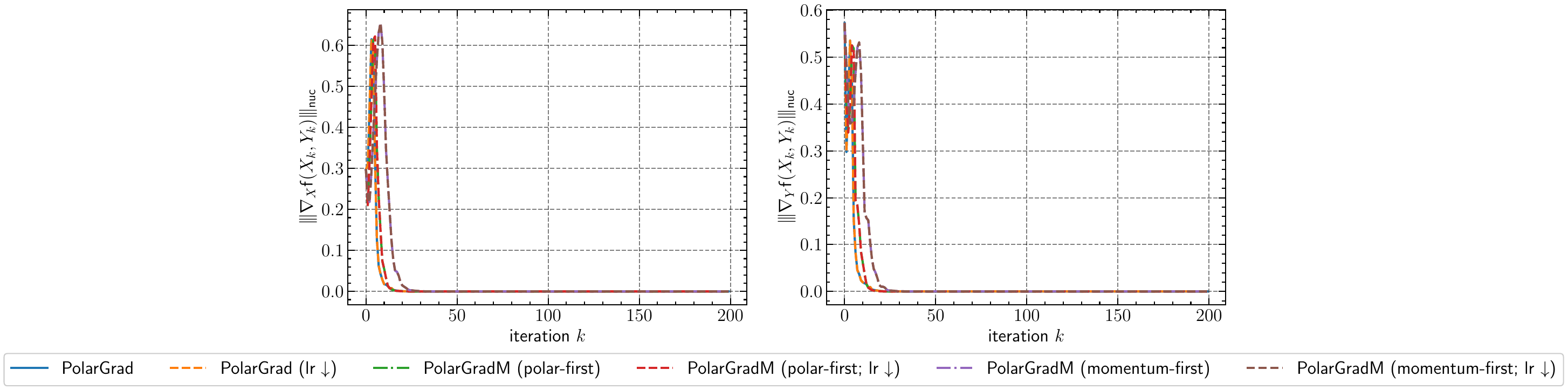}
            	\caption{Gradient nuclear norms of low-rank matrix completion with momentum-first and polar-first \PolarGradM.}
            	\label{fig:low_rank_mat_comp_mom_1_3}
            \end{figure}

            We compare the two types of EMA momentum and provide the hyperparameter setting of both momentum-first and polar-first \PolarGradM in the following table. Their plots are given in \Cref{fig:low_rank_mat_comp_mom_1,fig:low_rank_mat_comp_mom_1_3}. 
            \begin{table}[h!]
            	\centering
            	\caption{Optimizer hyperparameters for low-rank matrix completion}
            	\begin{tabular}{cccc}
            		\toprule
            		Optimizer & $\gamma_0$ & $\beta$ & inner steps  \\
            		\midrule
            		\PolarGradM (polar-first) & $15$ & $0.5$ & $2$ \\
            		\PolarGradM (polar-first; lr $\downarrow$) & $15$ & $0.5$ & $2$ \\
            		\PolarGradM (momentum-first) & $7.5$ & $0.5$ & $2$ \\
            		\PolarGradM (momentum; lr $\downarrow$) & $7.5$ & $0.5$ & $2$ \\
            		\bottomrule
            	\end{tabular}
            \end{table}

            We observe that we need to use a relatively small momentum in this nonconvex problem and are only able to recover comparable or even worse performance than vanilla \PolarGrad. Therefore, the use of momentum might not accelerate convergence in this problem. A thorough theoretical justification is left for future work. 
            
            \subsection{Qwen2.5 Pre-Training}
            The modified version of Qwen2.5 \citep{qwen2025qwen_full} is pre-trained on the OpenWebText-100k dataset\footnote{Available at \url{https://huggingface.co/datasets/Elriggs/openwebtext-100k}.} for one epoch, based on the toy example in the GitHub repository of \citep{liu2025muon}. 
            Qwen2.5 is chosen due to its more recent architecture, incorporated with many architectural design advances. 
            It only has 12 hidden layers and 16 heads, but without tie embedding (i.e., the embedding and classification head weight matrices are separate parameters) as we want to train both the embedding and head layers with \PolarSGDM. Its tokenizer has a vocabulary size of 151,936 (about three times that of GPT-2). This rather large vocabulary size indeed poses challenges to model training and leads to potential training instability. The implementation of \PolarSGDM is based on the QDWH algorithm. The model specifications (including those of GPT-2 Small and Medium in \Cref{sec:add_expt}), training hyperparameters and optimizer hyperparameters are provided in the following tables. Weight decay is not used for \Muon and \PolarSGDM.

            \begin{table}[h]    
            	\centering
            	\caption{Specifications of language models}
            	\label{table:settings}
            	\begin{tabular}{rrrr}
            		\toprule
            		Model & Qwen2.5 & GPT-2 Small 124M & GPT-2 Medium 350M \\
            		\midrule
            		$n_{\mathrm{params}}$ & 540,865,536 &  275,742,772 & 454,496,336 \\
            		$d_{\mathrm{model}}$ & 1024 & 768 & 1024 \\
            		$n_{\mathrm{layers}}$ & 12 & 12 & 6 \\
            		$n_{\mathrm{heads}}$ & 16 & 6 & 8 \\
            		$d_{\mathrm{head}}$ & 64 & 128 & 128 \\
            		vocab size & 151936 & 50304 & 50257 \\
            		layer norm & RMSNorm & RMSNorm & RMSNorm \\
            		activation & SiLU & ReLU$^2$ & ReLU$^2$ \\
            		\bottomrule
            	\end{tabular}
            \end{table}
            
            \begin{table}[h]
            	\centering
            	\caption{Training hyperparameters for Qwen2.5 pre-training}
            	\label{table:hyperparams_qwen}
            	\begin{tabular}{lc}
            		\toprule
            		Model & Qwen2.5 on OpenWebText-100k \\
            		\midrule
            		Training steps & 13299 \\
            		Sequence length & 512 tokens \\
            		Learning rate decay ratio (training steps) & $40\%$ \\
            		Batch size & 16 sequences \\
            		Precision & \texttt{bfloat16} \\         
            		Data-parallel size & 1 \\   
            		\bottomrule
            	\end{tabular}
            \end{table}
            
            The learning rate schedule for \AdamW is linear warmup (100 steps) and cosine decay to $0$, while the learning rate schedule for the other two optimizer combinations is linear decay from $\gamma_0$ to $0$ for the last $40\%$ of training steps. We use a weight decay of $0.1$ for \AdamW and no weight decay for \Muon and \PolarSGDM. 
            \begin{table}[H]
            	\centering
            	\caption{Optimizer hyperparameters for Qwen2.5 pre-training}
            	\label{table:optim_hyperparams_qwen}
            	\small
            	\begin{tabular}{cccccc}
            		\toprule
            		Optimizer & $\gamma_0$ & $\beta_{\Muon}$ & $\beta_{\PolarSGDM}$ & $(\beta_1, \beta_2)$  & inner steps \\
            		\midrule
            		\AdamW & $0.001$ & N/A & N/A & $(0.9, 0.95)$ & N/A \\
            		\Muon $+$ \AdamW & $(0.001, 0.001)$ & $0.95$ & N/A & $(0.9, 0.95)$ & $5$ (\Muon) \\
            		\Muon $+$ \PolarSGDM & $(0.001, 0.001)$ & $0.95$ & $0.5$ & N/A & $5$ (\Muon and QDWH) \\
            		\bottomrule
            	\end{tabular}        
            \end{table}
            
            It turns out that \PolarSGDM works better with a small momentum, probably due to the inclusion of the nuclear norm scaling term. 
            
            We also plot the gradient nuclear norms of the embedding and the head weight matrices, which can be viewed as indicators of convergence. 
            \begin{figure}[H]
            	\centering
            	\includegraphics[width=0.7\textwidth]{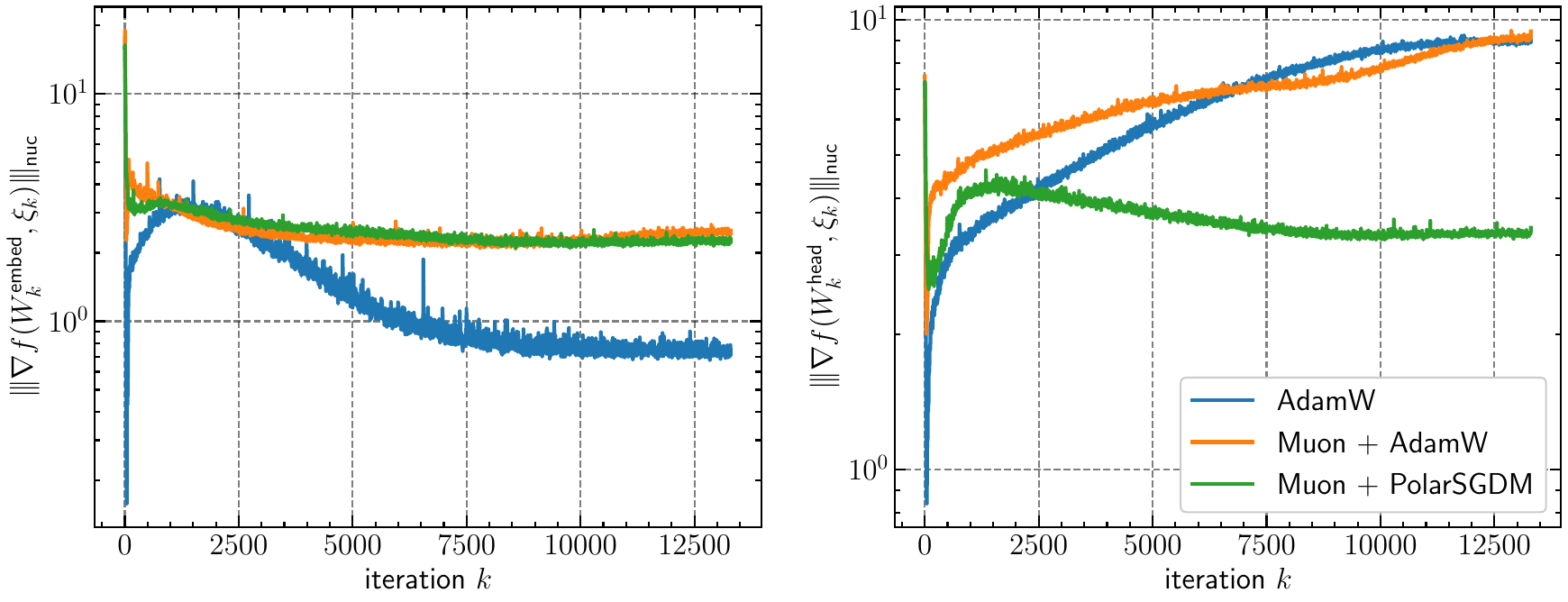}
            	\caption{Gradient nuclear norms of Qwen2.5 pre-training.}
            	\label{fig:qwen_2}
            \end{figure}
            We observe that the gradient nuclear norm of the head weight matrix is actually growing without converging when trained with \AdamW (blue and orange lines), indicating that \AdamW might not be appropriate for training such layers
            
            \subsubsection{Momentum-First and Polar-First \PolarSGDM}
            We now compare the two possible types of EMA momentum, momentum-first (which is similar to \Muon) and polar-first. The optimizer hyperparameters are the same as those in \Cref{table:optim_hyperparams_qwen}. 
            
            \begin{figure}[H]
            	\centering
            	\includegraphics[width=\textwidth]{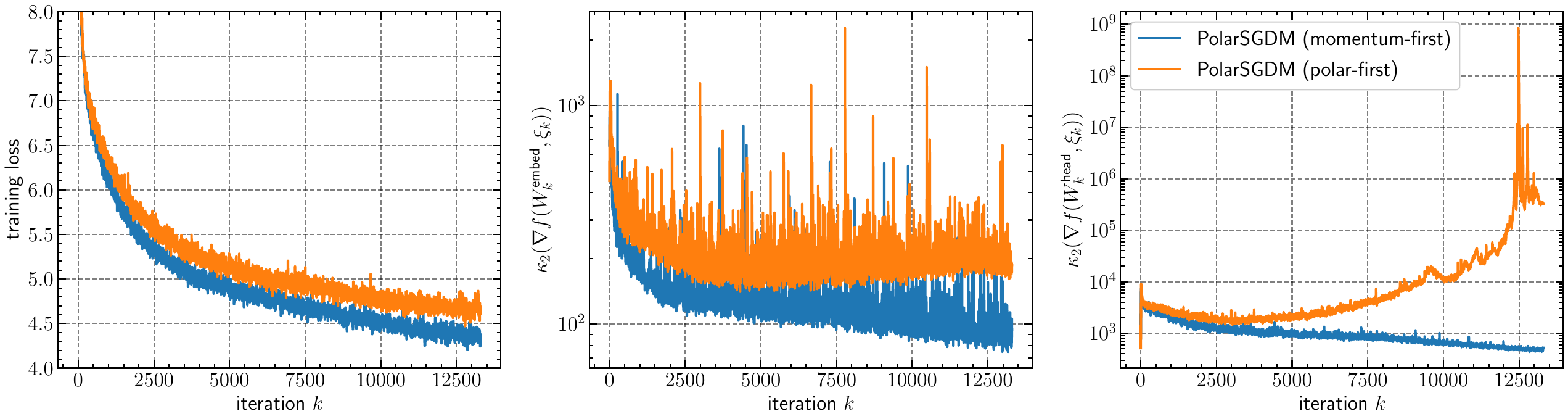}
            	\caption{Training losses and gradient condition numbers of Qwen2.5 pre-training with momentum-first and polar-first \PolarSGDM. }
            	\label{fig:qwen_3}
            \end{figure}
            
            \begin{figure}[H]
            	\centering
            	\includegraphics[width=0.7\textwidth]{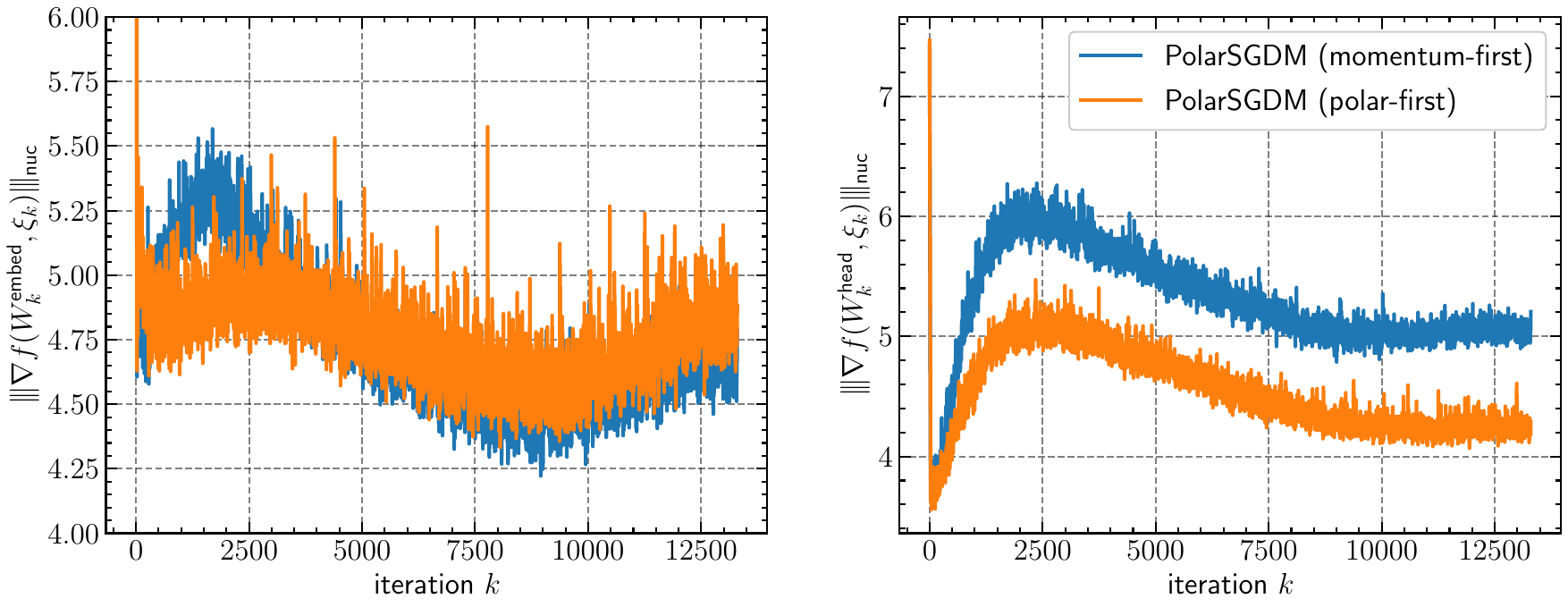}
            	\caption{Gradient nuclear norms of Qwen2.5 pre-training with momentum-first and polar-first \PolarSGDM.}
            	\label{fig:qwen_4}
            \end{figure}
            
            We see that the polar-first momentum is less desirable in terms of training loss convergence, and the gradient condition number of the head weight matrix also grows throughout training and has a strong spike at the end of training, although we do not tune its momentum parameter thoroughly in this experiment. This might indicate that the momentum-first momentum is more preferred for \PolarSGDM similar to \Muon, but we need more ablation studies to draw a definite conclusion here.

            \subsection{GPT-2 Small 124M Pre-Training}   
            We give the training and optimizer hyperparameters in \Cref{table:hyperparams_gpt2,table:hyperparams_opt_gpt2}.

            \begin{table}[h!]
            	\centering
            	\caption{Training hyperparameters for GPT-2 Small 124M pre-training}
            	\label{table:hyperparams_gpt2}
            	\begin{tabular}{lc}
            		\toprule
            		Model & GPT-2 Small 124M on FineWeb \\
            		\midrule
            		Training steps & 5000 \\
            		Sequence length & 1024 tokens \\
            		Learning rate schedule & linear decay from $\gamma_0$ to $0$ \\
            		Learning rate decay ratio (training steps) & $40\%$ \\
            		Global batch size & 1024 \\
            		Local batch size & 128 \\
            		Precision & \texttt{float32} for embedding; \texttt{bfloat16} otherwise \\         
            		Data-parallel size & 8 \\   
            		\bottomrule
            	\end{tabular}
            \end{table}
            
            \begin{table}[h!]
            	\centering    	
            	\caption{Optimizer hyperparameters for GPT-2 Small 124M pre-training}
            	\label{table:hyperparams_opt_gpt2}
            	\begin{tabular}{ccc}
            		\toprule
            		Hyperparameters & \Muon $+$ \Adam & \Muon $+$ \PolarSGDM \\
            		\midrule
            		$\gamma_0^\text{scalar}$ & $0.04$ & $0.04$ \\
            		$\gamma_0^\text{hidden}$ & $0.05$ & $0.05$ \\
            		$\gamma_0^\text{embed}$ & $0.6$ & $5$ \\
            		$\gamma_0^\text{value\_embed}$ & $0.6$ & $50000$ \\
            		$\gamma_0^\text{head}$ & $0.008$ & $0.02$ \\
            		$\beta_{\Muon}$ & $0.95$ & $0.95$ \\
            		$\beta_{\PolarSGDM}$ & N/A & $0.5$ \\
            		$(\beta_1, \beta_2)$ & $(0.8, 0.95)$ & N/A \\
            		$\varepsilon$ & $10^{-10}$ & N/A \\
            		inner steps & $5$ & $5$ (\Muon); $5$ (QDWH) \\
            		\bottomrule
            	\end{tabular}     
            \end{table}

            \section{Additional Numerical Experiments}
            \label{sec:add_expt}
            We also provide additional numerical experiments on GPT-2 Medium pre-training in the section.

            \subsection{GPT-2 Medium 350M Pre-Training}
            We now move on to the GPT-2 Medium track of the Modded-NanoGPT repository on the FineWeb dataset, making use of the setting of the 04/22/25 record. 
            We also keep the same optimizer choices as GPT-2 Small. We give the training and optimizer hyperparameters in \Cref{table:hyperparams_gpt2m,table:hyperparams_opt_gpt2m}. 
            
            \begin{table}[h]
            	\centering
            	\caption{Training hyperparameters for GPT-2 Medium 350M pre-training}
            	\label{table:hyperparams_gpt2m}
            	\begin{tabular}{lc}
            		\toprule
            		Model & GPT-2 Medium 350M on FineWeb \\
            		\midrule
            		Training steps & 5960 \\
            		Sequence length & 1024 tokens \\
            		Learning rate schedule & linear decay from $\gamma_0$ to $0$ \\
            		Learning rate decay ratio (training steps) & $70\%$ \\
            		Global batch size & 512 \\
            		Local batch size & 64 \\
            		Precision & \texttt{bfloat16} \\         
            		Data-parallel size & 8 \\   
            		\bottomrule
            	\end{tabular}
            \end{table}

            \begin{table}[h!]
            	\centering    	
            	\caption{Optimizer hyperparameters for GPT-2 Medium 350M pre-training}
            	\label{table:hyperparams_opt_gpt2m}
            	\begin{tabular}{ccc}
            		\toprule
            		Hyperparameters & \Muon $+$ \Adam & \Muon $+$ \PolarSGDM \\
            		\midrule
            		$\gamma_0^\text{scalar}$ & $0.015$ & $0.015$ \\
            		$\gamma_0^\text{hidden}$ & $0.025$ & $0.025$ \\
            		$\gamma_0^\text{embed}$ & $0.3$ & $2.5$ \\
            		$\gamma_0^\text{value\_embed}$ & $0.3$ & $25000$ \\
            		$\gamma_0^\text{head}$ & $1/320$ & $0.015$ \\
            		$\beta_{\Muon}$ & $0.95$ & $0.95$ \\
            		$\beta_{\PolarSGDM}$ & N/A & $0.5$ \\
            		$(\beta_1, \beta_2)$ & $(0.8, 0.95)$ & N/A \\
            		$\varepsilon$ & $10^{-10}$ & N/A \\
            		inner steps & $5$ & $5$ (\Muon); $2$ (QDWH) \\
            		\bottomrule
            	\end{tabular}     
            \end{table}
            
\begin{figure}[h!]
            	\centering
            	\includegraphics[width=\textwidth]{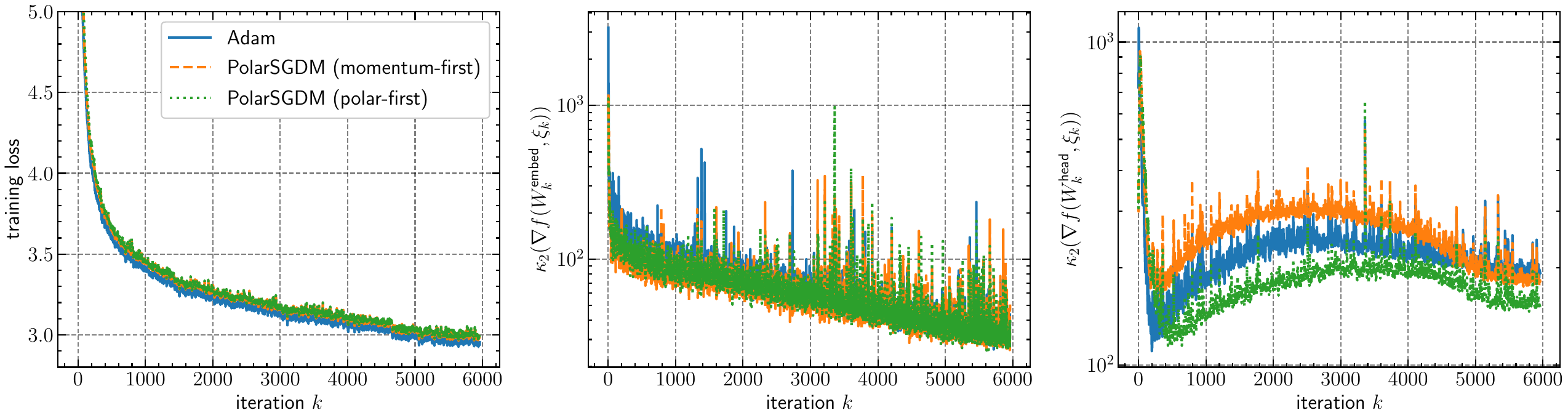}
            	\caption{Training losses and gradient condition numbers of GPT-2 Medium 350M pre-training.}
            	\label{fig:gpt2_medium}
            \end{figure}

\begin{figure}[h!]
            	\centering
            	\includegraphics[width=\textwidth]{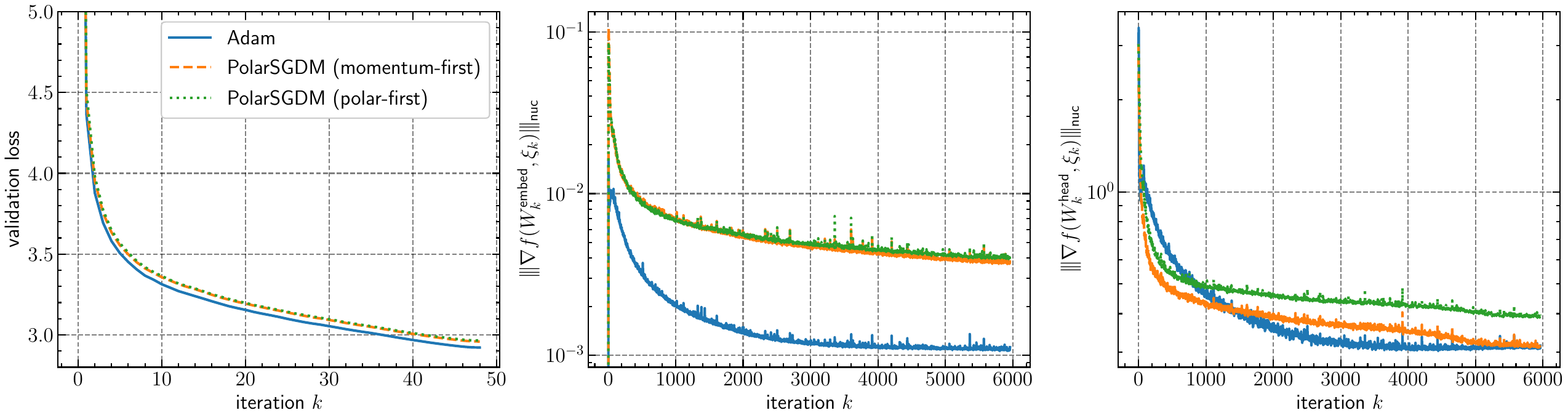}
            	\caption{Validation losses and gradient nuclear norms of GPT-2 Medium 350M pre-training.}
            	\label{fig:gpt2_medium_2}
            \end{figure}

            From \Cref{fig:gpt2_medium,fig:gpt2_medium_2}, we are able to make similar takeaways as the GPT-2 Small experiments. 

    \end{document}